\documentclass[12pt]{article}
\title{The Calculus and Gauge Integrals, \newline
by Ralph Henstock}
\author{Edited by P.~Muldowney}
\date{}

\newtheorem{theorem}{Theorem}

\newtheorem{example}{Example}

\newcommand{\da}{\delta}

\newcommand{\vt}{\vspace{5pt}\\}
\newcommand{\R}{\mathbf{R}}

\newcommand{\ve}{\varepsilon}

\newcommand{\proof}{\noindent\textbf{Proof. }\noindent}
\newcommand{\nproof}{\hfill$\mathbf{\bigcirc}\vspace{12pt}$ }

\date{}
\begin{document}
\maketitle

\section{Introduction}\label{Introduction}
\textsf{In a letter to Cambridge University Press, dated 18 October 1993, Ralph Henstock (1923--2007) proposed \emph{``a book that begins with the calculus integral and ends with the properties of the gauge (Riemann-complete, generalized Riemann, or Kurzweil-Henstock) integral, without having to define measure nor Lebesgue integration, and written for a beginning analyst, or a scientist interested in mathematics. This project has interested me for years, and has been revit\-alised by the enclosed review (of my fourth book)}\footnote{A review of \emph{The General Theory of Integration, R.~Henstock (1991)}, by Robert Bartle, appeared in the Bulletin of the American Mathematical Society, Volume 29, Number 1, July 1993, pages 136--139;  implying, in the words of Henstock,  ``the need of such a more elementary book.'' }.''}

\textsf{A reply (20 October 1993) by Roger Astley recommended an elementary treatment of the subject \emph{``whilst maintaining an appropriate level of rigour.''}
In the event, material  for about four chapters---Sections \ref{s1} to \ref{s4} of this document---was produced, but remained unpublished in Henstock Archive (2007)}. 

\textsf{The manuscript starts with Section \ref{s1} below which contains a review of basic operations of differentiation and integration. Section \ref{s2} introduces the properties of real numbers by means of infinite decimals.\footnote{This may benefit a beginner who  struggles with Dedekind sections, but understands that $\sqrt{2}$ is not rational and thinks of irrational real numbers as  non-terminating decimals.} This is followed by a thorough examination of the key points of traditional Riemann integration, framed in such a way as to lead naturally in Section \ref{s3} to an exposition of Riemann-complete integration---that is, the basic gauge integral of Henstock and Kurzweil.}

In ``\emph{Beyond dominated convergence: newer methods of integration}'', \textsf{Muldowney (2016) argues that the post-Riemann development of integration theory was driven, {not} specifically by problems of integrability of strange and exotic functions, but by the quest for rules and conditions  for taking limits under the integral sign, such as integrating series of functions; a quest which appeared to culminate in Lebesgue's dominated convergence theorem.}

\textsf{This aspect of integration is what motivates this previously unpublished manuscript of Henstock.  Theorems 3.1.3, 3.1.4 and 3.1.5 of Henstock (1991) give necessary and sufficient conditions for taking limits under the integral sign.\footnote{See also Section 4.13 (pp.~174--178) of Muldowney (2012), and Muldowney (2016).} }

\textsf{Section \ref{s4} of this manuscript presents these convergence criteria in a lucid, clear and careful manner, and then sets the traditional monotone and dominated convergence theorems of Lebesgue in this new context. This section alone warrants making this manuscript more widely available. Tonelli's theorem is presented in these terms, and there are indications (see (\ref{eq 3.1.4}), for example) that Henstock had in mind some parallel development of Fubini's theorem.}

\textsf{Readers who are familiar with basic integration theory can go direct to the convergence theorems in Section \ref{s4}. Example \ref{ex 4.1.1} is of particular interest.}

$\;\;\;\;\;\;\;\;\;\;\;\;\;\;\;\;\;\;\;\;\;\;\;\;\;\;\;\;\;\;\;\;\;\;\;\;\;\;\;\;\;\;\;\;\;\;\;\;\;\;\;\;\;\;\;\;\;\;\;\;\;\;\;\;\;\;\;
$
\textsf{--  \emph{P.~Muldowney, July 5 2016}}


\section{The Calculus}\label{s1}
\subsection{The Rules for Differentiation and Integration}\label{s1.1}
In this first chapter I have deliberately avoided the notation and exactness of modern analysis, in order to look at what formulae the calculus gives, before defining real numbers and a simple integration process in Chapter 2. Here we consider integration as an inverse operation to differentiation and assume the existence of what we need when necessary.

Differentiation and integration were first systematized by I.~Newton  (1642 --1727) and G.W.~Leibnitz (1646--1717). Over the years the notations have changed to become what we now use. In the simplest form of the first operation, given numbers $a<b$ and a function $F(x)$ defined in the range $a\leq x\leq b$, written $[a,b]$ and called a \emph{closed interval}, we could draw the graph of $y=F(x)$ there. We use
\begin{equation}\label{eq 1.1.1}
\frac{F(x+h)-F(x)}{h},
\end{equation}
the gradient of the chord of the curve from $x$ to $x+h$, naturally for $h \neq 0$, or else we would have the ridiculous fraction $0/0$. Taking the limit as $h \rightarrow 0$, $h$ positive and negative, $x$ fixed, and $x$, $x+h$ in the range, we obtain the gradient of the tangent to the curve at $x$, a value $f(x)$ depending on $x$, and so a function called the \emph{derivative} of $F$, and written $dy/dx$ or $dF/dx$ or $F'(x)$. Note that at $a$ and $b$ the derivative is one-sided.

\begin{example}
\label{ex 1.1.1}
If $y=F(t)$ represents the distance at time $t$ of a car driving along a road from some starting point, then, replacing $x$ by $t$, $f(t)$ represents its instantaneous speed. If $y=F(t)$ is its speed, $f(t)$ is its acceleration. Assuming that the car has constant mass, the resultant of forces acting on the car is equal to the mass times the acceleration, one of Newton's laws. Considering rockets instead of cars, the mass decreases as the propellant burns, and the force is the derivative of the mass times the speed in rectilinear flight. Such examples explain why Newton had to develop the calculus. 
\end{example}
In using the calculus, various rules were given which we can now write as follows.
\begin{equation}
\label{eq 1.1.2}
\mbox{If } F \mbox{ is constant throughout }[a,b], \mbox{ then } f=0 \mbox{ there.}
\end{equation}
For $F(x+h)=F(x)$ there, and the fraction (\ref{eq 1.1.1}) is $0$.

Next, if the derivatives at $x$ of $F$, $G$ are $f(x)$, $g(x)$, respectively, and if $m$, $n$ are constants, then {the derivative of } $mF+nG$ { is } 
\begin{equation}
\label{1.1.3}
mf(x) + ng(x) 
\end{equation}
{ there, because}
\begin{eqnarray*}
&&\frac{\left\{mF(x+h) + nG(x+h)\right\} - \left\{mF(x) + nG(x)\right\}}{h}\vt 
&  &\;\;=\;\;
m \left(\frac{F(x+h)-F(x)}h\right) +n\left( \frac{G(x+h)-G(x)}h\right) \vt
&&\;\;\rightarrow \;\; mf(x) + ng(x).
\end{eqnarray*}
{In (\ref{1.1.3}), the derivative of } $FG$ { is }
\begin{equation}
\label{1.1.4}
 fG+Fg.
\end{equation}
$F(x+h)G(x+h) -F(x)G(x)$ can be written as
\[
\begin{array}{rll}
&&\{F(x+h) - F(x)\}G(x) + F(x) \{G(x+h)-G(x)\}\vt & +& \{F(x+h)- F(x)\}\{G(x+h) -g(x)\}.
\end{array}
\]
Multiplying by $1/h = (1/h)(1/h)h$, the expression tends to
\[
f(x)G(x) +F(x)g(x) + f(x)g(x)\times 0.
\]
If, over an $x$-interval $[u,v]$, ($u<v$), the function $K$ is defined with values in $[a,b]$, and if $dK/dx$ is the derivative of $K$, then
{the derivative of } $F(K(x))$ { is }
\begin{equation}
\label{1.1.5}
\frac{dF(K)}{dK} \cdot \frac{dK}{dx},
\end{equation}
normally. This fails where $K(x+h) - K(x)=0$ for certain $x$, $h$. But as $K(x+h) =K(x)$, the fraction on the left is $0$. The formula still holds since the derivative of $K$ is given to exist, and so is $0$ if such $h$ are arbitrarily small for a particular $x$.

The rule for differentiating fractions can, with care, be deduced from these. First, to differentiate $1/G$ we have to assume $G(x) \neq 0$, and $G(x+h) \neq 0$ for $h$ near $0$. Then

\begin{eqnarray*}
G(x+h) - G(x) &=& \left( \frac{G(x+h) -G(x)}h\right)h \vt
&\rightarrow & g(x) \cdot 0, \vt
G(x+h) & \rightarrow & G(x), \vt
\frac{ \frac 1{G(x+h)} - \frac 1{G(x)} }h &=& 
\frac{G(x) -G(x+h)}{hG(x)G(x+h)} \vt
&\rightarrow & \frac{-g(x)}{G^2(x)}.
\end{eqnarray*}
This with (\ref{1.1.4}) gives the result: If, in \ref{1.1.3}, $[a,b]$ is an interval over which $G \neq 0$, then
{the derivative of } $\frac FG$ { is }
\begin{equation}
\label{1.1.6}
 \frac{fG-Fg}{G^2}.
\end{equation}
To deduce this from (\ref{1.1.4}) using $(F/G)\cdot G =F$, needs $F/G$ differentiable. Our proof avoids this.

Known differentiations begin with the identity function $F(x)=x$ with (\ref{eq 1.1.1}) equal to $1$.

If $n$  is { a constant, and } if $x \neq 0$ whenever the constant $ n\leq 0$,
{then the derivative of } $x^n$ { is }
\begin{equation}
\label{1.1.7}
 nx^{n-1}
\end{equation}
This is true when $n=1$. For other positive integers $n$ we use mathematical induction. Suppose the result true for the integer $n$. Then by (\ref{1.1.4}) the derivative of $x^{n+1}$, $=x^n\cdot x$, is
\[
\left(nx^{n-1}\right)\cdot x + x^n \cdot 1 = (n+1) x^n,
\]
and the result is true for $n+1$. Hence true for $1+1=2$, and so for $2+1=3$, $\ldots$. We step along the integers and (with patience) we arrive at any given positive integer.

Next, by an index law, if $x^0$ has any meaning then
\[
x^1 \cdot x^0 = x^{1+0} =x^1.
\]
When $x=x^1 \neq 0$ this gives $x^0=1$. When $x=0$ we have $0\cdot 0^0$, which gives nothing. Thus we prove the case when $n=0$ with $x \neq 0$.

For $q$ a positive integer let $x = y^q >0$. Then there is a $y>0$ with $y=x^{1/q}$, and by (\ref{1.1.5}),
\begin{eqnarray*}
1&=&\frac{dy}{dx}=\frac{dx}{dy} \cdot \frac{dy}{dx}, \vt
\frac{dy}{dx} &=& \left(\frac 1q\right)y^{1-q} = 
\left(\frac 1q\right)x^{\frac 1q-1}.
\end{eqnarray*} 
Thus (\ref{1.1.7}) is true for $n=1/q$, $x>0$. Using (\ref{1.1.5}) and a positive integer $p$,
\[
\frac{dx^{\frac pq}}{dx} = \frac{dy^p}{dx} = py^{p-1} \frac{dy}{dx} = \frac pq\cdot y^{p-q} = \frac pq \cdot x^{\frac pq -1},
\]
giving (\ref{1.1.7}) with $n=p/q$, $x>0$. For $n=-p/q$, where $p,q$ are positive integers, and $x>0$,
(\ref{1.1.6}) gives (\ref{1.1.7}) true for $n=-p/q$, $x>0$, since
\[
x^n = \frac 1{x^{\frac pq}},\;\;\;\;\;\;\;\;\;\;\;\;\frac{dx^n}{dx} =
\left(-\frac pq \right) \left(\frac{x^{\frac pq -1}}{x^{\frac{2p}q}}\right) = \left(-\frac pq\right)\left( x^{-\frac pq -1}\right).
\]
A ratio of two integers with positive denominator is called a \emph{rational} number; other numbers are \emph{irrational}. We thus see that (\ref{1.1.7}) is true for rational numbers $n$ with $x>0$. Later we will define $x^n$ ($x>0$) for $n$ irrational and again prove (\ref{1.1.7}).

It may happen that for some $a<b$, $u<v$, and some functions $y=G(x)$ in $[a,b]$ with values in $[u,v]$, and $x=F(y)$ in $[u,v]$ with values in $[a,b]$, we have $x=F(G(x))$ for all $x$ in $[a,b]$. Then $F,G$ are called \emph{inverse} functions. We have already had the example $x=y^q$, $y=x^{1/q}$. If both $F$ and $G$ have derivatives throughout their respective ranges, then by (\ref{1.1.5}),
\begin{equation}
\label{1.1.8}
\frac{dF}{dy} \cdot \frac{dG}{dx} = \frac{d(F(G(x))}{ds} =\frac{dx}{dx} =1,\;\;\;\;\;\;\;\;\;\;\;\;\frac{dG}{dx} = \frac 1{\frac{dF}{dy}}.
\end{equation}
This implies that if the derivatives of $F,G$ both exist, then neither can be $0$ in their respective ranges.

Note that if fairly transparent graph paper is used for the graph of $y=G(x)$ with the $y$-axis vertical and the $x$-axis horizontal and to the right, as is normal practice, then turning over the graph paper with the $y$-axis vertical and the $x$-axis horizontal and to the left, and then turning through 90 degrees, the paper will have the $x$- and $y$-axes interchanged and showing the graph of $x=F(y)$. A horizontal part of $y=G(x)$ in the range becomes a vertical line of the $x=F(y)$ graph, a jump or discontinuity.

The inverse process to the process of differentiation is called the \emph{calculus} (or \emph{Newton}) \emph{indefinite integral}. A function $F$ of points $x$ in $[a,b]$ is the \emph{primitive} or \emph{Newton indefinite integral} $\int f\,dx$ of a function $f$ in $[a,b]$, if $F'(x)=f(x)$ throughout $[a,b]$, using one-sided derivatives at $a$ and $b$. Hence from the rules for differentiation we find rules for Newton integration, which more general integrals also obey. First, by (\ref{1.1.3}) and the simple (\ref{eq 1.1.2}), {if }$f$ { is the derivative of }$F${ and }$c${ is a constant, then }
\begin{equation}
\label{1.1.9}
f\mbox{ is the derivative of }F+c.
\end{equation}
The Newton integral is not uniquely defined as we can add any constant to it. Also, by (\ref{1.1.3}) and the functions mentioned,
\begin{equation}
\label{1.1.10}
\int\left(mf+ng\right) dx = m\int f\,dx + n\int g\,dx + \mbox{constant}.
\end{equation}
Assuming that all integrals exist in (\ref{1.1.3}), (\ref{1.1.4}), we have the formula for integration by parts,
\begin{equation}
\label{1.1.11}
\int Fg\,dx 
=FG - \int Gf\,dx + \mbox{constant}.
\end{equation}
From (\ref{1.1.5}) and $h(K) = dF/dK$ we have
\[\int \frac{dF}{dK}\,dK =F(K(x))=\int \frac{dF}{dK}\cdot \frac{dK}{dx}\,dx +\mbox{ constant},\]
and, putting $h$ for $dF/dK$, we get
\begin{equation}
\label{1.1.12}
\int h(K)\,dK =\int h(K(x))\frac{dK}{dx} \,dx + \mbox{constant},
\end{equation}
which is the formula for \emph{integration by substitution}.

In calculus integration, to evaluate each integral it is normally treated individually by using the above rules in a systematic way, just as it has been calculated for centuries. We cannot give here details of all methods used, but we can show why various functions have appeared in the lists of calculated integrals.

We begin with integer powers of $x$. The integral of $x^n$ ($n \neq -1$) is found by differentiating $x^{n+1}$, giving $x^{n+1}/(n+1)\; +$ constant for the integral. This cannot give the integral of $x^{-1}$, so we look at its graph. There is an infinity at $x=0$, and the function is steadily increasing as $x$ rises in $x>0$, and in $x<0$. It appears intuitively that if we keep away from $x=0$, the area, under the curve from $x=1$ to $x=a>1$, is a well-defined function $A(a)$ of $a$, and similarly for the area under the curve from $x=1$ to $x=a$ in $0<a<1$. If we move to $x=a+h$, $a>1$, $A(a+h) -A(a)$ lies between $h$ times the least value of $x^{-1}$ and $h$ times the greatest value, both between $a$ and $a+h$,
\[
\frac h{a+h} \leq A(a+h) -A(a)\leq \frac ha ; \mbox{ and similarly }\; \frac ha\leq A(a) - A(a-h) \leq \frac h{a-h}.
\]
Dividing by $h$ and letting the positive $h$ tend to $0$, we have
\[
\frac{dA}{dx} = \frac 1a\;\;\;\;(x=a).
\]
To justify this argument we have to construct in some way the area under the curve. More generally, if we can construct the area under the curve $y=f(x)$ as $F(x)$, with $f$ having a special property called \emph{continuity}, then $dF/dx = f(x)$ and $F$ is an integral. Assuming the construction for $f(x) = x^{-1}$ we write $F(x)$ as $\log_e x$ or $\ln x$,
\[
\log_e a = \ln a = \int_1^a \frac 1 x\, dx  \;\;\;\;\;\;\;\;(a>1),
\]
where the symbols $\int_1^a x^{-1}\,dx$ represent (the value of the integral at $x=a$) minus (the value of the integral at $x=1$), the usual calculus notation for the \emph{definite integral}.

In (\ref{1.1.12}), if $K(a)=u$, $K(b)=v$, then for definite integrals,
\begin{equation}
\label{1.1.13}
\int_u^v h(K)\,dK = \int_a^b h(K(x))\frac{dK}{dx}\,dx.
\end{equation}
Keeping to $x\geq 0$, if, for some constant $a>1$ we have $K(x) = ax$, then for $b>1$,
\[
\int_a^{ab} \frac 1K\,dK = \int_1^b \frac 1{ax} \cdot a\,dx = \int_1^b \frac 1x\,dx.
\]
As the variable of integration is irrelevant when the integral is an area, we have
\begin{equation}
\label{1.1.14}
\log_e(ab) - \log_e(a) = \log_e b \;\;\;\;\;\;(a>1, b>1).
\end{equation}
Taking $\log_e 1$ to be the integral from $1$ to $1$, zero area, we define $\log_e 1$ to be $0$. Then when $0<a \leq 1$ or $0<b \leq 1$ or both, we again have (\ref{1.1.14}), and in particular,
\[
\log_e \frac 1a + \log_e a = \log_e \left( \frac 1a \cdot a\right) = \log_e 1 =0,
\]
\begin{equation}
\label{1.1.15}
\log_e \frac 1a = - \log_e a;\;\;\;\;\;\;\;\;
\log_e a^n = n \log_e a,
\end{equation}
for positive integers $n$ by using (\ref{1.1.14}) repeatedly, and for negative integers $n$ by also using the first part of (\ref{1.1.15}). For $n$-th roots $a^{\frac 1n}$ we have 
\[
n \log_e a^\frac 1n = \log_e {\left(a^\frac 1n\right)}^n =\log_e a \;\;\;\;\;\;(n \neq 0)
\]
and (\ref{1.1.15}) is true for $1/n$ replacing $n$. For $p/q$ replacing $n$ and $p,q$ non-zero integers,
\[
\log_e a^\frac pq = \log_e {\left(a^p\right)}^\frac 1q
= \frac 1q \log_e a^p = \frac pq \log_e a. 
\]
If $y = \log_e x$ we write $x = \exp y$ (the exponential function). From (\ref{1.1.14}), (\ref{1.1.15}),
\begin{equation}
\label{1.1.16}
\exp a \cdot \exp b = \exp(a+b),\;\;\;\;\exp 0 =1,
\end{equation}
\begin{equation}
\label{1.1.17}
\frac 1{\exp a} = \exp(-a),\;\;\;\;\;\;\;\;
(\exp a)^n = \exp (na)\;\;\;\;\;\;\;\;(n \mbox{ rational}).
\end{equation}
We write $e=\exp 1$, so that, for $n$ rational, $a>0$, $b$ rational,
\begin{equation}
\label{1.1.18}
\exp n = \left(exp 1\right)^n = e^n,
\end{equation}
\begin{equation}
\label{1.1.19}
a^b = \left(\exp(\log_e a)\right)^b = \exp(b \log_e a).
\end{equation}
For $b$ irrational or complex, and a wider definition of $\exp$, we define $a^b$ as $\exp(b \log_e a)$.

For differentiation, (\ref{1.1.8}) gives, from $y=\log_e x$,
\[
\frac d{dy}\exp y = \frac 1{\frac 1x} = x = \exp y,
\;\;\;\;\;\;\;\;\;\;\frac{d\left(a(x)^{b(x)}\right)}{dx} = 
\]
\begin{equation}
\label{1.1.20}
 =\;\;
 \frac d{dx} \exp (b(x) \log_e a(x)) 
=
a(x)^{b(x)} \left(b'(x) \log_e a(x) + b(x) \frac{a'(x)}{a(x)} \right).
\end{equation}
For those who know that the logarithm and exponential functions can be extended to functions of complex values with the above properties, and who know the exponential forms of $\sin x$ and $\cos x$, we can (with $\iota = \sqrt{-1}$) write 
\begin{eqnarray*}
\frac 1{1+x^2} &=& \frac 1{(x+\iota)(x-\iota)} \;\;=\;\; \frac  1{2\iota} \left(\frac 1{x-\iota} - \frac 1{x+\iota} \right), \vt
\int_0^a\frac 1{1+x^2}dx &=&\frac \iota 2 \left[ \log_e(x+\iota) - \log_e (x-\iota) \right]_0^a \vt
&=& \frac \iota 2 \log_e \left( \frac{a+\iota}{a-\iota} \cdot \frac{-1}\iota\right) \;\;=\;\; \frac \iota 2 \log_e \left( \frac{1-\iota a}{1+\iota a}\right) \;\;=\;\;y,\vt
1-\iota a &=& e^{-2\iota y} (1+\iota a), \vt
a&=& \frac{e^{\iota y}-e^{-\iota y}}{\iota \left(e^{\iota y}+e^{-\iota y}\right)} \;\;=\;\;\frac{\sin y}{\cos y} \;\;=\;\;\tan y,\vt
\int_0^a\frac 1{1+x^2}dx &=& \arctan a.
\end{eqnarray*}
Thus we can now integrate all fixed powers of $x$ and all polynomials (finite sums of constants, called \emph{coefficients}, times non-negative integer powers of $x$). For ratios of polynomial $p(x)$ and $q(x)$, these having no common factors, we factorise the denominator $q(x)$,
\[
q(x) = x^r(x-a)^s(x-b)^t \cdots, \;\;\;\;\;\;(r\geq 0,\;\;s>0,\;\;t>0, \ldots)
\]
where $r,s,t, \ldots$ are integers and where $a,b, \ldots $ are the different real or complex roots of $q(x)=0$. If the coefficients in $q(x)$ are all real, the complex roots occur in conjugate pairs $a, \bar a$, i.e.~for real $m,n$ and $a=m+\iota n$, then $\bar a = m-\iota n$, and the power s of $x-a$ is the same as the power of $x-\bar a$. The quadratic term with real coefficients $-2m$, $m^2+n^2$,
\[
(x-a)(x-\bar a) = x^2 - (a+\bar a) x + a \bar a,
\]
is raised to the power $s$. We now split up $p(x)/q(x)$ into partial fractions, a polynomial in $x$ plus fractions like (constant)/$(x-b)^t$, plus fractions like (constant)/$(x^2 + 2ux +v)^w$ for real constants $b,u,v$ and positive integers $t,w$, with $v-u^2>0$ as the quadratic has no real roots. The quadratic can be written $(x+u)^2 +(v-u^2)$. The polynomial integrates easily, and so do the first fractions, involving logarithmic terms when $t=1$, while the second fractions involve $\arctan$ terms.

Functions involving $\sqrt{1-y^2}$ with $-2<y<1$, for various $y$, can usually be evaluated using $y=\sin z$ or $\cos z$. And so on.
Texts on the calculus show how to evaluate many integrals; sometimes it is necessary to go further, to \emph{elliptic functions} and \emph{elliptic integrals} involving the square roots of
\[
ax^3+bx^2+cx +u,\;\;\;\mbox{ and of }\;\;\;
ax^4+bx^3+cx^2 +ux+v.
\]
The first type can be reduced to the second type by algebraic substitutions, and the second type gives rise to Legendre's normal forms, elliptic integrals of the first, second, and third kind, respectively:
\begin{eqnarray*}
&&\mbox{(i) }\;\;\int_0^x \frac{dx}{\sqrt{(1-x^2)(1-k^2x^2)}},\vt
&&\mbox{(ii }\;\;
\int_0^x \sqrt{\frac{1-k^2x^2}{1-x^2}} dx,\vt
&&\mbox{(iii) }\;\;\int_0^x \frac{dx}{(x^2-a)\sqrt{(1-x^2)(1-k^2x^2)}}.
\end{eqnarray*}
All letters apart from $x$ denote constants, $k$ being the \emph{modulus}, and usually we take $0<k<1$. For $w$ the integral (i) we write 
\[
x=\mbox{sn} (k,w)\;\mbox{ or }\;\mbox{sn}\, w,\;\;\;
\sqrt{1- \mbox{sn}^2\, w} = \pm \mbox{cn}\, w,\;\;\;\sqrt{1-k^2 \mbox{sn}^2\, w} =\pm \mbox{dn}\,w,
\]
the signs $\pm$ being chosen so that sn$\,0 =0$, cn$\,0=1=$dn$\,0$, 
and that cn$\,w$ and dn$\,w$ are differentiable everywhere. Writing $x=\sin \phi$ in (i), we have
\[
w = \int_0^\phi \frac{d\phi}{\sqrt{1-k^2 \sin^2 \phi}}
\]
and we write $\phi =$am$\,w$, the \emph{amplitude} of $w$. Thus
\[
\mbox{sn}\,w = \sin(\mbox{am}\,w),\;\;\;\;\;
\mbox{cn}\, w = \cos(\mbox{am}\,w).
\]
it is easily shown that if $1=$ sn$\,K$ then 
\[
\mbox{sn}(w+2K)= -\mbox{sn}\,w,\;\;\;
\mbox{cn}(w+2K)= -\mbox{cn}\,w,\;\;\;
\mbox{dn}(w+2K)= \mbox{dn}\,w.
\]
If $k' = $ dn$\,K = \sqrt{1-k^2}$, $\;\;\;1=$ sn$(k',K')$, then
\[
\mbox{sn}(w+2\iota K')= \mbox{sn}\,u,\;\;\;
\mbox{cn}(u+2\iota K')= -\mbox{cn}\,u,\;\;\;
\mbox{dn}(u+2\iota K')= -\mbox{dn}\,u,
\]
and sn, cn, dn are doubly periodic. Further,
\[
\frac{d(\mbox{sn}\,u)}{du} = \mbox{cn}\,u\,\mbox{dn}\,u,
\;\;\;\;
\frac{d(\mbox{cn}\,u)}{du} = -\mbox{sn}\,u\,\mbox{dn}\,u,
\;\;\;\;
\frac{d(\mbox{dn}\,u)}{du} = -k^2 \mbox{sn}\,u\,\mbox{cn}\,u.
\]
\section{Simple Definite Integration}\label{s2}
\subsection{Infinite Decimals and Real and Complex Numbers}
\label{s2.1}
(Omit this section if you are happy with your
 definition of real and complex numbers.) Here
 we assume known the relations $m=n$, $m>n$, $m \geq n$, and the operations $m+n$, $m-n$ (when $m \geq n$), and $mn$ on them, and their properties.
\begin{theorem}
\label{2.1.1}
A non-empty collection $S$, however large, of non-negative integers contains a minimum, say $m$, with the property that $n\geq m$ for all $n$ in $S$. If there is an integer $p$ such that $n \leq p$ for all $n$ in $S$, then $S$ contains a maximum, say $q$, with $n \leq q$ for all $n$ in $S$.
\end{theorem}
\proof
In the second part of the theorem $p$ might be in $S$. If not, then $p-1$ might be in $S$. If not, then $p-2$ might be in $S$. And so on. Eventually we arrive at the maximum $q$, after at most $p$ examinations of the contents of $S$. Similarly for the minimum, for which we begin at $0$ and proceed upwards.  

The practical man calculates a real number using decimals to base 10 to his required accuracy. We use the idea to define an \emph{infinite decimal} $x$, a sequence $x_0, x_1, x_2, \ldots $, or $(x_j)$, of non-negative integers and a set of rules, writing $x$ as
\[
x=x_0.x_1x_2x_3 \cdots .
\]
The first dot is the decimal point, $0\leq x_j \leq 9$ ($j \geq 1$), and $x_0$ is one of the sequence $0,1,2, \ldots , 105, \ldots $. As the base is 10 we can expect that
\[
10 x = \left(10 x_0 + x_1\right) . x_2x_3 \cdots \;.
\]
Thus for the special repeated integer decimal
$
y=0.999 \cdots$,
\[
10 y = 9.999 \cdots = 9+y ,
\]
so $y=1$, and we have evaluated a special important decimal. Again, for example, we ought to have $0.246999\cdots \;\;=$
\[
=\;\; 0.246000\cdots + 0.000999\cdots = 0.246000 \cdots +\frac{0.999\cdots}{10^3} = 0.247000\cdots \;.
\]
This is the background of the following rule. For infinite decimals $x=x_0.x_1x_2x_3 \cdots$ and $y=y_0.y_1y_2y_3 \cdots$ we write $x=y$ and $y=x$, if $x_j=y_j$ (all $j$) or of, for a particular integer $J\geq 1$, $x_j=9$ and $y_j=0$ (all $j \geq j$) while $y_{J-1}= x_{J_1}+1$ (so $x_{J-1} <9$ if $J \geq 2$), and $x_j=y_j$ (all  $j<J-2$, if any $j$ is in the range). Thus $x$ has a recurring 9 and $y$ a recurring 0, after a certain stage. We call $y$ a \emph{finite decimal}.

If neither of these two cases occur we write $x \neq y$. For the given $x$ let $I_n(x)$ be the integer $x_0$ written in the usual way (e.g.~105) followed by the integers $x_1,x_2, \ldots, x_n$ in that order, so that we could say that $I_n(x)$ is the integer part of $10^n x$. For example, $I_2(1.234\cdots)=123$. If $x_0=12$ then $I_2(x_0.345 \cdots) = 1234$. $I_2(0.0123 \cdots) =001$, written $1$, and $I_2(0.00123 \cdots) = 000$, written $0$. The use of $I_n$ saves many rules.

We write $x>y$ and $y<x$ if, for some integer $n$,
\[
I_j(x) = I_j(y) \;\;(j=0,1,2, \ldots , n-1)
,
\;\;\; \mbox{ but }\;\;\;I_n(x)>I_n(y).
\]
Then $I_m(x) >I_m(y)$ for all $m\geq n$. It follows that if $x<y$ and $y<z$ then $x<z$. If we do not have $x=y$ nor $x>y$ then, for some integer $J$,
\[
I_j(x) = I_j(y)\;\;\;(0\leq j<J),\;\;\;\;\;I_j(x) <I_j(y).
\]
By definition $y>x$, $x<y$. Thus if $x,y$ are two infinite decimals, either $x=y$ or $x<y$ or $x>y$. We write $x\geq y$, $y\leq x$, if $x=y$ or $x>y$. So if $x<y$ is false then $x\geq y$.
\begin{theorem}
\label{2.1.2}
Between any two unequal infinite decimals there is a finite decimal.
\end{theorem}
\proof
Let
$
x=x_0.x_1x_2x_3 \cdots \;\;\;<\;\;\;y=y_0.y_1y_2y_3 \cdots
$
where neither $x$ nor $y$ has repeated 9's form. Then, for some integer $J$,
\[
I_{J-1}(x) = I_{J-1}(y),\;\;\;\;\;I_J(x) < I_J(y),
\]
omitting the $I_{-1}$ when $J=0$. If
$
I_J(x) +2 \leq I_J(y)$, then
$I_J(z)=I_J(x)+1$, $z_j=0$ ($j>J$), will give a finite decimal  $z$ between $x$ and $y$. Otherwise $I_J(x)+1 =I_J(y)$ and we look at the $x_j$, $y_j$ ($j>J$). We can assume the most difficult case, when $y$ itself is a finite decimal. As $x$ does not have repeated 9's form there is a first $K>J$ with $x_K < 9$, and then we take the finite decimal $z$ between $x$ and $y$ with
\[
z_j = x_j\;\;\;(j<K),\;\;\;\;\;z_K=x_K +1 \leq 9,\;\;\;\;\; z_j=0 \;\;\;(j>K).
\]
\begin{theorem}
\label{2.1.3}
If $x,y$ are infinite decimals such that $x\leq y+10^{-n}$ for all positive integers $n$, then $x \leq y$.
\end{theorem}
\proof
If $x$ or $y$ or both have repeated 9's, we first change to the equal repeated 0's decimals. Then
\[
x=x_0.x_1x_2\cdots >y = y_0.y_1y_2 \cdots
\]
implies that there is an integer $N>0$ for which
\[
x_j=y_j\;\;\;(0\leq j<N),\;\;\;\;\;\;x_N>y_N.
\]
For some integer $J>N$, $y_J \neq 9$ and then
\[
x \leq y+10^{-J}=y_0.y_1y_2 \cdots y_{J-1}(y_J+1)y_{J+1} \cdots\;,
\]
contradicting $x_N>y_N$. Hence $x \leq y$.  \nproof

A collection or set $S$ of infinite decimals is \emph{empty} if no infinite decimal is in $S$, a rather trivial idea. If one or more infinite decimals lie in $S$, we say that $S$ is \emph{non-empty}. Such a set is said to be \emph{bounded above} by an infinite decimal $y$, and $y$ is an \emph{upper bound} of $S$, if $x\leq y$ for all $x$ in $S$. If no such $y$ exists we say that $S$ is \emph{unbounded above}. $S$ is said to have a \emph{supremum} (or \emph{least upper bound}) $u$ if $u$ is an upper bound of $S$ and if no nfinite decimal $v<u$ is an upper bound, i.e.~at least one member $x$ of $S$ lies in $v<x\leq u$. 

Going downwards instead of upwards, we say that a non-empty set $S$ of infinite decimals is 
\emph{bounded below} by an infinite decimal $w$, and $w$ is a \emph{lower bound} of $S$, if $x \geq w$ for all $x$ in $S$. $S$ is said to have an \emph{infimum} (or \emph{greatest lower bound}) $r$ if $r$ is a lower bound, and if no infinite decimal $s>r$ is a lower bound, i.e.~an $x$ of $S$ lies in $s>x\geq r$. We are not symmetrical as all $x \geq 0=0.000\cdots$. Later we have symmetry on including negative infinite decimals as well.

\begin{example}\label{ex 2.1.1}
If $y$ is an upper bound of $S$ and if $t>y$, then $t$ is an upper bound of $S$. If $w$ is an lower bound of $S$ and if $z<w$, then $z$ is an lower bound of $S$
\end{example}

\begin{example}\label{ex 2.1.2}
The set $(1,2,3,4)$ is bounded above with upper bounds $4,5,75.23$ etc. If $v<4$, one member $x$ of $S$ is in $v<x\leq 4$, namely $x=4$. Thus $4$ is the supremum. Similarly $1$ is the infimum. Generally, if $S$ is a finite set of infinite decimals, the supremum is the maximum, the greatest $x$ in $S$, and the infimum is the minimum, the least $x$ in $S$.
\end{example}

\begin{example}\label{ex 2.1.3}
For infinite decimals $a<b$, the set $S$ of all infinite decimals $x$ in $a \leq x \leq b$, is called a \emph{closed interval}, written $[a,b]$. As $x \leq b$, $b$ is an upper bound of $[a,b]$. So is $c$ if $c>b$. If $u<b$, an element $x$ of $[a,b]$ is in $u<c \leq b$, namely $x=b$. So $b$ is the supremum; similarly $a$  is the infimum, and both lie in $[a,b]$. There are three other kinds of intervals from $a$ to $b$, depending on whether or not each of $a,b$ lies in the interval. The set of all infinite decimals $x$ in $a<x<b$ is called an \emph{open interval}, written $(a,b)$. Here neither $a$ nor $b$ lies in the set, so we use Theorem \ref{2.1.2}. Again, $b$ is an upper bound, while if $u\leq a$ there is a finite decimal $x$ in $a<x<b$ with $u<x$, and $u$ cannot be an upper bound. If $a<u<b$ then a finite decimal $x$ lies in $u<x<b$, so $a<x<b$ and $u$ cannot be an upper bound. Thus $b$ is the supremum. Two half-open (or half-closed) intervals $a\leq x<b$, written $[a,b)$, and $a<x\leq b$, written $(a,b]$, also have supremum $b$ and infimum $a$, by using the above arguments appropriately.
\end{example}
\begin{example}\label{ex 2.1.4}
The set of even positive integers $2,4,6,8,\ldots$ is not bounded above. For if $u=u_0.u_1u_2 \cdots$ is a supposed upper bound, $u_0+1 \geq u$, and so is an upper bound. Either $u_0+2$ or $u_0+3$ is even, equal to $2n$ for an integer $n$, and $2n>u$. So $u$ is not an upper bound, and the set is unbounded above. 
\end{example}
\begin{theorem}
\label{2.1.4}
A non-empty set $S$ of infinite decimals has an infimum $r$. If $S$ is bounded above then $S$ has a supremum $u$. Each of $r,u$ can be in $S$ or out of it (see Example \ref{2.1.3}).
\end{theorem}
\proof
We prove the second part first, changing all repeated 9 decimals to the equal finite decimals. If, for some infinite decimal $y$, $x\leq y$ for all $x$ in $S$, then $x_0 \leq y_0$. By Theorem \ref{2.1.1} the collection of integers $x_0$, for all $x$ in $S$, has a maximum $u_0$ with $0\leq u_0\leq y_0$. We can now forget about $y$, it has served its purpose. Let $S_0$ be the set of all $x$ in $S$ with $x_0=u_0$. Then those $x$ have $0 \leq x_1\leq 9$, and by Theorem \ref{2.1.1} the $x_1$ have a maximum $u_1$ in $0\leq u_1 \leq 9$. Let $S_1$ be the set 
of all $x$ in $S_0$ with $x_1=u_1$. And so on. (This is a concealed induction.) In this way we define an infinite decimal $u=u_0.u_1u_2 \cdots$. Such a $u$ with repeated 9's is kept as it is until the proof ends. By construction $x \leq u$ for all $x$ in $S$, and $u$ is an upper bound. Let $v=v_0.v_1v_2 \cdots <u$. If $v_0<u_0$, each $x$ in $S_0$ has $v_0<x_0 = u_0$, $v<x\leq u$, and $v$ is not an upper bound. If $v_j =u_j$ ($0 \leq j <J$) but $v_J <u_J$, then each $x$ in $S_J$ has $v_j = x_j =u_j$ ($0 \leq j <J$), $v_J <x_J =u_J$, $v<x\leq u$, and $v$ is not an upper bound. Thus $u$ is the supremum. If necessary, we now transfer from repeated 9's to repeated 0's. Replacing $y$ by $0$ and going upwards, a similar construction gives the infimum.  \nproof

We now add, subtract, multiply and divide infinite decimals. If $x,y$ are infinite decimals let $(x+y)_n$ and $(xy)_n$ be the infinite decimals with
\begin{eqnarray*}
I_n(u) &=& I_n(x) + I_n(y),\vt
u_j &=&0\;\;(j>n),\vt
I_{2n}(v) &=&I_n(x)I_n(y),\vt
v_j&=&0\;\;(j>2n).
\end{eqnarray*}
Then $(x+y)_n$ and $(xy)_n$ rise with increasing $n$,  and both have suprema since
\[
(x+y)_n \leq x_0+y_0+2,\;\;\;\;\;(xy)_n \leq (x_0+1)(y_0+1)
\]
and by Theorem \ref{2.1.4}. Thus in an obvious notation we define
\[
x+y = \sup_n (x+y)_n,\;\;\;\;\;(xy)_n = \sup_n(xy)_n.
\]Note that $(x+y)_n=(y+x)_n$, $(xy)_n=(yx)_n$, so that $x+y=y+x$, $xy=yx$.

If $x=y$ we define $x-y$ to be $0$. If $x>y$ we define $(x-y)_n$ to be
\[
w=w_0.w_1w_2 \cdots,\;\;\;\;I_n(w) =I_n(x) -I_n(y),\;\;\;\;w_j=x_j\;\;(j>n).
\]
If $x>0$, there is a least integer $N\geq 0$ with $I_N(x)>0$ and we write
\[
(1/x)_n = \frac{10^{n+N}}{I_{n+N}(x)},
\]
obtained by division in the usual way of arithmetic, as a finite or recurring infinite decimal, recurring since after a finite stage the remainders are taken from the non-zero numbers 1 to $I_{n+N}(x)-1$, and we follow the remainder by $0$ and divide again. At that stage, if two remainders are the same, these and the other remainders between the first and second recurrence will be repeated constantly. The $(x-y)_n$, $(1/x)_n$ both fall with increasing $n$, bounded below by $0$, so that we take $x-y$ and $1/x$ to be the corresponding infima.

\begin{example}\label{ex 2.1.5}
The numbers $1, 1/2, 1/3, 1/4, \ldots$ can be written as a set of recurring infinite decimals that decrease steadily as $j$ increases. Prove that $0$ is the infimum.
\end{example}
The infinite decimals for $x=10^{-n}$ ($j=10^n, n=1,2, \ldots$) are members of the given set. By Theorem \ref{2.1.3}, $0$ is the infimum of the set, and as every $x \geq 0$, $0$ is the infimum of the $1/j$.

There are various algebraic rules concerning addition, subtraction, multiplication and division, that could be considered at this point, such as
\[
(x-y)+y =x,\;\;\;\;\;x+(y+z) = (x+y)+z,
\]
etc. They could be put as exercises.

The most constructive way to define negative infinite decimals is probably to use ordered pairs $\{x,y\}$ of infinite decimals $x,y$, so that $\{x,y\}$ is different from $\{y, x\}$. We define a \emph{real number} $R$ to be a set of all ordered pairs $\{x,y\}$ for which $\{u,v\}$ in in $R$ if and only if $x+v=y+u$ for some $\{x,y\}$ in $R$. If $x>y$ then $u=x-y$, $v=0$ satisfy the equation, and we say that $R$ is \emph{positive}. If $x>y$ then $u=0$, $v=y-x$ satisfy the equation and we say that $R$ is \emph{negative}. If $x=y$ then $u=0=v$ satisfy, and by the usual abuse of mathematical symbolism we write $R$ as $0$. If $x \geq y$ then $R$ is positive or $0$, and we call $R$ \emph{non-negative}. Similarly, if $x \leq y$, $R$ is negative or $0$, and we call $R$ \emph{non-positive}. For $R$ the set of ordered pairs $\{x,y\}$ and $T$ the set of ordered pairs $\{u,v\}$, we write $R>T$, $R<T$, $R=T$, respectively, if and only if
\[
x+v > y+u,\;\;\;\;x+v<y+u,\;\;\;\;x+v=y+u.
\]
We have avoided the use of $x-y$ (infinite decimal subtraction) when $x<y$. The terms ``bounded above (below)'', ``supremum'', ``infimum'', can now be defined for collections of real numbers just as they were defined for infinite decimals.
\begin{theorem}
\label{2.1.5}
A non-empty set $S$ of real numbers, bounded above by a real number, has a supremum. A non-empty set $S$ of real numbers, bounded below by a real number, has an infimum.\footnote{The conclusions in Theorem \ref{2.1.5} are called the \emph{Dedekind property} of the real numbers, named after the German mathematician Richard Dedekind (1831--1916) who discovered it. It is sometimes given as an axiom. Here we prove it from the construction.} 
\end{theorem}
\proof
We denote a supremum of $S$, a set of real numbers or of infinite decimals, by $\sup(S)$, an infimum by $\inf(S)$. For the first part, if there are non-negative real numbers in $S$ let $S^*$ be the set of all infinite decimals $u$ with $\{u,0\}$ in $R$, $R$ in $S$, and use Theorem \ref{2.1.4}, second  part.  If $\{\sup(S^*),0\}$ is in $T$ then $T=\sup(S)$. Next, if $R<0$ for all $R$ in $S$, we use the set $S^{**}$ of all infinite decimals $w$ with $\{0,w\}$ in $R$, $R$ in $S$, with Theorem \ref{2.1.4}, first part. If $\{0,\inf(S^{**}\}$ is in $T$ then $T=\sup(S)$. Similarly for Theorem \ref{2.1.5}, second part, or change every $\{x,y\}$ to $\{y,x\}$ and use Theorem \ref{2.1.5}, first part.
\nproof

Sometimes we use the convention that if $S$ is unbounded above then the supremum is $+\infty$; if $S$ is unbounded below it has the infimum $-\infty$.

We now have to deal with sequences. A \emph{sequence} is a function of the positive integers $1,2, \ldots , N$ for some integer $N$, called a \emph{finite sequence}, or a function of all positive integers, called an \emph{infinite sequence}. Thus $2,5,4,4$ is a finite sequence, while 
\[
1, \frac 12, \frac 13, \frac 14, \ldots , \frac 1n, \ldots
\]
is an infinite sequence. The general member of a sequence $s$ is often written $s_n$, and the sequence itself is written $(s_n)$.

\begin{theorem}
\label{2.1.6}
The rational numbers (ratios of integers) can be put as a disjoint sequence.
\end{theorem}
\proof
Every positive rational number occurs in the following sequence,
\[
\frac 11,\frac 12.\frac 21,\frac 13,\frac 22, \frac 31,\frac 14,\frac 23,\frac 32,\frac 41,\frac 15, \ldots \;.
\]
There are groups of fractions $p/q$, this particular fraction lying in the $(p+q-1)$-th group, 
which contains $p+q-1$ members. In the group we begin with $1/(p+q-1)$ and continue by raising the numerator by $1$ and lowering the denominator by $1$ for each successive fraction, till we reach $(p+q-1)/1$. In this first sequence every $p/q$ is repeated by $np/nq$ for $n=2,3, \ldots$, hence we omit from the sequence every $p/q$ for which $p$ and $q$ have a common factor greater than $1$, obtaining the sequence $(s_n)$. The final required sequence is then
\[
0,\;\;s_1,\;\;-s_1,\;\;s_2,\;\;-s_2, \ldots\;.
\]
\begin{theorem}
\label{2.1.7}
The real numbers in $[0,1)$ cannot all be put in sequence.
\end{theorem}
\proof
Using the representation by repeated decimals that do not have the repeated 9 form, let $(x_n)$ be a sequence of infinite decimals
\[
x_n\;\;=\;\;0.x_{n1}x_{n2}x_{n3} \cdots\;.
\]
Using Cantor's diagonal process let $y$ be given by the rule
\begin{eqnarray*}
y&=&0.y_1y_2y_3\cdots ,\vt
y_n&=&\left\{ \begin{array}{ll}x_{nn}+1&(0\leq x_{nn} \leq 7),\vt
1&(x_{nn}=8,9)\end{array}\right. \;\;\;\;(n=1,2, \ldots).
\end{eqnarray*}
Then no $y_n$ is $9$, $y$ cannot be a repeated 9 decimal, and $y$ cannot be in the sequence as it differs from $x_n$ in the $n$th decimal place. If the original sequence contained all of the real numbers in $[0,1)$, we would have a contradiction. \nproof

By induction we can find a second sequence disjoint from the first, by replacing the original sequence by $y,x_1,x_2,x_3, \ldots $ to find a $z$ disjoint from the new sequence; and so on. We have shown in some sense that the set of real numbers in $[0,1)$ is larger than the set of all rationals. In the former set closure we have a covering theorem, the \emph{Heine-Borel theorem}, also called \emph{Borel's covering theorem}, or the \emph{Borel-Lebesgue theorem} as Borel dealt with the case of a sequence of open sets and Lebesgue the general case.

\begin{theorem}
\label{2.1.8}
For $b>a$ let $[a,b]$ lie in the union of a family $\Gamma$ of open intervals. Then the union of a finite number of intervals of $\Gamma$ contains $[a,b]$.
\end{theorem}
\proof
First, $a$ lies in an interval of  $\Gamma$ Let $s$ be the supremum of the set $C$ of all points $c$ in $a<c \leq b$ such that $[a,c]$ lies in the union of a finite number of intervals of $\Gamma$. By Theorem \ref{2.1.5}, $s$ exists, and $a<s\leq b$. So there is an interval of $\Gamma$ that contains $s$, say $(u,v)$, and $u<s<v$, so that there is a $c$ in $C$ and in $(u,v)$. To the finite number of intervals of $\Gamma$ covering $[a,c]$ we add $(u,v)$, so that $[a,v)$ is covered by a finite number. By definition of $s$ we cannot have $s<b$, so that $s=b$.  \nproof

We now go a stage further, with ordered pairs $z=\{x,y\}$ of real numbers which satisfy the usual algebraic laws of addition and multiplication, together with a multiplication
\[
\{a,b\}\cdot\{c,d\} = \{ac-bd, ad+bc\},\;\;\;\;\;
\{0,1\}\cdot\{0,1\} = \{-1,0\}.
\]
We call $z=\{x,y\}$ a \emph{complex number} with \emph{real part} $x$ and \emph{imaginary part} $y$. We write $x$ for 
$\{x,0\}$, $\iota$ for $\{0,1\}$, $\iota y$ for $\{0,y\}$, and $x+\iota y$ for $z$. We define $\{x,y\}=\{u,v\}$ by $x=u$, $y=v$. It is enough to define $\{x,0\} =0 + \iota 0$, written $0$; then
\[
0 = (x + \iota y)\cdot (x-\iota y) = \{x^2+y^2,0\},\;\;\;\;\;
x^2+y^2 =0,\;\;\;\;\;x=0=y.
\]
Now let $z=x+\iota y$, $w=u+\iota v$, for $x,y,u,v$ real numbers. Then
\[
(xu+yv)^2 =x^2u^2+y^2v^2+2xyuv = (x^2+y^2)(u^2+v^2) - (xv-yu)^2.
\]
Writing $|z|$ as the non-negative square root of the non-negative number $x^2+y^2$, called the \emph{modulus} of $z$, we obtain
\[
|xu+yv| \leq |z| \cdot |w|,
\]
with equality when $z=0$ or $w=0$, or $z/w$ is real (from $xv=yu$). Then
\begin{eqnarray*}
|z+w|^2 &=& |(x+u) + \iota (y+v)|^2 \;\;=\;\; (x+u)^2 + (y+v)^2\vt
&=& |z|^2 +|w|^2 + 2(xu+yv) \;\;\leq \;\;|z|^2 +|w|^2 + 2|z|\cdot |w|,
\end{eqnarray*}
\begin{equation}
\label{eq 2.1.1}
|z+w| \;\;\leq \;\; |z|+|w|
\end{equation}
with equality when $xu+yv$ is non-negative and the other conditions are true, thus
\begin{equation}
\label{eq 2.1.2}\mbox{with equality when } z=0 \mbox{ or } \frac zw \mbox{ is a positive real number.}
\end{equation}
For $t=z+w$, $||t|-|w|| \leq |t-w|$. Interchanging $t$ and $w$, and then replacing $t$ by $z$, we have
\begin{equation}
\label{eq 2.1.3}
\left|^{ }_{ }{}|z|_{ }^{ }-|w|^{}_{}\right| \;\;\leq \;\; |z-w|.
\end{equation}
We can represent $z=\{x,y\}$ by a point with co-ordinates $(x,y)$ using rectangular co-ordinates in the plane. If $\theta$ is the suitably measured angle between the line joining $(0,0)$ and $(x,y)$ with the $x$-axis then $x=r \cos \theta$, $y=r\sin \theta$ with $r=|z|$, and $(r,\theta)$ the polar co-ordinates of $(x,y)$. The angle $\theta$ is known as the \emph{amplitude} or \emph{argument} of $z$.

The question now arises: of what use is this polar representation? So we look at multiplication. Also let $w=u,v)$ with polar representation $(s,\phi)$. Then 
\begin{eqnarray*}
zw &=& rs \left(\cos \theta \cos \phi - \sin \theta \sin \phi + \iota \left( \sin \theta \cos \phi + \cos \theta \sin \phi  \right)\right) \vt
&=& rs \left(\cos (\theta + \phi) + \iota \sin (\theta + \phi) \right),
\end{eqnarray*}
and the polar co-ordinates of $zw$ are $(rs, \theta + \phi)$. We multiply the moduli and add the amplitudes, and the second result gives a seeming connection with logarithms.

Further developments will need sums of infinite series, to which we now turn. Sometimes $s_n$ is the sum of the first $n$ terms of an infinite sequence
$(a_n)$,
\[
s_n = a_1+a_2+ \cdots +a_n,
\]
for example, the sum of a geometric sequence with 
\begin{eqnarray*}
a_j &=& r^j,\vt
rs_n &=& r^2 +r^3 + \cdots + r^{n+1}\;\; =\;\;s_n-r+r^{n+1}, \vt
(r-1)s_n &=& r(r^n-1), \vt
s_n &=& \left\{
\begin{array}{ll} \frac{r(r^n -1)}{r-1}, & (r \neq 1),\vt
n, & (r=1).
\end{array} \right.
\end{eqnarray*}

\textbf{GAP IN MANUSCRIPT}

\noindent
\textbf{Notes:} \textit{Infinite decimals were mentioned in T.J.~Bromwich (1926) pp.~394--401 and were used by P.~Dienes\footnote{Dienes was  Henstock's Ph.D.~supervisor. - P.M.} (1931) pp.~1--18. The results in Theorem \ref{2.1.5} are called the \emph{Dedekind property} of the real numbers, named after the German mathematician Richard Dedekind (1831--1916) who discussed it. Some authors give the result as an axiom. Here we prove it from the infinite decimal construction. Suprema and infima were discussed by R.~Dedekind (1909). The rationals were put in sequence by G.~Cantor (1875) who (1874 a) showed that all real numbers could not be put in a single sequence. E.~Borel (1895) gave his covering theorem for a sequence of open sets. See Heine (  ). Lebesgue (  ) gave the general case. So Theorem \ref{2.1.8} is called the Heine-Borel theorem, or Borel's covering theorem, or the Borel-Lebesgue covering theorem.}

\subsection{Cauchy, Riemann and Darboux}\label{s2.2}

A.L.~Cauchy (1789--1857) systematized in Cauchy (1821) the earlier definitions and theorems of the constructive calculus integral. His simplest construction was for a closed bounded (or compact) interval $[a,b]$ of real numbers $x$ in $a\leq x \leq b$ where $a<b$. A \emph{partition} of $[a,b]$ is a finite set of numbers $x_0, x_1, x_2, \ldots , x_n$ satisfying
\[
a=x_0<x_1< \cdots < x_{n-1}<x_n = b.
\]
The \emph{mesh} of the partition is the greatest of the $x_j-x_{j-1}$ ($j=1, \ldots ,n$). For a function $f$ defined everywhere in $[a,b]$ he took the limit, as the mesh tends to $0$, of the sum
\[
\sum_{j=1}^n f(x_j)(x_j-x_{j-1}) = f(x_1)(x_1-x_0) + f(x_2)(x_2-x_1) + \cdots +f(x_n)(x_n-x_{n-1}).
\]
In 1854 G.F.B.~Riemann (1826--1866) gave a slight but significant generalization, see Riemann (1868), replacing $f(x_j)$ by $f(\xi_j)$ where $\xi_j$ is arbitrary in $[x_{j-1},x_j]$ ($j=1, \ldots ,n$) giving
\begin{equation}
\label{eq 2.2.1}
\sum_{j=1}^n f(\xi_j)(x_j-x_{j-1}),
\end{equation}
a sum very near to that of Cauchy, since if
\begin{equation}
\label{eq 2.2.2}
x_{j-1}<\xi_j < x_j,\;\;\; f(\xi_j)(x_j-x_{j-1})
=f(\xi_j)(\xi_j-x_{j-1})+f(\xi_j)(x_j-\xi_{j}),
\end{equation}
$f$ being evaluated at an end of the interval concerned, not always at the right-hand end. Applying (\ref{eq 2.2.2}) wherever necessary, we have a partition with more points. The sum (\ref{eq 2.2.1}) is called a \emph{Riemann sum} and its limit is called the \emph{Riemann integral}. 

When the values of $f$ are real, J.B.~Darboux (1842--1917) in Darboux (1875) gave a definition which most calculus integration texts now use. He replaced $f(\xi_j)$ by the supremum $M(x_{j-1},x_j)$ of the values of $f$ in $[x_{j-1},x_j]$ to get \emph{upper sums}, and replaced $f(\xi_j)$ by the corresponding infimum $m(x_{j-1},x_j)$ to get \emph{lower sums}. For the graph of a positive $f$, the upper sum is the sum of areas of rectangles with bases the $[x_{j-1},x_j]$ and heights just enough to include the graph.  For the lower sum the rectangles lie just below or on the graph. As the mesh tends to $0$, if upper and lower sums tend to the same finite limit, it is the \emph{Riemann-Darboux integral}.
\begin{equation}
\label{eq 2.2.3}
\mbox{If the Riemann-Darboux integral exists then $f$ is bounded.}
\end{equation}
An upper sum and a lower sum have to be finite, so that the suprema and infima for the smaller intervals have to be finite. The maximum of the finite number of suprema is the supremum of $f$ on $[a,b]$, and the minimum of the infima is the infimum of $f$, and $f$ is bounded.
\begin{equation}
\label{eq 2.2.4}
\mbox{There is a bounded function with no Riemann-Darboux integral.}
\end{equation}
Let $f$ be $1$ for rational numbers, $0$ otherwise. Finite decimals are rational numbers and lie in every interval in the positive real line by Theorem \ref{2.1.2}, and in the negative real line by multiplying by $-1$. Hence the supremum of $f$ is $1$ in every interval. A non-zero rational multiple
of $\sqrt 2$ is an irrational number. So for $p<q$ we find a non-zero rational $r$ in $(p/\sqrt 2, q/\sqrt 2)$, $r/\sqrt 2$ is in $(p,q)$ and the infimum of $f$ in $[p,q]$ is $0$. All upper sums for $[a,b]$ are equal to $b-a$, all lower sums are $0$, and the Riemann-Darboux integral does not exist.

\begin{theorem}
\label{2.2.1}
For $f$ real-valued over $[a,b]$, the Riemann integral exists if and only if the Riemann-Darboux integral exists, with the same value, and then $f$ is bounded.
\end{theorem}
\proof
Let the Riemann-Darboux integral exist, so $f$ bounded there by (\ref{eq 2.2.3}). Then
\[
m(p,q) \leq f(x) \leq M(p,q)\;\;\;\;\;(a \leq p \leq x \leq q \leq b,\;\;\;\;p<q),
\]
every Riemann sum over $[a,b]$ lies between the corresponding upper and lower sums, and the Riemann integral exists with the same value. Conversely, let the Riemann integral exist with value $I$ over $[a,b]$. Then for small enough mesh all Riemann sums $S$ lie in $[I-1, +1]$. Then
\begin{eqnarray*}
S &\equiv & \sum_{j=1}^n f(\xi_j) (x_j-x_{j-1}) \;\;\;= \;\;\;f(\xi_k) (x_k - x_{k-1}) + T, \vt
f(\xi_k) & \leq & \frac{I+1-T}{x_k - x_{k-1}}, \vt
M(x_{k-1}, x_k) & \leq & \frac{I+1-T}{x_k - x_{k-1}},\;\;\;\;
(k=1,2, \ldots , n)
\end{eqnarray*}
finite, where we have kept $k$ and $T$ fixed and varied $\xi_k$. Similarly $m(x_{k-1}, x_k)$ is finite and $f$ is bounded in $[a,b]$. As no real number $v<M(p,q)$ can be an upper bound of the values of $f$ in $[p,q] \subseteq [a,b]$, there is a $\xi$ in $[p,q]$ with $f(\xi)>v$. Given $\ve>0$, we take
\begin{eqnarray*}
p&=& x_{j-1},\;\;\;q=x_j, \;\;\;v=M(x_{j-1}, x_j) - \ve,\;\;\;\xi=\xi_j\;\;\;(j=1,\ldots , n), \vt
S&>& \sum_{j=1}^n \left(M(x_{j-1}, x_j) - \ve(x_j-x_{j-1})\right) \;\;=\;\;\sum_{j=1}^n M(x_{j-1},x_j) - \ve(b-a).
\end{eqnarray*}
Similarly there are points $\eta_j$ with
\begin{eqnarray*}
&&x_{j-1} \leq \eta_j \leq x_j\;\;\;(j=1, \ldots ,n),\vt
&&S^* \equiv \sum_{j=1}^n f(\eta_j) (x_j-x_{j-1}) < \sum_{j=1}^n m(x_{j-1},x_j) + \ve(b-a).
\end{eqnarray*}
As the mesh tends to $0$, both $S$ and $S^*$ tend to $I$, so for small enough mesh the upper and lower sums lie in the closed interval centre $I$ and length $2\ve + 2\ve(b-a)$. As $\ve>0$ is arbitrary, the upper and lower sums tend to $I$, and the Riemann-Darboux integral exists equal to $I$, the value of the Riemann integral.  \nproof

From Theorem \ref{2.2.1} and (\ref{eq 2.2.4}) the Riemann integral cannot integrate all bounded functions. 

When $f$ is bounded and real-valued in $[a,b]$ we define the \emph{upper Riemann integral} of $f$ over $[a,b]$ to be the upper limit of Riemann sums (\ref{eq 2.2.1}) over $[a,b]$ as the mesh tends to $0$. As in Theorem \ref{2.2.1} this is the upper limit of upper sums over $[a,b]$, which can be called the \emph{upper Riemann-Darboux integral}, written $(R)\overline{\int}_a^b f\,dx$. Similarly the \emph{lower Riemann integral} $(R)\underline{\int}_a^b f,dx$ is the lower limit of (\ref{eq 2.2.1}) over $[a,b]$ as the mesh tends to $0$, and is the lower limit of lower sums, called the \emph{lower Riemann-Darboux integral}. Clearly the Riemann integral exists if and only if $f$ is bounded and the upper and lower Riemann integrals are equal.

\begin{theorem}
\label{2.2.2}
Let $f$ be real-valued and bounded in $[a,c]$ where $a<b<c$. Then
\begin{eqnarray}
\label{eq 2.2.5}
(R)\overline{\int}_a^c f\,dx 
&=&(R)\overline{\int}_a^b f\,dx  + (R)\overline{\int}_b^c f\,dx ,\nonumber \vt
(R)\underline{\int}_a^c f\,dx 
&=&(R)\underline{\int}_a^b f\,dx  + (R)\underline{\int}_b^c f\,dx ;
\end{eqnarray}
\begin{equation}
\label{eq 2.2.6}
D(a,c) \equiv
(R)\overline{\int}_a^c f\,dx 
-(R)\underline{\int}_a^c f\,dx  =D(a,b)+D(b,c);
\end{equation}
 As $b$ increases and as $a$ decreases:
\begin{equation}
\label{eq 2.2.7}
D(a,b)\mbox{ is non-negative and monotone increasing};
\end{equation}

\begin{equation}
\label{eq 2.2.8}
\mbox{If }
D(a,b)=0 \mbox{ then } D(p,q) =0 \mbox{ for all }[p,q] \subseteq [a,b].
\end{equation}
\end{theorem}
\proof
Given $\ve>0$, we take the Riemann sums $R_1, R_2$ over $[a,b]$, $[b,c]$ respectively, such that
\[
R_1> (R) \overline{\int}_a^b f\,dx - \frac \ve 3,\;\;\;\;\;
R_2> (R) \overline{\int}_b^c f\,dx - \frac \ve 3.
\]
As $R_1+R_2$ is a Riemann sum for $[a,c]$, and the meshes for $[a,b]$, $[b,c]$ and so $[a,c]$ are arbitrarily small, we have
\begin{eqnarray}\label{eq 2.2.9}
(R)\overline{\int}_a^c f\,dx &>& R_1+R_2 - \frac \ve 3 \;\;\;>\;\;\;(R)\overline{\int}_a^b f\,dx
+ (R)\overline{\int}_b^c f\,dx -\ve
 \nonumber \vt
\overline{\int}_a^c f\,dx &\geq &\overline{\int}_a^b f\,dx + \overline{\int}_b^c f\,dx.
\end{eqnarray}
For the opposite inequality take mesh $p$ less than $\ve/(8N)$ where $|f| \leq N$ in $[a,c]$, $0<v-u \leq p$, $u<b<v$, and points $\xi \in [u,v]$, $\eta \in [u,b]$, $\zeta \in [b,v]$, so 
\[
\left|f(\xi) (v-u) -f(\eta) (b-u) - f(\zeta) (v-b) \right|
\leq 2N(v-u) < \frac \ve 4.
\]
For the given $p$ and a suitable Riemann sum $R_3$, for $[a,c]$,
\begin{eqnarray}
(R)\overline{\int}_a^c f\,dx &<& R_3 + \frac \ve 4 
\;\;\;<\;\;\;
R_4+R_5 + \frac \ve 2 \nonumber \vt
&<&
(R)\overline{\int}_a^b f\,dx
+ (R)\overline{\int}_b^c f\,dx + \ve,
 \nonumber 
\end{eqnarray}
where $R_4$ and $R_5$ are the parts of $R_3$ given by intervals lying in $[a,b]$, $b,c]$ respectively, except that when an interval $[u,v]$ of the division over $[a,c]$ has $u<b<v$, we split into $[u,b]$ and $[b,v]$ with $f$ evaluated at $b$, and distribute to $R_4, R_5$ accordingly. Hence the inequality opposite to (\ref{eq 2.2.9}) and so (\ref{eq 2.2.5}), first part. Similarly the second part, then (\ref{eq 2.2.6}), (\ref{eq 2.2.7}), (\ref{eq 2.2.8}) since every $D(\cdot, \cdot)\geq 0$.

\begin{theorem}
\label{2.2.3}
The Riemann integral of $f$ over $[a,b]$ exists if and only if, given $\ve>0$, there is a $\delta>0$ such that, for every two divisions $D,D'$ of $[a,b]$ with mesh less than $\delta$,
\begin{equation}
\label{eq 2.2.10}
\left| (D)\sum f(t)(v-u) -(D')\sum f(t')(v'-u')\right| < \ve.
\end{equation}
\end{theorem}
\proof
Given $\ve>0$ and the Riemann integral $I$ of $f$ over $[a,b]$, let $\delta>0$ be such that for all divisions $D$ of $[a,b]$ with mesh less than $\delta$,
\[
\left| I - (D) \sum f(t)(v-u) \right| < \frac 12 \ve.
\]
For $D,D'$ divisions of $[a,b]$ with mesh less than $\delta$ and subtracting, we have (\ref{eq 2.2.10}).
Conversely, from (\ref{eq 2.2.10}) we prove integrability. When $\ve=1/n$ let $\delta = \delta_n$ and $s_n$ the sum over a division $D_n$ of $[a,b]$ with mesh less than $\delta_n$. Then for integers $m>n$ the meshes of $D_m$ and $D_n$ are less than $\delta_n$. Hence 
\[
\left|s_n-s_m\right| = \left|(D_n)\sum f(t)(v-u) - (D_m)\sum f(t')(v'-u') \right| < \frac 1n
\]
and $(s_n)$ is a fundamental, and so convergent, sequence, with limit $I$ say. As $m \rightarrow \infty$,
\[
\left| (D_n)\sum f(t)(v-u) -I \right| \leq \frac 1n.
\]
As $D_n$ is any division with mesh less than $\delta_n$, the Riemann integral exists with value $I$.  \nproof

\begin{equation}
\label{eq 2.2.11}
\mbox{\textbf{Not every calculus integral is a Riemann integral.}}
\end{equation}
We show an unbounded derivative. Let $F$ and its derivative $f$ be 
\begin{eqnarray*}
F(x) &=& x^2 \sin \left(\frac 1{x^2}\right) \;\;\;\;(x\neq 0),\;\;\;\;\;\;\;\;F(0)=0,\vt
f(x)&=& 2x \sin \left(\frac 1{x^2}\right) - \frac 2x \cos \left(\frac 1{x^2}\right)\;\;\;\;(x \neq 0).
\end{eqnarray*}
The first term tends to $0$ with $x$, but the second term and $f(x)$ oscillate unboundedly as $x \rightarrow \infty$.
\[
\frac{F(x) -F(0)}x = x \sin\left(\frac 1{x^2}\right) \rightarrow 0\;\;\;\;(x \rightarrow 0),\;\;\;\;\;\;\;f(0)=0.
\]
Thus the derivative exists for all $x$, unbounded in any interval containing $0$, and $f$ has no Riemann integral over such an interval.

\begin{theorem}
\label{2.2.4}
(Darboux) For $a<b$ let $F$ be differentiable everywhere in $[a,b]$ with derivative $f$. Let $f(a)<q<f(b)$. Then for some $\xi$ in $a<\xi<b$, $f(\xi) =q$
\end{theorem}
\proof
The differentiable $F$ is continuous, and therefore so is $G$ with
\[
G(x) = F(x) -qx,\;\;\;\;\;\;\;\;g(x) =f(x) -q
\]
its derivative. By Theorem (...?...) $G$ is bounded and obtains its infimum at a point $\xi \in [a,b]$. As $g(a)<0$, $G$ is falling at $a$, and as $g(b)>0$, $G$ is rising at $b$. So the infimum cannot be at $a$ nor at $b$, and $a<\xi<b$. As $G(\xi)$ is the infimum of the values of $G$,
\begin{eqnarray*}
&&
G(\xi+h)\geq G(\xi),\;\;\;\;\;\;G(\xi-h) \geq G(\xi),\vt && \vt
&&\frac{G(\xi +h) -G(\xi)}h \geq 0,\;\;\;\; 
\frac{G(\xi ) -G(\xi-h)}h \leq 0 \;\;\;\;(h \rightarrow 0), \vt &&\vt
&&g(\xi) =0
\end{eqnarray*}
as $g(\xi)$ is the limit of both fractions as $h \rightarrow 0+$, and $f(\xi) =q$. Similarly if $f(a)>q>f(b)$.
\begin{equation}
\label{eq 2.2.12}
\mbox{\textbf{Not every Riemann integral is a calculus integral.}}
\end{equation}
The function $f(x) =0$ ($x<\frac 12(a+b)$),
$f(x) =1$ ($x\geq \frac 12(a+b)$), has a Riemann-Darboux integral over $[a,b]$ equal to $b -\frac 12(a+b)$ $=\frac 12(b-a)$. But in Theorem \ref{2.2.4} with $q=\frac 12$, there is no $x$ with $f(x) = \frac 12$, and $f$ is not a derivative in the whole of $[a,b]$.

Riemann (1868), art.~5, shows that:

\noindent
\textbf{The necessary and sufficient condition for the Riemann integral of a bounded function to exist over $[a,b]$, is that the total length of the subintervals for which the oscillation is greater than any fixed positive number, is arbitrarily small.
}
\begin{equation}
\label{eq 2.2.13} 
\end{equation}

We can rewrite the condition as follows:

\noindent
\textbf{Given $\ve>0$, in every partition the sum of lengths of those intervals for which the oscillation is greater than $\ve$, is a value which tends to $0$ with the mesh.}
\begin{equation}
\label{eq 2.2.14} 
\end{equation}

For the interval $[u,v]$, the \emph{oscillation} is
\[
O(u,v) = M(u,v)-m(u,v) \geq 0,
\]
so that as the upper and lower sums tend to the same limit, their difference tends to $0$ with the mesh,
\begin{equation}
\label{eq 2.2.15}
\sum_{j=1}^n O(x_{j-1},x_j)(x_j-x_{j-1}) \rightarrow 0.
\end{equation}
Let $\sum'$ be summation over those $j$ for which $O(x_{j-1},x_j) >\ve$, and $\sum''$ over the rest. Then
\[\begin{array}{rll}
&&\ve \sum'(x_j-x_{j-1}) \rightarrow 0
\end{array}\]
with the mesh. Thus the condition is necessary. To show that it is sufficient, for $N$ the upper bound of $|f|$ in $[a,b]$, and given $\ve>0$, 
\[
\begin{array}{rll}
O(x_{j-1},x_j) &\leq &2N,\vt
0& \leq & \sum_{j=1}^n O(x_{j-1},x_j)(x_j-x_{j-1})\vt & \leq & 2N\sum^\prime(x_j-x_{j-1}) 
+\ve \sum^{\prime\prime} (x_j-x_{j-1}).
\end{array} 
\]
By the condition the first term tends to zero. The second term is less than or equal to $\ve(b-a)$. Since $\ve>0$ is as small as we please, as the mesh tends to zero, we have (\ref{eq 2.2.15}). Thus $f$ is Riemann-Darboux and so Riemann integrable. Note that we do not need the condition for all $\ve>0$, it is enough to use a sequence of $\ve>0$ decreasing to zero.

\begin{equation}
\label{eq 2.2.16}
\mbox{\textbf{There is a bounded derivative that is not Riemann integrable.}}
\end{equation}
This result goes beyond (\ref{eq 2.2.11}) but is far more complicated. Let $(x_n)$ be a sequence of numbers in $[0,1]$ such that every interval in $[0,1]$ contains an $x_n$, e.g.~the finite decimals or the rationals in $[0,1]$. Let $G_n$ be the open interval
\[
\left(x_n - \frac 1{2^{n+2}},\;\;\;x_n + \frac 1{2^{n+2}}\right)
\]
and let $G$ be the set of all points in all $G_n$ ($n=1,2, \ldots $). We call $G$ the \emph{union} of the $G_n$ and write
\[
G= \bigcup_{n=1}^\infty G_n.
\]
There is much overlapping, so that the points of $G$ do not have a simple formula, even when the $x_n$ are finite decimals. A point $x \in G$ is a point of $G_n$ for some $n$, and $(x,\;x_n + 2^{-n-2}) \subseteq G$. Let $y$ be the supremum of all $z$ with $(x,z) \supseteq G$. Now $x_n \leq 1$, so $z \leq 1 + 2^{-3}$ (\textsf{??? {Should be  $1+2^{-n-2}$? -- P.M.}}), $y$ is finite and $(x,y) \subseteq G$ while $y \notin G$. Similarly there is a $u \notin G$ with $(u,x) \subseteq G$. Thus $x \in (u,y) \subseteq G$, but we cannot enlarge $(u,y)$ in $G$. Beginning with $G_1$ let $H_1$ be the $(u,y) \supseteq G_1$.  If $G_n$ has a point in common with $H_1$ then by construction $G_n \subseteq H_1$. Let $j$ be the least integer such that $G_j \cap H_1$ is empty, and let $H_2$ be the $(u,y) \supset G_j$, $(u,y) \subseteq G$. Let $k$ be the least integer such that $G_k$ has no point in common with $H_1$ nor with $H_2$. And so on. Thus we construct a sequence $(H_m)$ of disjoint open intervals with union $G$. The $H_m$ are called \emph{connected components} of $G$. If $H_m = (u,y)$ then $[u+\ve, y-\ve]$ is an interval for small enough $\ve>0$, and lies in $H_m$, and so lies in the union of some of the $G_n$; which are open intervals. Thus the Heine-Borel theorem (Theorem \ref{2.1.8}) shows that $[u+\ve, y-\ve]$ is covered by a finite number of the $G_n \subseteq H_m$. Taking $\ve>0$ arbitrarily small, $y-u-2\ve$ and so $y-u$ are bounded above by the sum of the lengths of the $G_n \subseteq H_m$, and the sum of lengths of all the $H_m$ is not greater than the sum of lengths of all the $G_n$, namely $\frac 12$.

Let $0=a_0<a_1< \cdots <a_p=1$ be a partition of $[0,1]$, with $b$ the sum of the lengths of all intervals $(a_{j-1}, a_j)$ lying entirely within an $H_m$ or a union of abutting $H_m$ together with the common ends. Then $b \leq \frac 12$. The sum of the lengths of the other partition intervals, which we can call \emph{black intervals}, is not less than $1-\frac 12 = \frac 12$. We arrange an oscillation of $1$ on each of the $H_m$, so that (\ref{eq 2.2.14}) fails and $f$ is not Riemann integrable. 

For each $m$ let $I_m$ be the closest interval symmetrically at the centre of $H_m$ with length equal to the square of the length of $H_m$, and let the function $f$ be continuous, $1$ at the centre of $I_m$ and $H_m$, $0$ at the ends of $I_m$, always between $0$0 and $1$, for $m=1,2, \ldots$, and $0$ outside the $I_m$.

Then each black interval contains completely certain $I_m$ since the original $(x_n)$ has points in every interval.
To show that $f$ is a derivative we put 
\[
J(m,x) = I_m \cap [0,x],
\]
by convention taking the integral $0$ over the empty set. Let
\[
F(x) \equiv \sum_{m=1}^\infty \int_{J(m,x)} f(t)\,dt.
\]
If $x \in G$, then for some $p$, $x \in H_p$, and in a neighbourhood of $x$ all terms of the sum for $F$ are constant except possibly the term for $m=p$. Thus for $h>0$ and $h<0$,
\begin{eqnarray*}
\frac d{dx}\left(F(x)\right) &=& \frac d{dx} \int{J(p,x)} f(t)\,dt \;\;\;=\;\;\;\lim_{h \rightarrow 0} \frac 1 h \int_x^{x+h} f(t)\,dt \vt
&=& f(x) + \lim_{h \rightarrow 0} \frac 1 h \int_x^{x+h} \left(f(t)-f(x) \right)\,dt ,
\end{eqnarray*}
By continuity of $f$, given $\ve>0$ there is a $\delta>0$ such that
\begin{eqnarray*}
|f(t)-f(x)| &<& \ve\;\;\;\;\;(|t-x| < \delta), \vt
\left| \frac 1 h \int_x^{x+h} \left(f(t)-f(x) \right)\,dt\right| & \leq & \ve \;\;\;\;\;(|h| < \delta) ,
\end{eqnarray*}
and the derivative of $F$ is $f(x)$ since $\ve>0$ can be made arbitrarily small.

Finally $F'(x) =0$ when $x \notin G$. This is obvious when $x$ is an endpoint of two $H_m$. Otherwise let $x$ be in the interval $I \subseteq [0,1]$ with $I \cap I_m$ not empty for some $m$. Let $L(K)$ denote the length of the interval $K$ with $L(\emptyset) =0$. Put $S_m = L(H_m) \leq \frac 12$. As $I$ contains at least that part of $H_m$ that lies on one side of $I_m$ up to that end-point of $I_m$,
\begin{eqnarray*}
L(I \cap H_m) &\geq & \frac 12 \left(S_m -S_m^2\right) \;\;\; \; \;\;\;\geq \;\;\; \frac 14 S_m, \vt
L(I \cap I_m) &\leq & L(I_m) \;\;\;=\;\;\; S_m^2 \;\;\;\leq \;\;\; 16 L (I \cap H_m)^2.
\end{eqnarray*}
Summing over the $m$ with $I \cap I_m$ not empty,
\[
\sum L(I \cap I_m) \leq \sum 16 L(I \cap H_m)^2 \leq 16 L(I) \sum  L(I \cap H_m \leq 16 L(I)^2.
\]
Now $0 \leq f \leq 1$ and $f=0$ in $H_m$ except in $I_m$. If $I=(x,y)$ or $(y,x)$ and $K(m,x,y) = I \cap I_m$,
\begin{eqnarray*}
\left| \sum_{m=1}^\infty \int_{K(m,x,y)} f(t)\,dt\right|
& \leq & \sum_{m=1}^\infty L(K(m,x,y))\;\;\leq \;\;16L(I)^2 \;\;=\;\;16(y-x)^2;\vt
\left| \frac{F(y)-F(x)}{y-x}\right| &\leq & 16|y-x| \;\;\;\rightarrow \;\;\;0
\end{eqnarray*}
as $y \rightarrow x$, and $F'(x) =0$ ($x \notin G$). Thus $F'=f$ everywhere in $[0,1]$, finishing the proof. \nproof

In the language of Lebesgue and measure theory, which appeared 30 years later than (\ref{eq 2.2.13}), this condition is that the bounded function is continuous almost everywhere. (\ref{eq 2.2.16}) was first proved by Volterra (1881). Other examples followed, the simplest being the excellent construction by Goffman (1977), reproduced here.

In order to deal with even the calculus integration there is a need to go beyond Riemann integration. Darboux's definition, so useful in practice, is a barrier against generalizations such as the gauge integral, so that in the next two sections we consider definitions extending Riemann integration, before we deal with the gauge integral. Historically they appeared more than ten years after the gauge integral. 

\begin{example}
\label{ex 2.2.1}
Evaluate
\[
\lim_{j\rightarrow \infty} \left\{ \frac 1{j^2+1} + \frac 2{j^2+4} + \cdots + \frac j{j^2+j^2}\right\}.
\]
(University of Ulster, 1986, M112)
\end{example}
\textbf{Hint:} Write the expression inside $\{\cdots \}$ as a sum of
\[
\frac 1 j \cdot \frac rj \left\{ 1 + \left(\frac rj \right)^2 \right\}^{-1}
\]
for $r=1, \ldots , j$, to show that it is a Riemann sum to integrate $x(1+x^2)^{-1}$ from $0$ to $1$, and hence evaluate it.

\begin{example}
\label{ex 2.2.2}
Similarly find
\[
\lim_{j\rightarrow \infty} 
\sum_{k=1}^j \frac{(2k-1)^5}{2^5j^6}.
\]
(New University of Ulster, 1975, M111)
\end{example}

\begin{example}
\label{ex 2.2.3}
If $r$ is fixed in $0<r<1$ and if $f(0)=0$ and $f(x)=r^n$ ($r^n<x \leq r^{n-1}$), $n=1, 2, \ldots$, prove that $f$ is Riemann integrable over $0\leq x \leq 1$ to $r/(1+r)$, even though $f$ has an infinity of discontinuities.
(New University of Ulster, 1972, M112)
\end{example}

\begin{example}
\label{ex 2.2.4}
Let $(s_j)$ be a not necessarily monotone sequence of points in $[0,1]$ with infimum $s$ and let
\[
f(x) = \left\{
\begin{array}{lll}&\sum\left\{ \frac{s_j}{j^2}:\;{s_j<x}\right\} &(x>s), \vt
&0& (x\leq s). \end{array} \right.
\]
Prove that the Riemann integral of $f$ exists over $[0,1]$ with value
\[
\sum_{j=1}^\infty \frac{s_j(1-s_j)}{j^2}.
\]
(New University of Ulster, 1972, M213)
\end{example}

\subsection{Infinite Intervals}\label{s2.3}
Having defined the Riemann integral of a bounded real-valued function over a closed bounded interval, we extend it to infinite intervals beginning with $[a, \infty)$, the set of all real numbers $x \geq a$.
The obvious definition of the integral was used for centuries, being
\begin{equation}
\label{eq 2.3.1}
\int_a^\infty f\,dx = \lim_{b\rightarrow \infty} \int_a^b f\,dx,
\end{equation}
and systematized by Cauchy (1823), Lec.~24. See also de la Vall\'ee Poussin (1982 a,b). Proofs of some properties of such an integral, the limit of the limit of Riemann sums, so a double limit, are difficult. The limit process given here is easier to handle, while a more general limit process is given in section \ref{s3.1}. If
\begin{equation}
\label{eq 2.3.2}
a=x_0<x_1<x_2 < \cdots <x_{n-1}<x_n=b
\end{equation}
are real numbers, the finite collection of $[x_{j-1},x_j]$ ($j=1,2, \ldots ,n$) has been called a \emph{partition} $P$ of $[a,b]$, and the greatest of the $x_j-x_{j-1}$, the \emph{mesh} $|P|$ of $P$. A collection of interval-point pairs $([x_{j-1},x_j], \xi_j)$ is a \emph{division} $D$ of $[a,b]$ based on $P$ if
\begin{equation}
\label{eq 2.3.3}
x_{j-1} \leq \xi_j \leq x_j\;\;\;\;\;(j = 1,2, \ldots ,n).
\end{equation}
Given a function $f$ defined on $[a,b]$, the \emph{Riemann sum} for $D$ and $f$ is
\begin{equation}
\label{eq 2.3.4}
\sum_{j=1}^n f(\xi_j)(x_j-x_{j-1}) \equiv (D) \sum f(\xi)(v-u)
\end{equation}
where $[u,v]$ denotes the interval and $\xi$ the associated point in the second sum.

Let $\R$ be the class of all real-valued functions on $[a, \infty)$, Riemann integrable over $[a,b]$ for each $b>a$, such that the finite limit in (\ref{eq 2.3.1}) exists. A real-valued function $g>0$, strictly decreasing on $(0,b)$ for some $b>0$, is called a \emph{regulating function} if $g(x)$ increases without bound (tends to infinity) as $x$ tends to $0$ through positive values ($x$ \emph{tends to $0$ from above}, $x \rightarrow 0+$). A function $f$ on $[a, \infty)$ is said to be \emph{regulated by} $g$, if $g$ is a regulating function and if a number $I$ satisfies the following condition. For each $b_n>a$ of a sequence strictly increasing to infinity, each sequence $(P_n)$ of partitions, $P_n$ of $[a, b_n]$ with 
\[
|P_n| \rightarrow 0,\;\;\;\;\;\;0<\frac{b_n}{g(|P_n|)} \rightarrow 0
\]
as $n \rightarrow \infty$, each sequence $(D_n)$ of divisions, $D_n$ based on $P_n$, with Riemann sum $R_n$ for $D_n$ and $f$, then $R_n \rightarrow I$ as $n \rightarrow \infty$.

\begin{theorem}
\label{2.3.1}
If $f$ is regulated by $g$ for $[a,\infty)$, then $f$ is Riemann integrable over $[a,b]$ for each $b>a$, and the integrals tend to a finite limit as $b \rightarrow \infty$.
\end{theorem}
\proof
If $f$ is regulated by $g$ and if $|P'_n| \leq |P_n|$ ($n=1,2, \ldots$), then
\[
0<\frac{b_n}{g(|P'_n|)} \leq \frac{b_n}{g(|P_n|)}  \rightarrow 0
\]
as $n \rightarrow \infty$, and we can take 
\[
\ve_n \rightarrow 0+,\;\;\;\;\;\frac{b_n}{g(\ve_n)} \rightarrow 0,\;\;\;\;\;|P_n| \leq \ve_n.
\]
Taking $R_n$ within $2^{-n}$ of $(R)\overline{\int}_a^{b_n} f\,dx$, as $R_n \rightarrow I$,
\[
(R)\overline{\int}_a^{b_n} f\,dx \rightarrow I;\;\;\;\;\;\;
\mbox{ and similarly }(R)\underline{\int}_a^{b_n} f\,dx \rightarrow I,
\]
and the difference $D(a,b_n)\geq 0$ between the integrals tends to $0$ as $n \rightarrow \infty$. But by (\ref{eq 2.2.7}), $D(a,b_n)$ is monotone increasing in $n$. Hence every $D(a,b_n)=0$, and by (\ref{eq 2.2.7}) again, $D(a,b) = 0$ ($b>a$), $f$ is Riemann integable over $[a,b]$ for every $b>a$, and the Riemann integral for $b=b_n$ tends to $I$. As $(b_n)$ is an arbitrary strictly increasing sequence tending to infinity it follows that the Riemann integral over $[a,b]$ tends to $I$ as $b \rightarrow \infty$, which is (\ref{eq 2.3.1}) with $I$ as the value of the integral over $[a, \infty)$. Thus Theorem \ref{2.3.1} is proved. \nproof

Clearly we have similar results for the integral over $[a,b]$ when $a \rightarrow -\infty$. 

Haber and Shisha (1974), p.~3, states that every $f$ in $\R$, i.e.~real-valued $f$ Riemann integrable over $[a,b]$ for every $b>a$, with (\ref{eq 2.3.1}) true, is regulated by some regulating function $g$. I have not been able to prove this.

For some $f$ we need not use regulating functions. We say that $f$ is simply integrable over $[a, \infty)$ if there is a number $I$ such that, given $\ve>0$, there are numbers $B>a$, $\delta>0$, with the property that if $b>B$ and $P$ is a partition of $[a,b]$ with $|P|<\delta$, and $R$ is a Riemann sum for $f$, based on $P$, then $|R-I|<\ve$. Thus $f$ is simply integrable if the Riemann sums using partitions $P$ of $[a,b]$ tend to a unique finite limit when $b \rightarrow \infty$ and $|P| \rightarrow 0$ simultaneously. As above, $f$ is in $\R$. (Theorem \ref{2.3.1})

A more complicated definition involves the difference between the upper sum and the lower sum for a partition $P$ of an interval $[a,b]$, written $O(f,P)$.

A function $f$ is said to satisfy the \emph{uniform Riemann condition} if, given $\ve>0$, there is a $\delta(f)$ independent of $b>a$, such that if $P$ is a partition of $[a,b]$, then
\begin{equation}
\label{eq 2.3.5}
|P|<\delta(f) \mbox{ implies } O(f,P) < \ve.
\end{equation}
Given $\ve>0$, there are a $B(f)>a$ and a $\delta'(f)>0$ such that whenever $b'>b>B(f)$ and $P$ is a partition of $[b,b']$, $|P|< \delta'(f)$, and $R$ is any Riemann sum for $f$, based on $P$, then
\begin{equation}
\label{eq 2.3.6}
|R|<\ve.
\end{equation}
If we weakened this condition by letting $\delta$ to depend on $b$, and $\delta'$ to depend on $b, b'$, then (\ref{eq 2.3.5}) would give the Darboux condition for Riemann integrability over $[a,b]$ and all $b>a$, and (\ref{eq 2.3.6}) would then give the Cauchy convergence condition to ensure (\ref{2.3.1}), namely, that the difference between the integral over $[a,b]$ and the integral over $[a,b']$ tends to $0$ as $B \rightarrow \infty$.

\begin{theorem}
\label{2.3.2}
A function $f$ is simply integrable if and only if $f$ satisfies the uniform Riemann condition.
\end{theorem}
\proof
To show that a simply integrable function $f$ satisfies the uniform Riemann condition, 
given $\ve>0$, let $B$, $\delta$ be as in the first definition. By Theorem \ref{2.3.1}, $f$ is Riemann integrable on $[a,B]$, so that there is a $\delta'>0$ such that if $P$ is any partition of $[a,B]$ with $|P|<\delta'$ then $O(f,P) < \ve$. If $a<b<B$ and $P'$ any partition of $[a,b]$ with $|P'|<\delta'$, then $\delta'$ might not at first sight be small enough to ensure the required inequality. We use a simple special argument to obtain what we require. Let $R,R'$ be Riemann sums for $[a,b]$ and based on $P'$, and let $R''$ be a Riemann sum for $[b,B]$. Then $R+R''$ and $R'+R''$ are two Riemann sums for $[a,B]$ and so, for the corresponding partition $P''$ of $[a,B]$, arranged to have $|P''|<\delta'$,
\[
|R -R'| = |(R+R'')-(R'+R'')| \leq O(f,P''),\;\;\;\;\;O(f,P') \leq O(f,P'') < \ve.
\]
Let $\delta'' = \min(\delta, \delta') >0$. If $b>a$ and $P$ a partition of $[a,b]$ with $|P|<\delta''$ then $O(f,P) < \ve$ when $b \leq B$. When $b>B$, every Riemann sum $R$ using $f$ and based on $P$ satisfies
\[
\left| R - \int_a^\infty f, dx\right| < \ve,\;\;\;\;\;\;\;\; O(f,P)<2\ve,
\]
on taking $R$ near to the upper sum and then near to the lower sum, and we have (\ref{eq 2.3.5}) on replacing $\ve$ by $\ve/2$. Then $f$ in $\R$ gives (\ref{eq 2.3.6}).

Conversely, if the uniform Riemann condition holds, so does the weak condition and $f$ is in $\R$. Put
\[
I=\int_a^\infty f\,dx.
\]
Given $\ve>0$, let $B(f)$, $\delta(f)$ be as in the definition of the uniform Riemann condition with $\ve$ replaced by $\ve/3$, and $\delta(f)<1$. Put $B'=B(f)+1$. If $b>B$ let the partition $P$ of $[a,b]$ be given by
\[
a=x_0<x_1< \cdots <x_n=b,\;\;\;\;\;\;|P|<\delta(f) <1.
\]
Let $m$ be the greatest integer such that $x_m\leq B'$ and let $P_1$ be the partition of $[a,x_n]$ by points $x_0, x_1, \ldots , x_n$. For the Riemann sum $R$ based on $P$ let $R_1$ be the sum of the first $m$ terms of $R$. Since $|P|<1$, so that $x_m>B(f)$, (\ref{eq 2.3.6}) implies that
\[
\left|I-\int_a^{x_m} f\,dx\right| < \frac \ve 3,\;\;\;\;\;\;\left|R-R_1\right| < \frac \ve 3,
\]
and (\ref{eq 2.3.5}) implies that
\[
\left|R_1-\int_a^{x_m} f\,dx\right| < \frac \ve 3.
\]
Thus we prove the theorem from $|R-I|<\ve$. \nproof

Clearly all these results can be given for $(-\infty,b]$, the set of all $x$ in $x \leq b$, using
\[
\int_{-\infty}^b f\,dx = \lim_{a \rightarrow -\infty}\int_a^b f\,dx,
\]
and for $(-\infty, \infty)$ (the complete real line), using
\[
\int_{-\infty}^\infty f\,dx = \lim_{a \rightarrow -\infty}\lim_{b \rightarrow \infty}\int_a^b f\,dx.
\]

\subsection{An Infinity of a Function in a Bounded Interval}
\label{s2.4}
The next step in extending Riemann's integral is to integrate over a point where a function has an infinity. For example, the integral of $f(x) = 1/\sqrt x = x^{-\frac 12}$ in $x>0$ is $(x^{\frac 12}/\frac 12) +$constant,
\[
\int_\ve^1 f(x)dx = 2 - 2\ve^\frac 12
\;\;\;\;(0<\ve<1),
\]
avoiding the infinity at $x=0$. As $\ve$ gets smaller and smaller ($\ve$ shrinks) so does $\ve^\frac 12$, and the integral tends to $2$. Usually we define the integral 
\begin{equation}
\label{eq 2.4.1}
\int_{0+}^1 f\,dx \equiv \lim_{\ve \rightarrow 0+} \int_\ve^1 f\,dx
\end{equation}
if the limit exists. In the present case it gives $2$ for the integral over $[0,1]$.In other words we assume that the integral is continuous at $0$ where the infinity of $f$ occurs. This was systematised by Cauchy (1821), Lec.~25. We can replace $[0,1]$ by $[a,b]$ with $a<b$ real, and
\begin{equation}
\label{eq 2.4.2}
\int_{a+}^b f\,dx \equiv \lim_{\ve \rightarrow 0+} \int_{a+\ve}^b f\,dx,\;\;\;\;\;\;
\int_{a}^{b-} f\,dx \equiv \lim_{\ve \rightarrow 0+} \int_{a}^{b-\ve} f\,dx.
\end{equation}
In (\ref{eq 2.4.1}) and (\ref{eq 2.4.2}), as in (\ref{eq 2.3.1}), the integral is defined by a double limit. To avoid this, we first use the dominated integral of Osgood and Shisha (1976b, 1977) and Lewis and Shisha (1983) which is defined for the interval $[0,1]$. In Section \ref{s3.2} we will go further.

Let $f$ be a real or complex valued function on $(0,1]$. The \emph{dominated integral} of $f$ over $[0,1]$ is a number $I$ having the following property. Given $\ve>0$, there exist $\delta$, $\chi$ in $0<\delta<1$, $0<\chi <1$, such that

\begin{equation}
\label{eq 2.4.3}
\left| I - \sum_{j=1}^n f(\tau_j) (t_j - t_{j-1}) \right| < \ve
\end{equation}
whenever $0<  t_0<t_1< \cdots < t_n=1$, $t_0<\chi$, $t_{j-1}\leq \tau_j \leq t_j$, $t_{j-1}t_j^{-} > 1-\delta$, $j=1,2, \ldots ,n$.

\begin{theorem}
\label{2.4.1}
If $f$ is complex-valued on $(0,1]$ and has a dominated integral $I$ over $[0,1]$, then $\Re f$, the real part of the value of $f$, has a dominated integral $\Re I$; and $\Im f$, the imaginary part of the value of $f$, has a dominated integral $\Im I$.
\end{theorem}
\proof
These follow from (\ref{eq 2.4.3}) since, for complex $z$,
\[
\Re \left(I - \sum_{j=1}^n f(\tau_j) (t_j-t_{j-1}) \right) = \Re I - \sum_{j=1}^n \Re f(\tau_j)(t_j-t_{j-1}),
\]
$|\Re z| \leq |z|$, $|\Im z| \leq |z|$. \nproof

\begin{theorem}
\label{2.4.2}
Let $f$ have a dominated integral in $[0,1]$. Then $f$ satisfies the Riemann condition for $f$, namely that $f$ is defined on $(0,1]$ and bounded on $[\theta 1]$ and Riemann integrable, for each $\theta $ in $0<\theta <1$, and that, given $\ve>0$, there is a $\theta$ in $0<\theta<1$ such that, for $)(a,b)$ the oscillation of $f$ in $[a,b]$,
\begin{equation}
\label{eq 2.4.4}
OS_{j=1}^{n ???}(f;t_0,t_1, \ldots , t_n) \equiv
\sum_{j=1}^n O(t_{j-1},t_j) (t_j-t_{j-1}) < \ve
\end{equation}
whenever $0<t_0<t_1< \cdots < t_n=1$, $t_{j-1}t_j^{-1}>1-\theta$, $t_0<\theta$.
Conversely, if $f$ satisfies the Riemann condition for $f$, then $f$ has a dominated integral in $[0,1]$.

\end{theorem}

\proof
The only restraint on $\tau_j$ is that $t_{j-1}\leq \tau_j \leq t_j$. If $f$ is real-valued we choose $\tau_j$ so that $f(\tau_j)$ is as near as we like to its supremum, while keeping the other $\tau_k$ fixed, and still have the same inequality. Hence the supremum of the values of $f$ in $[t_{j-1},t_j]$ is finite. Being true for $j=1,2, \ldots ,n$ in turn, we see that $f$ is bounded above in $[t_0,1]$.
Similarly $f$ is bounded below in $[t_0,1]$ and so is bounded there, and $t_0$ can be arbitrarily small. Further, by the same argument we can replace each $f(\tau_j)$ in (\ref{eq 2.4.3}) by its supremum, and also by its infimum, on incurring the slight penalty for each, of replacing $<\ve$ by $\leq \ve$. Taking the difference between the two inequalities and using the oscillation $O(t_{j-1}, t_j)$ of section \ref{s2.2}, we have (\ref{eq 2.4.4}) with $\ve$ replaced by $\leq 2\ve$.

For complex-valued $f$ we first use Theorem \ref{2.4.1}. and so obtain (\ref{eq 2.4.4}) with $2\ve$ and with $O(t_{j-1}, t_j)$ the oscillation first of $\Re f$ and then of $\Im f$, respectively $O_1$ and $O_2$. For $a= \Re f(\tau_j) - \Re f(\tau_{j'})$, $b= \Im f(\tau_j) - \Im f(\tau_{j'})$,
\[
|a+ \iota b| \leq |a| + |b| = O_1 + O_2,\;\;\;\;\;\;
\sup \left| f(\tau_j) - f(\tau_{j'})\right| \leq O_1+O_2.
\]
The supremum is now the oscillation of $f$ in $[t_{j-1}, t_j]$, so that we have (\ref{eq 2.4.4}), except for $\leq 4\ve$ in place of $<\ve$. All we need to do now is to replace the original $\ve$ by $\ve/5$ to obtain (\ref{eq 2.4.4}) exactly.

To show that $f$ is Riemann integrable over $[\theta, 1]$ for each $\theta$ in $0<\theta <1$ we need only observe that the condition $t_{j-1}t_j^{-1} > 1-\delta$ is the same as
\begin{equation}
\label{eq 2.4.5}
t_{j-1} > T_j - \delta t_j,\;\;\;\;\;\;|t_j -t_{j-1}| < \delta t_j
\end{equation}
so that each partition of $[\theta, 1]$ with mesh less than $\delta t_j$ can be used with a suitable partition of $[0, \theta]$ that satisfies the conditions, and thus for two partitions of $[\theta, 1]$ of mesh less than $\delta$, and the same partition of $[0, \theta]$, the difference between the two values between the modulus signs in (\ref{eq 2.4.3}), has modulus less than $2\ve$. In this difference the two values of $I$ cancel, and so does that part of the sum from the partition of $[0, \theta]$, leaving sums over partitions of $[\theta, 1]$ alone, the modulus of such difference being less than $2\ve$. Hence the Riemann integral over $[\theta, 1]$ exists, see Theorem \ref{2.2.4}. Thus the first part of Theorem aref{2.4.2} is proved.

To prove that if $f$ satisfies the Riemann condition for $f$ then $F$ is dominated integrable, we note that (\ref{eq 2.4.4}) is given, subject to the conditions on $(t_j)$. \nproof

\section{Gauge Integration}\label{s3}
\subsection{Basic Ideas}\label{s3.1}
In Section \ref{s2.4}, (\ref{eq 2.4.5}) shows that in $|P|<\delta$, if the constant $\delta$ is replaced by
\begin{equation}
\label{eq 3.1.1}
\delta(x) = c(x-a),\;\;\;\;\;(a<x\leq b), \;\;\;\;\;\delta(a)>0,
\end{equation}
for some constant $c>0$, we have the \emph{dominated integral} of Osgood and Shisha (1976b, 1977) (see Lewis and Shisha (1983)). Thus when $f$ is unbounded near to $a$ we replace
\begin{equation}
\label{eq 3.1.2}
\int_{a+}^b f(x)\,dx \equiv \lim_{c \rightarrow a+}(R)\int_c^b f\,dx
\end{equation}
by a single limit. Similarly when trouble occurs on the left of $b$ in $[a,b]$, we use
\begin{equation}
\label{eq 3.1.3}
\delta(x) = c(b-x),\;\;\;\;\;(a\leq x< b), \;\;\;\;\;\delta(b)>0,
\end{equation}
for some constant $c>0$, to replace the double limit
\begin{equation}
\label{eq 3.1.4}
\int_{a}^{b-} f(x)\,dx \equiv \lim_{c \rightarrow b-}(R)\int_a^c f\,dx
\end{equation}
by a single limit.

There could be many such difficulties in $[a,b]$, making a suitable $\delta$ very elaborate. To deal with such a case we look at the recovery of a function $F$ (the \emph{primitive}) from its known derivative $f$. Such an $F$ is not unique since for any constant $c$ the derivative of $F+c$ is also $f$. As $F$ is otherwise unique (Theorem ???), $c$ follows from the value of $F$ at a single point.   

Looking at the definition of a derivative, if at a point $x$ we have an error less than $\ve>0$ in the computation of the derivative $f$ of $F$, then we have the inequality
\begin{equation}
\label{eq 3.1.5}
\left|\frac{F(x+h)-F(x)}h - f(x)\right| < \ve.
\end{equation}
For $h$ positive and negative and tending to $0$, but never $0$, the fraction tends to $f(x)$, so that for some range of values of $h$, say $0<|h| < \delta$, the inequality is true. For some points $x$, sometimes large values of $\delta>0$ will suffice, while at other points
 $\delta$ has to be very small. So we take $\delta$ as a function of $x$, as in (\ref{eq 3.1.1}, \ref{eq 3.1.3}). Naturally $\delta$ also depends on $\ve>0$, but we do not write $\delta(x, \ve)$ usually, since $\delta$ and $\ve$ are normally not far apart.\footnote{One could call $\delta$ and $\ve$ the ``heavenly twins of analysis''!}

Rewriting (\ref{eq 3.1.5}) with $h=u-x$,
\begin{equation}
\label{eq 3.1.6}
\left|\frac{F(u)-F(x)}{u-x} - f(x)(u-x)\right| < \ve|u-x|, \;\;\;\;\;(0<|x-u| < \delta(x).
\end{equation}
Thus $F(u) - F(x)$ and $f(x)(u-x)$ differ by a value with modulus less than $\ve|u-x|$. If (\ref{eq 3.1.6}) holds at all points $x$ of the interval $[a,b]$, with $u$ also in $[a,b]$, and if we can fit together a finite number of such intervals $[u,x]$, $[x,u]$ without overlapping, to cover $[a,b]$ exactly, we can build up $F(b)-F(a)$ by splitting it into terms $F(x)-F(u)$, $F(u)-F(x)$, replacing them by $f(x)(x-u)$, $f(x)(u-x)$ respectively, and have a ``Riemann sum'' of such terms. The error is at most the sum of the separate terms in (\ref{eq 3.1.6}), namely $\ve(b-a)$. 

These remarks lead  us to gauge integration, to which we now turn.

\begin{theorem}
\label{3.1.1}
Given real numbers $a<b$ and an arbitrary function $\delta(x)>0$ at each point $x$ of $[a,b]$, there are a finite number of points $a=u_0<u_1< \cdots <u_n=b$, called a \emph{partition} of $[a,b]$, and a finite number of points $x_1, \ldots , xn$, such that $x_j=u_j$ or $x_j=u_{j-1}$, and that $[u_{j-1},u_j]$ lies in \[(x_j - \delta(x_j),\;\;\;
 x_j + \delta(x_j)) \;\;\;\;\;\;(j=1,2, \ldots , n).\] Such an arrangement of the interval-point pairs $(u_{j-1}, u_j, x_j)$ is called a $\delta$-\emph{fine division} $D$ of $[a,b]$ based on $P$.
\end{theorem}
\proof
As $a,b$ are finite, by the Heine-Borel-Lebesgue covering theorem, Theorem \ref{2.1.8}, a finite number of open intervals 
\[
I(x)\equiv (x-\delta(x),\;x+\delta(x))
\]
covers $[a,b]$, i.e.~each point of $[a,b]$ lies in one or more of the finite number of intervals $I(x)$. Some points might be covered by three or more of the $I(x)$, causing difficulties in finding a division, so we first remedy the situation, taking the largest of the $I(x)$ with equal $x$. We can arrange that 
\begin{equation}
\label{eq 3.1.7}
\mbox{each point of $[a,b]$ lies in at most two of the $I(x)$.}
\end{equation}
First suppose that with $a \leq u<v<w\leq b$, $I(u), I(v), I(w)$ have a common point. When $v-\delta(v) \leq u-\delta(u)$,
\[
\delta(v) \geq \delta(u) + v-u > \delta(u),\;\;\;\;\;\;
v+\delta(v) > v+\delta(u) > u+\delta(u),
\]
$I(u)$ is contained in $I(v)$, and omitting $I(u)$ will not remove the cover. Similarly, if $v+\delta(v) \geq w+\delta(w)$ then $I(w)$ is contained in $I(v)$ and can be omitted. Otherwise, as all three intervals have a common point, if
\[
u-\delta(u) < v-\delta(v)<v+\delta(v) < w+ \delta(w),
\]
the union of $I(u)$ and $I(w)$ contains $I(v)$, which can be omitted without breaching the cover. Putting the intervals in a finite sequence and dealing in turn with each interval which meets two other intervals, we have (\ref{eq 3.1.7}), ending with the centres satisfying
\[
a \leq x_1 < x_2< \cdots <x_n \leq b,
\]
the $I(x_j)$ being consecutive, with $a$ in $I(x_1)$, $b$ in $I(x_n)$, and $I(x_j)$ overlapping with $I(x_{j+1})$, the intersection containing points of $(x_j, x_{j+1})$. All three intervals being open, their intersection is an open interval in which we choose a point $y_j$ ($j=1,2, \ldots , n-1$). The points
\[
a \leq x_1<y_1<x_2<y_2< \cdots y_{n-1} <x_n \leq b
\]
divide $[a,b]$ into non-empty $\delta$-fine closed intervals with the points $x_r$ the centres of the open intervals of the cover, giving a $\delta$-fine division of $[a,b]$, as required. \nproof

By Theorem \ref{3.1.1} we can now define the gauge (Riemann-complete, generalized Riemann or Kurzweil-Henstock) integral $F(a,b)$ of the finite-valued function $f$ in $[a,b]$ ($a<b$, finite).

A number $F(a,b)$ is the \emph{gauge integral} of $f$ on $[a,b]$ when, for each $\ve>0$ there is a function $\delta>0$ (the \emph{gauge}) on $[a,b]$ duch that for every $\delta$-fine division $D$ of $[a,b]$,
\begin{equation}
\label{eq 3.1.8}
\left| (D)\sum f(x)(v-u) - F(a,b)\right| < \ve.
\end{equation}
Here, $(D)\sum$ denotes the sum over the $D$ with general term $([u,v],x)$. By \ref{eq 3.1.6}), if $f$ is the finite derivative of a function $F$ on $(a,b)$, with a finite right-hand derivative at $a$ and a finite left-hand derivative at $b$, we have gauge integration of $f$ with integral $F(a,b)=F(b)-F(a)$, thus solving the problem of integrating finite derivatives on $[a,b]$. $f$ can be real-valued or complex-valued. It does not cause any problem that (\ref{eq 3.1.6}) finishes with an error term over $D$ of $\ve(b-a$, for we can replace the original $\ve>0$ by $\ve/(b-a) >0$, so that we end with the error $\ve>0$.

The integration of finite derivatives is the first example of the power of the gauge integral. As Theorem \ref{3.1.1} is the central theorem in gauge integration, we examine it further.

\begin{theorem}
\label{3.1.2}
The existence of $\delta$-fine divisions of $[a,b]$ for each function $\delta>0$ on $[a,b]$, is equivalent to the Heine-Borel-Lebesgue covering theorem for $[a,b]$.
\end{theorem}
\proof
By Theorem \ref{3.1.1} we need only show that the existence implies the covering theorem. Let $[a,b]$ be covered by a collection $C$ of open intervals with centres $x$ in $[a,b]$. For each such $x$ let $(x-k, x+k)$ be one of the intervals in $C$ and let the function $\delta>0$ on $[a,b]$ be $k$ at $x$. If $x$ in $[a,b]$ is not the centre of an interval of $C$, then as $C$ covers $[a,b]$, $x$ lies in an interval of $C$, and so is in an open interval $(x-k,x+k)$ lying in that interval. Here we take $\delta(x)$ to be this $k$, so defining $\delta$ in $[a,b]$. By the existence of $\da$-fine divisions of $[a,b]$ there is such a division consisting of $(I_j,x_j)$, say, and for $J_j\equiv (x_j - \da(x_j), x_j+\da(x_j))$, each $J_j$ is a $K_j$ of $C$ or lies in an interval $K_j$ of $C$. By the $\da$-fineness the $K_j$ form a cover of $[a,b]$ and they form a finite subset of $C$. Hence the covering theorem. \nproof

The kind of proof of Theorem \ref{2.1.8} (Borel) proves Theorem \ref{3.1.1} easily.

\noindent
\textbf{Second proof of Theorem \ref{3.1.1}.}
Let $t$ in $a<t\leq b$ be such that $[a,t]$ has a $\da$-fine division, and let $s$ be the supremum of such $t$. Then for $a<t<a+\da(a)$, $[a,t]$ has $([a,t],a)$ as a $\da$-fine division, to which we add $([t,s],s)$ if $t<s$, and $([s,w],s)$ for some $w$ in $s<w<s+\da(s)$ if $s<b$. In this case $[a,w]$ has a $\da$-fine division with $w>s$, a contradiction implying $s=b$, and then the $\da$-fine division of $[a,b]$. \nproof

Given two gauges $\da_1(x)>0, \da_2(x)>0$ on $[a,b]$, we often replace them by a third gauge $\da_3(x)>0$ so that every $\da_3$-fine division has the properties of a $\da_j$-fine division ($j=1,2$). Thus
\[
\da_3(x) \leq \da_1(x),\;\;\;\;\;\;\da_3(x) \leq \da_2(x)
\]
and we can use $\;\;\;\;\;\;\da_3\equiv \min(\da_1, \da_2),\;\;\;\;\;$
$\da_3(x) = \min(\da_1(x), \da_2(x))>0$.

The first use of this simple device is to show the integral unique.

\begin{theorem}
\label{3.1.3}
Let $f$ be finite on $[a,b]$ with $a<b$. If $G$ and $H$ are two values of the gauge
integral of $f$ on $[a,b]$ then $G=H$.
\end{theorem}
\proof
In (\ref{eq 3.1.8}) let the gauges for $F(a,b)=G,H$, be $\da_1,\da_2$, with $\da_3=\min(\da_1,\da_2)$. Then
\begin{eqnarray*}
|G-H| &=&\left|\left((D)\sum f(x)(v-u) -H\right)- \left((D)\sum f(x)(v-u) -G\right)\right| \vt
&\leq &\left|(D)\sum f(x)(v-u) -H\right|+ \left|(D)\sum f(x)(v-u) -G)\right|\vt
&<& 2\ve 
\end{eqnarray*}
for a $\da_3$-fine division $D$ of $[a,b]$ as both results are true. By Theorem \ref{2.1.3} and arbitrary $\ve>0$, $|G-H|=0$, $G=H$, and $F(a,b)$ is unique whenever it exists. \nproof

\begin{theorem}
\label{3.1.4}
For real numbers $a<b$ let $f,g$ be defined on $[a,b]$, with constants $\alpha, \beta$. If $f,g$ have respective gauge integrals $F,G$ on $[a,b]$ then $\alpha f+ \beta g$ has gauge integral $\alpha F+\beta G$ there.
\end{theorem}
\proof
As in (\ref{eq 3.1.8}) let the gauges for $f,g$ be $\da_1>0, \da_2>0$ with $\da_3=\min (\da_1,\da_2)$. Then
$\left|\left((D)\sum \left(\alpha f(x)+ \beta g(x)\right)(v-u) -(\alpha F + \beta G\right)\right| $ =
\[
=\left|\alpha\left((D)\sum f(x)(v-u) -F\right)+ \beta \left((D)\sum f(x)(v-u) -G\right)\right|\leq  (|\alpha|+ |\beta|)\ve ,
\]
for a $\da_3$-fine division of $[a,b]$. Hence the result by choice of $\ve>0$, $\da_1>0$, $\da_2>0$. \nproof

\begin{theorem}
\label{3.1.5}
In Theorem \ref{3.1.4} for $f,g$ real-valued and $f\leq g$ everywhere in $[a,b]$, then $F\leq G$.
\end{theorem}
\proof
In Theorem \ref{3.1.4}, for every $\ve>0$ and a suitable $D$ depending on $\ve$,
\[
F-\ve<(D)\sum f(x)(v-u) \leq (D)\sum g(x)(v-u) < G+\ve,\;\;\;\;F-G < 2\ve.
\]
By Theorem \ref{2.1.3} we have $F-G \leq 0$, $F\leq G$. \nproof

When $f$ is complex-valued we have a partial converse of Theorem \ref{3.1.4}.

\begin{theorem}
\label{3.1.6}
For real numbers $a<b$,
and $f$ complex-valued and integrable over $[a,b]$, then so are its real and imaginary parts integrable, and conversely.
\end{theorem}
\proof
Let $f=g+\iota h$ be integrable to $G+\iota H$ over $[a,b]$ with $g(x)$, $h(x)$, $G$, $H$ real-valued. Then
$\;\;\;\;\;\;|(D) \sum g(x)(v-u) -G| \;\;\leq $
\begin{eqnarray*}
 &\leq &\left|\left((D)\sum g(x)(v-u) -G\right)+\iota \left((D)\sum h(x)(v-u) -H\right)\right| \vt
& = &\left|(D)\sum f(x)(v-u) -(G + \iota H\right|,
\end{eqnarray*}
and similarly for $h,H$; and $g,h$ are integrable to $G,H$ respectively. The converse follows from Theorem \ref{3.1.4} with $\alpha =1$, $\beta=\iota$.  \nproof

Corresponding to fundamental sequences, we have similarities here.

\begin{theorem}
\label{3.1.7}
For real numbers $a<b$,
and given a finite function $f$ on $[a,b]$ and, for each $\ve>0$, a gauge $\da>0$ depending on $\ve$, such that every two $\da$-fine divisions $D,D'$ of $[a,b]$ satisfy
\begin{equation}
\label{eq 3.1.9}
\left| (D)\sum f(x)(v-u) - (D')\sum f(x')(v'-u')\right| < \ve.
\end{equation}
then $f$ is integrable over $[a,b]$.
\end{theorem}
\proof
Taking $\ve=2^{-j}$ let $\da$ be $\da^*_j>0$, $\da_j = \min(\da_1^*, \da_2^*, \ldots , \da_j^*)>0$.
Then a $\da_j$-fine division of $[a,b]$ is a $\da_j^*$-fine division and a $\da_k$-fine division if $k<j$. For some $\da_j$-fine division $D_j$ of $[a,b]$ let
\[
S_j \equiv (D_j)\sum f(x) (v-u).
\]
Then in (\ref{eq 3.1.9}) we can take
\[
\ve = \frac 1{2^j},\;\;\;\;\;\;D=D_j,\;\;\;\;\;D'=D_k\;\;\;(k>j),\;\;\;\;\;|s_j -s_k| < \frac 1{2^j},
\]
and $(S_J)$ is a fundamental sequence, so convergent to some limit, say $S$. As $k \rightarrow \infty$,
\begin{equation}
\label{eq 3.1.10}
\left|S_j -S\right| < \frac 1{2^j}.
\end{equation}
In (\ref{eq 3.1.9}) we now take $D'=D_j$ and $D$ an arbitrary $\da_j$-fine division of $[a,b]$, and (\ref{eq 3.1.10});
\[
\left|(D)\sum f(x)(v-u) -S\right| < \frac 1{2^{j-1}}.
\]
Given $\ve>0$, there is always an integer $j$ with $2^{1-j} \leq \ve$, and so $f$ is integrable over $[a,b]$. \nproof

\begin{theorem}
\label{3.1.8}
If $[p,q] \subseteq [a,b]$ ($a<b$) with $f$ integrable over $[a,b]$, then $f$ is integrable over $[p,q]$ to $F(q)-F(p)$ where $F(p) \equiv \int_a^p f\,dx$ ($p>a$), $F(a)=0$.
\end{theorem}
\proof
$[a,b] \setminus [p,q] = [a,p)\cup (q,b]]$, where $[a,p)$ or $(q,b]$ or both can be empty, and there is no division over an empty set. Given $\ve>0$, let $\da>0$ on $[a,b]$ be such that
\begin{equation}
\label{eq 3.1.11}
\left| (D)\sum f(x)(v-u) - F(b)\right| < \ve
\end{equation}
for all $\da$-fine divisions $D$ of $[a,b]$. For $D_1,D_2,D_3,D_4$ $\da$-fine divisions of $[a,p]$, $[q,b]$, $[p,q]$, $[p,q]$, respectively (taking $D_j$ empty if over an empty set), then $D_1 \cup D_2 \cup D_3$ and $D_1 \cup D_2 \cup D_4$ are $\da$-fine divisions of $[a,b]$, so that by (\ref{eq 3.1.1}),
\begin{eqnarray}
\label{eq 3.1.12}
&&\left| (D_3)\sum f(x)(v-u) - (D_4)\sum f(x')(v'-u')\right| \nonumber \vt
&=&
\left| (D_1 \cup D_2 \cup D_3)\sum f(x)(v-u) -F(b) \right. \;\;-  
\nonumber \vt
&&\;\;\;\;\left. +\;\; F(b) - (D_1 \cup D_2 \cup D_4)\sum f(x')(v'-u'))\right| \nonumber \vt
&\leq & 2\ve.
\end{eqnarray}
Hence by Theorem \ref{3.1.7}, $f$ is integrable over $[p,q]$, and in (\ref{eq 3.1.12}) we can let 
\begin{eqnarray}
\label{eq 3.1.13}
&&(D_4) \sum f(x) (v-u) \;\;\rightarrow \;\;\int_p^q f\,dx, \nonumber \vt
&&
\left| ( D_3)\sum f(x)(v-u) - \int_p^q f\,dx\right| \leq 2\ve.
\end{eqnarray}
This last is the same for all $[p,q] \subseteq [a,b]$ and all $\da$-fine divisions $D_3$ of $[p,q]$, a uniformity. Hence if $a \leq p<q<r \leq b$, as a $\da$-fine division $D_3$ of $[p,q]$ and a $\da$-fine division $D_5$ of $[q,r]$
give a $\da$-fine division $D_3\cup D_5$ of $[p,r]$ while the sums add, 
\begin{eqnarray*}
&& \left| \int_p^r f\,dx - \int_p^q f\,dx - \int_q^r f\,dx \right|\;\;\;= \vt
&=&
\left| \int_p^r f\,dx 
- (D_3 \cup D_5)\sum-  \sum \int_p^q f\,dx + (D_3) \sum - \int_q^r f\,dx +(D_5) \sum  \right| \vt
& \leq & 6 \ve.
\end{eqnarray*}
As the values of the integrals do not depend on $\ve>0$, the first modulus is $0$ and 
$
\int_p^r f\,dx  = \int_p^q f\,dx  + \int_q^r f\,dx $.
\nproof

We say that the integral is \emph{finitely additive}. For $p=a$, the last integral is $F(r) - F(q)$. Clearly we have a vital theorem, going inwards to intervals inside $[a,b]$. Going outwards,
\begin{theorem}
\label{3.1.9}
If $a<b<c$ with $f$ integrable over $[a,b]$ and $[b,c]$, then $f$ is integrable over $[a,c]$ to the sum of the other two integrals.
\end{theorem}
\proof
Let us take $\da^*(x), \da^*_1(x), \da^*_2(x)$ such that 
$\da^*(x) = \frac 12 |x-b|$ ($x \neq b$), $\da^8(b)=1$. Then
\begin{eqnarray*}
 \left|(D_1 )\sum f(x)(v-u)- F_1   \right| 
& < & \frac 12 \ve,\vt
 \left|(D_2 )\sum f(x)(v-u)- F_2   \right| 
& < & \frac 12 \ve,\vt
F_1 &\equiv & \int_a^b f\,dx,\;\;\;\;\;\;
F_2\;\;\;\equiv \;\;\; \int_b^c f\,dx,
\end{eqnarray*}
for every $\da_1^*$-fine division $D_1$ of $[a,b]$
and every $\da_2^*$-fine division $D_2$ of $[b,c]$.
We take
\[
0<\da(x) \leq \left\{
\begin{array}{ll}
\min(\da^*(x), \da_1^*(x)) & (a \leq x<b), \vt
\min(\da^*(b), \da_1^*(b), \da_2^*(b)) & (x=b, \vt
\min(\da^*(x), \da_2^*(x)) & (b < x\leq c).
\end{array}
\right.
\]
By Theorem \ref{3.1.1} there is a $\da$-fine division $D$ of $[a,c]$. By construction of $\da^*$, the intervals $[u,v]$ with $b$ in $[u,v]$, and $([u,v],x)$ in $D$, have $x=b=u$ or $v$, as no other interval can reach $b$. If $v=b$ then $[u,b] \subseteq [a,b]$ and $([u,b],b)$ is $\da_1^*$-fine.
If $u=b$ then $[b,v] \subseteq [b,c]$ and $([b,v],b)$ is $\da_2^*$-fine.
 Thus $D$ can be split up into a $\da_1^*$-fine $D_1$ over $[a,b]$ and a $\da_2^*$-fine $D_2$ over $[b,c]$ , and
\begin{eqnarray*}
&& \left|(D )\sum f(x)(v-u)- (F_1 +F_2)  \right| \;\;\;=\vt
 &=&
 \left|(D_1 )\sum f(x)(v-u)- F_1 +(D_2 )\sum f(x)(v-u)-F_2)  \right| \;\;\;<\;\;\; \ve.
 \end{eqnarray*}
Being true for every $\ve>0$ and every suitable $\da>0$, $f$ is integrable over $[a,c]$ to $F_1+F_2$. \nproof

\begin{theorem}
\label{3.1.10}
Let intervals $I_j \equiv [u_{j-1},u_j]$ ($j=1,2, \ldots ,n$), where $a=u_0 <b=u_n$, form a partition $P$ of $[a,b]$. Then there is a gauge $\da$ on $[a,b]$ such that every $\da$-fine division $D$ of $[a,b]$ is a \emph{refinement} of $P$, i.e.~if $(I,x)$ is in $D$ then $I \subseteq I_j$ for some $j$ in $1, 2, \ldots , n$.
\end{theorem}
\proof
We take
\[
\begin{array}{l}
0<2\da(x) \leq \left\{
\begin{array}{l}
\min(u_j-x, x-u_{j-1}) \;\;\; (u_{j-1} < x< u_j,\;j=1,2, \ldots , n), \vt
\min(u_j - u_{j-1}, u_{j+1}-u_j) \;\;\; (x=u_j,\; j=2,3, \ldots , n-1),
\end{array}
\right. \vt
0 <2\da(a)  \leq \;\;u_1-a, \vt
0< 2\da(b)  \leq \;\;b-u_{n-1}.
\end{array}
\]
By Theorem \ref{3.1.1} there is a $\da$-fine division $D$ of $[a,b]$ formed of interval-point pairs $([u,v],x)$ with $x=u$ or $x=v$. If $u_{j-1}<x<u_j$ then $[u,v]$ cannot include $u_{j-1}$ nor $u_j$, by the inequalities, and $[u,v] \subseteq I_j$. If $x=u_j$ then $\da(x)$ is smaller than the distances from $u_j$ to $u_{j-1}$ (if any) and to $u_{j+1}$ (if any), so that $[u,v]$ lies in $I_j$ or in $I_{j+1}$. Hence the theorem. \nproof

People call the next theorem the Henstock lemma or the Saks-Henstock lemma. We reserve the second name for (\ref{eq 3.1.15}) and call (\ref{eq 3.1.16}) the Kolmogorov-Saks-Henstock lemma.

\begin{theorem}
\label{3.1.11}
Let $f$ integrable over $[a,b]$ ($a<b$) and,given $\ve>0$, let $\da(x)>0$ on $[a,b]$ be such that, for all $\da$-fine divisions $D$ of $[a,b]$,
$[b,c]$, then $f$ is integrable over $[a,c]$ to the sum of the other two integrals.
\begin{equation}
\label{eq 3.1.14}
\left| (D)\sum f(x)(v-u) - F(a,b)\right| < \ve,\;\;\;\;\;\;
F(u,v) \equiv \int_u^v f\,dx.
\end{equation}
If $D_1 \subseteq D$ then
\begin{equation}
\label{eq 3.1.15}
\left| (D_1)\sum f(x)(v-u) - F(u,v)\right| \leq \ve,
\end{equation}
\begin{equation}
\label{eq 3.1.16}
 (D)\sum \left|f(x)(v-u) - F(u,v)\right| \leq 4 \ve.
\end{equation}
\end{theorem}
\proof
Theorem \ref{3.1.8} shows that for $F(a,a) =0$, $F(a,x) = \int_a^x f\,dx$, $F(u,v)$ exists equal to $F(a,v) - F(a,u)$, and we can write (\ref{eq 3.1.14}) in the form
\begin{equation}
\label{eq 3.1.17}
 \left|(D)\sum \left(f(x)(v-u) - F(u,v)\right)\right| <\ve.
\end{equation}
For those $([u,v],x)$ in $D \setminus D_1$ we can replace each $([u,v],x)$ by a $\da$-fine division $D^*_u$ of $[u,v]$ and still have a $\da$-fine division of $[a,b]$ to replace $D$ in (\ref{eq 3.1.14}). Letting each such
\[
(D^*_u)\sum f(x')(v'-u') \rightarrow F(u,v),
\]
we replace each such $f(x)(v-u)-F(u,v)$ in (\ref{eq 3.1.17}) by $0$ and have \ref{eq 3.1.15}). In the original $D$ with $f$ real-valued we take $D_1$ the terms in $D$ with $f(x)(v-u)-F(u,v)\geq 0$, and $D\setminus D_1$ the rest. Then (\ref{eq 3.1.15}) twice gives (\ref{eq 3.1.16}) with $2\ve$ on the right. If $f$ is complex-valued, by Theorem \ref{3.1.6} the result holds with the real part of $f$ and with the imaginary part of $f$, giving $4\ve$ on the right of (\ref{eq 3.1.16}). 
\nproof

The latter result shows how near the integral is to some of its Riemann sums, in a rather deeper way than (\ref{eq 3.1.14}), and it leads to the idea of \emph{variation}.

Divisions are finite sets of interval-point pairs $(I,x)$. Let $h(I,x)$ be a function of such interval-point pairs; for example,
\begin{equation}
\label{eq 3.1.18}
h(I,x) = f(x)mI - F(I) \;\;\;\;\;\;(I=[u,v],\;\;\;mI=v-u).
\end{equation}
Given a gauge $\da$ on $[a,b]$ let
\[
V(h;\da;[a,b]) = \sup_D (D) \sum |h(I,x)|,
\]
the supremum of the sums over $\da$-fine divisions $D$ of $[a,b]$. The \emph{variation} $V(h;[a,b])$ of $h$ in $[a,b]$ is defined to be
\[
V(h;[a,b]) \equiv \inf_\da V(h;\da;[a,b])
=\limsup_{\da \rightarrow 0+} (D) \sum |h(I,x)|,
\]
the infimum being taken over all gauges $\da$ on $[a,b]$. If $V(h;[a,b])$ is finite, we say that $h$ is of \emph{bounded variation} over $[a,b]$. If $V(h;[a,b])=0$ we say that $h$ is of \emph{variation zero} over $[a,b]$.

\textbf{Thus (\ref{eq 3.1.16}) states that, for the given conditions, (\ref{eq 3.1.18}) is of variation zero over $[a,b]$.}

Two interval-point functions $h(I,x)$, $k(I,x)$ are \emph{variationally equivalent} if $h-k$ is of variation zero. Thus, in (\ref{eq 3.1.16}), $f(x)mI$ and $F(I)$ are variationally equivalent.

If $X$ is a set on the real line with indicator $\chi(X;x)$ (equal to $1$ if $x$ is in $X$, and to $0$ if $x$ is not in $X$) we write, respectively, $V(h;\da;[a,b];X)$ and $V(h;[a,b];X)$ for
\[
V(h.\chi(X;\cdot);\da;[a,b]),\;\;\;\;\;\;\;\;\;
V(h.\chi(X;\cdot);[a,b]).
\]
Do not confuse $V(h;[a,b];X)$ with $V(h;\da;X)$. These have $[a,b]$ in different places, and the second $V$ has the Greek letter $\da$.

If $V(h;[a,b];X)$ is finite we say that $h$ is of \emph{bounded variation in} $X$, relative to $[a,b]$,
while if $V(h;[a,b];X)=0$ we say that $h$ is of \emph{variation $0$ in $X$, relative to $[a,b]$}. We say that a property is true \emph{$h$-almost everywhere ($h$-a.e.)} if it is true except in a set $X$ with $V(h;[a,b];X)=0$ Such a set $X$ is said to be of \emph{$H$-variation zero}. Sometimes p.p.~(\emph{presque partout}) replaces a.e. If $h([u,v], x)=v-u$ we usually omit the $h$- from $h$-variation, $h$-a.e.~etc. 

\begin{theorem}
\label{3.1.12}
If $h,k$ are variationally equivalent in $[a,b]$, then for each set $X$ on the real line, even if one side is $+\infty$,
\begin{equation}
\label{eq 3.1.19}
V(h:[a,b];X)=V(k;[a,b];X).
\end{equation}
If $V(h;[a,b])=0$ and if $[u,v]\subseteq [a,b]$, then $V(h;[u,v])=0$. If $V(h;[a,b];X)$ $=0$ and if $Y \subseteq X$, then
\begin{equation}
\label{eq 3.1.20}
V(h;[a,b];Y)=0.
\end{equation}
Let $f(x)$ be integrable over $[a,b]$ with integral $F(I)$ ($I \subseteq [a,b]$). Then $|f(x)|$  is integrable over  $[a,b]$
\begin{equation}\label{eq 3.1.21}
 \mbox{ if and only if } F \mbox{ is of bounded variation over }[a,b].
\end{equation}
\end{theorem}
\proof
In (\ref{eq 3.1.19}), $h-k$ has variation zero. Hence $|h-k|\chi(X;x)$ has variation zero since the sum over a division $D$, of $|h-k|\chi(X;x)$, has the same or fewer terms than the sum for $|h-k|$. Using the inequalities
\[
|h| \leq |h-k| + |k|,\;\;\;\;\;\;|k| \leq |h-k| + |h|,
\]
and given $\ve>0$, there is a gauge $\da$ on $[a,b]$ with 
\begin{eqnarray*}
V(h-k;\da;[a,b];X)&<&\ve, \vt
V(h;[a,b];X) & \leq & V(h;\da;[a,b];X) \vt
&\leq & V(h-k;\da;[a,b];X) + V(k;\da;[a,b];X) \vt
&<& \ve \;\;+\;\;V(k;\da;[a,b];X) , \vt
V(h;[a,b];X) &\leq & V(h;[a,b];X) 
\end{eqnarray*}   
as $\da \rightarrow 0+$ and so $\ve \rightarrow 0+$, even when one side is $+\infty$. Interchanging $h$ and $k$ gives equality. In (\ref{eq 3.1.20}) we have the same or fewer terms in sums over $[u,v]$ as for over $[a,b]$, since the closure of $[a,b] \setminus [u,v]$ is covered by $\da$-fine divisions; and similarly for $Y$ replacing $X$,
\[
V(h; \da; [u,v]) \leq V(h; \da; [a,b]),\;\;\;\;\;\;
V(h; \da; [a,b];Y) \leq V(h; \da; [a,b];X).
\]
In (\ref{eq 3.1.21}) let $|f|$ be integrable to $K$ over $[a,b]$. Then, given $\ve>0$, there is a gauge $\da$ on $[a,b]$ such that for all $\da$-fine divisions $D$ of $[a,b]$ with $h(I,x) =f(x)mI$,
\[
\begin{array}{rlcll}
K-\ve &<& (D)\sum |f(x)| (v-u) &<&K+\ve,\vt
K-\ve &<& V(h;\da;[a,b] )&\leq & K+\ve,\vt
K-\ve &\leq & V(h;[a,b] )&\leq & K+\ve,
\end{array}
\]
so $\int_a^b |f|\,dx = K = V(h;[a,b])$, and $K$ is finite.

By (\ref{eq 3.1.19}), $V(F;[a,b]) = V(h;[a,b])$, finite. Conversely, if $f$ is integrable to $F$ of bounded variation over $[a,b]$, (\ref{eq 3.1.16}) implies that for any $\da$-fine divisions $D$ of $[a,b]$, $(D)\sum |F(I)|$ tends to a finite limit in $[a,b]$, to the same limit. But $F$ is finitely additive, so that if $D'$ is a division of $[u,v]$ then
\[
|F([u,v])| = \left|(D')\sum F(I)\right| \leq (D') \sum |F(I)|
\]
and sums of $|F|$ rise or stay the same under subdivision. As $F$ is of bounded variation in $[a,b]$, 
sums of $|F|$ over $[a,b]$ are bounded by $V(F;[a,b])$. If $\ve>0$ there is a division $D$ with
\begin{equation}
\label{eq 3.1.22}
V(F;[a,b])-\ve <(D)\sum |F(I)| \leq V(F;[a,b]).
\end{equation}
By Theorem \ref{3.1.10}, for some gauge $\da$ in $[a,b]$, every $\da$-fine division of $[a,b]$ is a refinement of this $D$, and so its sum of $|F|$ also lies within the same bounds as in (\ref{eq 3.1.22}). As $\ve>0$ is arbitrary, $V(F;[a,b])$ is the limit of sums $F$ as $\da \rightarrow 0+$, and $|f|$ is integrable to $V(F;[a,b])$ over $[a,b]$.  
\nproof

\begin{theorem}
\label{3.1.13}
If $h$ is an interval-point function on $[a,b]$ and if $(X_j)$ is a sequence of sets on the real line with union $X$, then
\begin{equation}
\label{eq 3.1.23}
V(h;[a,b];X) \leq \sum_{j=1}^\infty V(h;[a,b];X_j).
\end{equation}
If $V(h;[a,b];X)=0$ with $f$ a finite-valued point function on $[a,b]$, then $Vfh;[a,b];X)=0$. Conversely, if $V(fh;[a,b];X)=0$ and if $X_0$ is the set where $f\neq 0$, then
\begin{equation}
\label{eq 3.1.24}
V(h;[a,b];X\cap X_0)=0.
\end{equation}
\end{theorem}
\proof
In (\ref{eq 3.1.23}), if some $X_j$ have points in common, and if we remove the common points from the $X_j$ with the larger $j$, obtaining disjoint $X_j^*$, the union $X$ is the union of the $X_j^*$, and we need only prove (\ref{eq 3.1.23}) for the $X_j^*$. Replacing the $X_j^*$ by the $X_j$ will not lower the sum, so that we can assume the $X_j$ disjoint. Further, we take the tight-hand side finite or there is nothing to prove.

Given $\ve>0$, let the gauge $\da_j$ on $[a,b]$ be such that
\begin{equation}
\label{eq 3.1.25}
V(h;\da_j;[a,b];X_j)\leq V(h;[a,b];X_j) + \frac \ve{2^j}\;\;\;\;\;\;(j=1,2, \ldots).
\end{equation}
We take the gauge $\da$ to satisfy
\begin{equation}
\label{eq 3.1.26}
\da(x) = \left\{ \begin{array}{ll}
\da_j(x) & (x \in X_j,\;\;\;j=1,2, \ldots),\vt
1& (x \in \setminus X).
\end{array}\right.
\end{equation}
If $D$ is a $\da$-fine division of $[a,b]$ and if $D_j$ is the subset of the $(I,x)$ in $D$ with $x$ in $X_j$, then as $D$ is only a finite set there is an integer $m$ with $D_j$ empty for $j>m$. By (\ref{eq 3.1.25}), 
\begin{eqnarray*}
(D)\sum |h(I,x)|\chi(X;x) &=& \sum_{j=1}^m (D_j) \sum |h(I,x)| \vt
&\leq &\sum_{j=1}^m V(h;\da_j;[a,b];X_j) \vt
&<& \sum_{j=1}^\infty V(h;[a,b];X_j) + \ve, \vt
V(h;\da;[a,b];X) &\leq & \sum_{j=1}^\infty V(h;[a,b];X_j) + \ve.
\end{eqnarray*}
As $\da \rightarrow 0+$ and $\ve \rightarrow 0+$ we have (\ref{eq 3.1.23}). Thus if each $X_j$ has $h$-variation zero, so has $X$.

If $X$ has $h$-variation zero and if $f$ is a finite-valued point function let $X_j$ be the set of all $x$ in $X$ with
\[
j-1 < |f(x)| <j\;\;\;\;\;(j=1,2, \ldots).
\]
Then by (\ref{eq 3.1.20}), $X_j$ has $h$-variation zero and
\begin{eqnarray*}
V(fh;[a,b];X_j) & \leq & jV(h;[a,b];X_j) \;\;\;=\;\;\;0,\vt
V(fh;[a,b];X) & \leq & \sum_{j=1}^\infty jV(h;[a,b];X_j) \;\;\;=\;\;\;0
\end{eqnarray*}
Conversely, $f^{-1}$ exists in $X_0$, so that if $V(fh;[a,b];X)=0$ then 
\[
V\left(h;[a,b];X \cap X_0\right) =
V\left(f^{-1}(fh);[a,b];X \cap X_0\right)=0
\]
by (\ref{eq 3.1.20}) and the first part.  \nproof

\begin{theorem}
\label{3.1.14}
The integral $F([u,v])$ of $fh$ over $[u,v]$, is $0$ for all $[u,v]\subseteq [a,b]$ ($b>a$), or for all $[u,v] \subseteq [a,b]$, or for $[a,b]$ alone when $fh \geq 0$, if and only if $f=0$ in $[a,b]$ except for a set $X$ in $[a,b]$ with $V(h;[a,b];X)=0$.
\end{theorem}
\proof
If $F([a,v])=0$ for all $v$ in $a<v\leq b$ then
\[
F([u,v]) = F([a,v) - F([a,u]) =0\;\;\;\;(a<u<v\leq b).
\]
By (\ref{eq 3.1.16}), given $\ve>0$, there is a gauge $\da$ such that for all $\da$-fine divisions $D$ of $[a,b]$,
\[
(D) \sum |f(x) h(I,x)| \leq 4\ve,\;\;\;\;V(fh;\da;[a,b]) \leq 4\ve,\;\;\;\;V(fh;[a,b])=0.
\]
(\ref{eq 3.1.16}) finishes the proof. 

Conversely, $V(fh;[a,b])=0$, so that for $\ve>0$ a gauge $\da>0$ exists such that for all $\da$-fine divisions $D$ of $[u,v]\subseteq [a,b]$,
\[
\left|(D)\sum fh -0\right| \leq (D)\sum |fh| \leq \ve,
\]
and $F([u,v])$ exists and is $0$. In case $fh \geq 0$ with $F([a,b])=0$, then for $\ve>0$ there is a gauge $\da$ such that for all $\da$-fine divisions $D$ of $[a,b]$,
\[
(D) \sum |fh| = \left|(D)\sum fh -0\right| \leq \ve,\;\;\;\;\;\;V(fh;[a,b]) =0,
\]
and the rest of the proof follows. \nproof

\begin{theorem}
\label{3.1.15}
If $[u,v]\subseteq [a,b]$ and $h(I,x)$ is defined over $[a,b]$, then
\begin{equation}
\label{eq 3.1.27}
V(h;[u,v]) \leq  V(h;[a,b]).
\end{equation}
For $a<b<c$ with $h(I,x)$ defined over $[a,c]$, then
\begin{equation}
\label{eq 3.1.28}
 V(h;[a,b]) +  V(h;[b,c])=  V(h;[a,c]).
\end{equation}
\end{theorem}
\proof
As $[a,b]\setminus [u,v]$ is one or two intervals we take divisions over their closure. Adding to a division of $[u,v]$ we have a division of $[a,b]$. Sums over the various divisions not of $[a,b]$ add together to give a sum over $[a,b]$, which is not greater than $V(h;\da;[a,b])$ for $\da$-fine divisions. As $\da \rightarrow 0$ we prove (\ref{eq 3.1.27}). Thus in (\ref{eq 3.1.28}), if either of $V(h;[a,b])$, $V(h;[b,c])$ is $+\infty$ the result is true. So we can assume both variations finite, and, given $\ve>0$, and each $\da_1^*.)$, we take some $\da_1^*$-fine division $D_1$ over $[a,b]$,
\[
(D_1)\sum |h(I,x)| > V(h;[a,b]) -\ve,
\]
and similarly for $[b,c]$ and $\da_2^*$ and $D_2$. As in the proof of Theorem \ref{3.1.9} we construct $\da$ over $[s,b]$ and can then replace $\da_1^*, \da_2^*$ by $\da$. For some $\da$-fine $D_1,D_2$ we have
\begin{eqnarray*}
(D_1 \cup D_2)\sum |h| &>& V(h;[a,b]) +V(h;[b,c]]) \;\;-\;\;2\ve,\vt
V(h;[a,c]) &\geq & V(h;[a,b]) + V(h;[b,c]).
\end{eqnarray*}
To prove the opposite inequality and so (\ref{eq 3.1.28}), given $\ve>0$ and $\da_0>0$, there is always a $\da_0$-fine
division $D$ of $[a,b]$ with
\begin{equation}
\label{eq 3.1.29}
(D)\sum |h| >
V(h;[a,c]) -\ve.
\end{equation}
As in Theorem \ref{3.1.9} we construct $\da$ from $\da_0 = \da_1^*$ in $[a,b]$, and $\da_0 = \da_2^*$ in $[b,c]$, with a $\da$-fine division $D$ satisfying (\ref{eq 3.1.29}), which by the construction can be split up into a $\da$-fine division $D_1$ of $[a,b]$ and a $\da$-fine division $D_2$ of $[b,c]$, with
\[
V(h;\da;[a,b]) +V(h;\da;[b,c]]) \geq (D_1)\sum |h|+(D_2)\sum |h| >
V(h;[a,c])-\ve.
\]
As $\ve \rightarrow 0$ we prove the opposite inequality and so (\ref{eq 3.1.28}).  \nproof

Thus the variation is finitely additive. It is also an outer measure in Lebesgue theory (Theorem \ref{3.1.13}), the next property needing a \emph{regular} outer measure, though not needed here.

\begin{theorem}
\label{3.1.16}
Let $(X_n)$ be a monotone increasing sequence of sets in $[a,b]$ with $h(I,x)$ defined there. Then
\begin{equation}
\label{eq 3.1.30}
\lim_{n\rightarrow \infty}V(h;[a,b];X_n) =  V(h;[a,b]; \lim_{n\rightarrow \infty}X_n).
\end{equation}
\end{theorem}
\proof
As $X\equiv \lim_{n \rightarrow \infty}X_n$ contains $X_j$,
\begin{equation}
\label{eq 3.1.31}
V(h;[a,b];X_j) \leq   V(h;[a,b];X),
\;\;\;\;\; 
\lim_{n\rightarrow \infty}V(h;[a,b];X_n) \leq V(h;[a,b];X).
\end{equation}
Thus if the limit is a conventional $+\infty$, the result is true. If the limit is finite we prove the opposite inequality to (\ref{eq 3.1.31}). For let the gauges $\da_n(x), \da(x)$ in $[a,b]$ satisfy
\begin{eqnarray}
\label{eq 3.1.32}
V(h;\da_n;[a,b];X_n)& \leq &  V(h;[a,b];X_n) +\frac \ve{2^n}\;\;\;(\ve>0, \;\;n=1,2,3, \ldots ) ,\end{eqnarray}
\begin{eqnarray}
\da(x) &=& \left\{
\begin{array}{ll}
\da_n(x)&(x \in X_n \setminus X_{n-1},\;\;n=1,2, \ldots ,\;\; X_0 \mbox{ empty}), \nonumber \vt
 1 &( x \in \setminus X).\nonumber
 \end{array} \right. \nonumber
\end{eqnarray}
If $D$ is a $\da$-fine division of $[a,b]$ and if $D', D_n'$ are the subsets of $D$ with the $x$ in $X$ and in $X_n \setminus X_{n-1}$, respectively ($n=1,2,3, \ldots$), there is a greatest integer $m$ depending on $D$, such that $D_m'$ is not empty. (Remember that $D$ contains only a finite number of $(I,x)$.)
Let $E_j$ be the union of the $[u,v]$ with $([u,v],x) \in D_j'$
and let $F_j$ be the unio of the $[u,v]$ with $([u,v],x) \in D\setminus D_j'$. Then by Theorem \ref{3.1.15}, (\ref{eq 3.1.28}),
\begin{eqnarray*}
V(h;\da;E_j;X_j) + V(h;\da;F_j;X_j) &\leq & 
V(h;\da;[a,b];X_j), \vt
V(h;\da;E_j;X_j) &\leq & V(h;[a,b];X_j) +\frac \ve{2^j} -V(h;F_j;X_j) \vt
&=& V(h;E_j;X_j) +\frac \ve{2^j},\vt
(D) \sum |h| \chi(X;x) &=& (D') \sum |h| 
\;\;\;=\;\;\;\sum_{j=1}^m (D')\sum |h|\vt 
&\leq &\sum_{j=1}^m  V(h;\da;E_j;X_j) \vt
&<& \sum_{j=1}^m  \left(V(h;E_j;X_j) + \frac \ve{2^j} \right) \vt
&<& \sum_{j=1}^m V(h;E_j;X_m) +\ve \vt &\leq &V(h;[a,b];X_m) + \ve,
\end{eqnarray*}
\[
V(h;[a,b];X) \leq V(h;\da;[a,b];X) \leq 
\lim_{j\rightarrow \infty} V(h;[a,b];X_j) +\ve,
\]
giving the opposite inequality to (\ref{eq 3.1.31}) as $\ve \rightarrow 0$, and hence (\ref{eq 3.1.30}). \nproof

\begin{example}
\label{ex 3.1.1}
If $f=g$ almost everywhere in $[a,b]$ ($b>a$) and the integral $F$ of $f$ exists over $[a,b]$, then $g$ is integrable to $F$ over $[a,b]$. (\emph{\textbf{Hint:} Consider the integral of $f-g$ and use Theorem \ref{3.1.14}.})
\end{example}

\begin{example}
\label{ex 3.1.2}
If $F([u,v]) = \int_u^v f\,dx$ then $\int_a^b g\,dF = \int_a^b gf\,dx$ provided one side exists.
\end{example}
To ensure that the integral using $dF$ is of the type defined in (\ref{eq 3.1.8}) we assume that $F([a,x])$ is strictly increasing in $x$, taking all values between its minimum and maximum in $[a,b]$> Or we can define the integral of $h(I,x)$ by replacing $f(x)(v-u)$ in (\ref{eq 3.1.8}) by $h([u,v],x)$. This gives the \emph{integration by substitution} of the calculus.

The second proof of Theorem \ref{2.1.1} is in Henstock (1955) p.~277, with extra reference to rationals. Gauge integration (unnamed) first appeared there, followed independently by Kurzweil (1957), Lemma 1.1.1 p.~423. Here the proof is the first proof of Theorem \ref{3.1.1}. A proof for rectangles was found by J.~Mawhin in P.~Cousin (1895) p.???, but with no reference to integration.

Much of this section's theory first appeared in Henstock (1963), with a 25-year update in Henstock (1988).

\subsection{The Cauchy and Harnack Extensions, and Integration over Infinite Intervals}
\label{s3.2}
If the singularities of a function form a complicated set, constructions of the type given in section \ref{s2.4} become very intricate. However, if one extends the gauge integral similarly, the resulting integral is still the gauge integral; we have not found anything new. 

This is shown in Theorem \ref{3.2.1} for the case of a singularity at $c$, the gauge integral existing over $[a,b]$ for all $b$ in $a<b<c$. For a singularity at $a$, with the function integrable over $[b,c]$ for all $b$ in $a<b<c$, there is a similar theorem with the integral over $b,c]$ tending to a limit as $b \rightarrow a+$. If the singularity occurs between $a$ and $b$ we split the interval at the singularity and use both results. 

Effectively, we are illustrating the continuity of the integral with respect to its endpoints of the range of $x$.

So far, the gauge integral has been defined over finite intervals $[a,b]$. We need a slightly different definition for an interval of infinite length, and a proof like that of Theorem \ref{3.2.1} is needed to connect with Cauchy's definition. 
Such results can be called \emph{Cauchy extensions}, remembering Cauchy (1823), Lecs.~24, 25. In his very detailed construction of the integral of a derivative, Denjoy used repeatedly the Cauchy extensions,together with an extension of Harnack, and this section finishes with a proof that the Harnack extension also is contained in the gauge integral.

\begin{theorem}
\label{3.2.1}
Let $f$ be finite over $[a,c]$ for real numbers $a<c$, and let the gauge integral $F(b)$ of $F$ over $[a,b]$ exist for all $b$ in $a<b<c$. If $F(b) \rightarrow G$ as $b \rightarrow c-$ then $f$ is gauge integrable over $[a,c]$ with value $G$. Conversely, if $F(c)$ exists, so does $F(b)$ for $a<b<c$, and $F(b) \rightarrow F(c)$ as $b \rightarrow c-$.
\end{theorem}
\proof
Let $b_j=c-(c-a)2^{-j}$ ($j=0,1,2, \ldots$). We are given that for $j>0$, $F(b_j)$ exists, so that by Theorem \ref{3.1.8} the integral $F(b_j)-F(b_{j-1})$ exists over $B_j \equiv [b_{j-1}, b_j]$. Hence given $\ve>0$, there is a gauge $\da_j(x)$ in $B_j$ for which every $\da_j$-fine division $D_j$ over $B_j$ satisfies
\[
\left| (D_j) \sum f(x)(v-u) -F(b_j) + F(b_{j-1})\right| < \frac \ve{2^{j+1}}.
\]
If $b_{J-1} <w<b_j$ with $D^*$ a $\da_J$-fine division of $[b_{J-1},w]$, then by Theorem 3.1.11 (\ref{eq 3.1.15}), the Saks-Henstock lemma,
\begin{eqnarray*}
\left| (D^*)\sum f(x)(v-u) -F(w) + F(b_{J-1})\right| &\leq & \ve{2^{J+1}}, \vt
\left|F(x) -G\right| &<& \frac \ve{2}\;\;\;\;\;(c-\da(c)<x<c)
\end{eqnarray*}
for small $\da(c)>0$. Also let
\begin{eqnarray*}
\da(x) &=& \da_j(x)\;\;\;\;\;(b_{j-1}<x<b_j), \vt
\da(b_j) &=& \min\left(\da_j(b_j),\;\;\da_{j+1}(b_j)\right)\;\;\;\;\;(j=1,2, \ldots), \vt
\da(c)&=& \da_1(b_0).
\end{eqnarray*}
By this construction we have in $[a,c]$ divisions $D_1, \ldots , D_{j-1}, D^*$, and interval-point pair $([w,c], c)$, forming a 
$\da$-fine division $D$ of $[a,c]$ for which
\begin{equation}
\label{eq 3.2.1}
\left|(D)\sum f(x)(v-u) - G\right| < \ve.
\end{equation}
To ensure that every $\da^*$-fine division of $[a,c]$ is a division of type $D$ we only have to arrange that the division refines the partition 
\[
B_1, \ldots , B_{J-1}, [b_{J-1}, w], [w,c],
\]
by Theorem \ref{3.1.10}, with $\da^* \leq \da$. Thus by (\ref{eq 3.2.1}), $f$ is gauge integrable over $[a,c]$ to $G$.

For the converse we use Theorems \ref{3.1.8}, \ref{3.1.11} (???), (\ref{eq 3.1.15}) with $D_1$ consisting of a single interval-point pair $([w,c],c)$.  \nproof

Turning to the definition of a division over $x\geq a$, for some real number $a$, the infinite integral being written $[a, +\infty)$, we begin by taking some real number $b>a$ and $\da$-fine divisions of $[a,b]$. The controller of the $\da$-fineness of divisions of $[a,b]$ is a neighbourhood of each $x$ in $[a,b]$. 
\[
([u,v],x)\mbox{ is } \da\mbox{-fine if }
[u,v] \subseteq (x-\da(x),\;x+ \da(x)),\;\;\;\;x=u\mbox{ or } x=v.
\]
This hints at how we may define $\da$-fineness over $[a, +\infty)$, using a suitable neighbourhood of the conventional $+\infty$, remembering that no real number is greater than $+\infty$. 

For some real number $b>a$ we use a one-sided open interval with $x>b$, written $(b, +\infty)$, $[B,+\infty)$ being $\da$-\emph{fine } if $B>b$. The rest of the division of $[a, \infty)$ is a $\da$-fine division of $[a,B]$. 

The \emph{Riemann sum} for the $\da$-fine division of $[a, + \infty)$ and for a point function $f$ over $[a, +\infty)$, is obtained by taking the conventional $f(+\infty)$ and $f(+\infty)(+\infty -B)$ to be $0$, to give a continuous integral ``at $+\infty$'', and then adding the Riemann sum over $[a,B]$. Using these Riemann sums over finite intervals, we can have a direct definition of the integral over $[a,+infty)$.

The interval of $x$ satisfying $x \leq b$, for some real number $b$, written $(-\infty, b]$, has a similar definition of the integral over the interval. One simple way is to change every $x$ to $-x$. Similarly for $-\infty, +\infty)$, the whole real line. 

We could put together the definitions for $[0,+\infty)$ and $(-\infty, 0]$, and so for each division we would use real numbers $u<v$ and $(-\infty, u)$, $(v, +\infty)$. 

The $u,v$ are independent. If, for example, we took $u=-v$ we would have an integral useful in complex variable theory, but not one amenable to the approach of this section as one could arrange an infinity on one side of $0$ that is cancelled by a negative infinity on the other side, so that the integral over $(-\infty, u]$ need not exist.

The integrals of this section, except for the integrals over $(-\infty, +infty)$ with dependence of $u$ on $v$, have a theory like that of section \ref{s3.1}, together with analogues of Theorem \ref{3.2.1}, the last analogues showing that the integrals are the same as those defined by the Cauchy extensions, but without needing double limits.

These results can be used as examples to test the comprehension of the reader.

\section{Limits Under the Integral Sign} \label{s4}
\subsection{Introduction and Necessary and Sufficient Conditions} \label{s4.1}
Early in the history of the calculus it was found necessary to take limits under the integral sign, for example, when can we integrate an infinite series term by term?
\begin{equation}
\label{eq 4.1.1}
\int\left(\sum_{n=1}^\infty a_n(x) \right)dx = \sum_{n=1}^\infty \left(\int a_n(x)\,dx \right),\;\;\;\int \lim_{n \rightarrow \infty} s_n(x)\,dx = \lim_{n \rightarrow \infty} \int s_n(x) \,dx,
\end{equation}
where $s_n(x)$ is the sum to $n$ terms of the series of the $a_n(x)$. Two limits are involved, the limit that gives the integral, and the limit as $n \rightarrow \infty$ of the sum $s_n(x)$ to $n$ terms, and the two limiting operations might interact and produce strange results.

One often quoted example involves the sequence $(r_n)$ of rationals (see Theorem \ref{2.1.6}). Let $s_n(x)$ be the indicator of the first $n$ rationals, the function that is $0$ except when $x=r_j$ ($1 \leq j \leq n$), when we get $s_n(r_j) =1$. he Riemann integral of $s_n(x)$ over any interval is $0$. But $\lim_{n\rightarrow \infty}s_n(x)$ is the indicator of all the rationals, and so cannot be Riemann integrable as its lower Darboux sum is $0$ (there are irrationals in every interval) and its upper Darboux sum is the length of the interval (there are rationals in every interval).

When Lebesgue published his magnificent paper in 1902, that particular difficulty was resolved. The Lebesgue integrals of this $s_n(x)$ and its limit are both $0$.

The problem of calculus integration of derivatives was studied very carefully, and after centuries of effort Denjoy (1912) gave a most complicated construction to integrate all derivatives. A variety of other definitions followed, culminating in the gauge integral, the main subject of this book.

In evaluating some integrals, differentiation under the integral sign is a great help. Here, for an integrable function $s(x,y)$ of $x$,
\begin{eqnarray}
\label{eq 4.1.2}
\frac d{dy} \int s(x,y)\,dx & \equiv & \lim_{h \rightarrow \infty} \int \frac{s(x, y+h) -s(x,y)}h dx \nonumber \vt
&=& \int
\lim_{h \rightarrow \infty}  \frac{s(x, y+h) -s(x,y)}h dx
\nonumber \vt
&\equiv & \int \frac{\partial s(x,y)}{\partial y} dx.
\end{eqnarray}
Again we have the inversion of two limit operations.
De la Vall\'{e}e Poussin (1892) quotes a request of C.~Jordan, which is translated from the French as follows: \emph{Give a rigorous theory of differentiation under the integral sign of definite integrals, with precise conditions which limit Leibnitz's rule, principally for unbounded regions of integration or unbounded functions, and particularly many celebrated definite integrals.}

Clearly C.~Jordan wished to study (\ref{eq 4.1.2}) in all its generality, and after a century, gauge integration enables us to find necessary and sufficient conditions for (\ref{eq 4.1.1}) and (\ref{eq 4.1.2}).
  
An even more complicated inversion is that of the order of repeated integrals, see Tonelli (1924). For intervals $I$ of $x$, $J$ of $y$,
\begin{eqnarray}
\label{eq 4.1.3}
\int_J \left(\int_I s(x;y) ,dx\right) dy & \equiv & \lim_{\da \rightarrow 0}{}^* 
(D) \sum \left( \int_I s(x;y) \,dx\right) (\beta - \alpha)
\nonumber \vt
&=& \int_I \left( \lim_{\da \rightarrow 0} {}^*
(D) \sum s(x;y) (\beta-\alpha) \right) dx
\nonumber \vt
&= & \int_I
\left( \int_J s(x;y)\, dy\right)dx,
\end{eqnarray}
where here, $\lim^*_{\da \rightarrow 0}$ denotes the gauge integral limit relative to $y$. Two conditions together, are both necessary and sufficient. If now we interchange the two processes of integration, the two conditions are interchanged, as will be seen later.

We begin with a sequence $(s_n(x))$ of real- or complex-valued functions integrable on a real interval $[a,b]$, that tend \emph{pointwise} (i.e.~at each point) to a function $f(x)$ on $[a,b]$ as $n \rightarrow \infty$.

This pointwise convergence is not a great limitation. It can be relaxed to become $s_n(x) \rightarrow f(x)$ a.e., that is, the (finite) limit occurs everywhere except in a set of variation zero (or measure zero). Here the indicator $\chi$ of the set where convergence does not occur, is such that $f(x)\chi(x)$ has integral $0$ for every finite-valued function $f$ (see Theorem \ref{3.1.14}) and in particular $s_n(x)$ and $s_n(x)(1-\chi(x))$ have the same value for their integrals, and $s_n(x)(1-\chi(x))$ converges everywhere in $[a,b]$.

Two properties are examined,
\begin{eqnarray}
&& \mbox{the integrability of the limit function }f\mbox{ over }[a,b], \label{eq 4.1.4} \vt
&&\lim_{n\rightarrow \infty}\int_a^b s_n(x)\,dx = \int_a^b f(x)\,dx \equiv \int_a^b \lim_{n\rightarrow \infty} s_n(x)\,dx, \label{eq 4.1.5}
\end{eqnarray}
and the two conditions result from the two properties.

In the statement and proof of Theorem \ref{3.1.11} with $f$ real- or complex-valued, (\ref{eq 3.1.15}) and (\ref{eq 3.1.16}) follow from (\ref{eq 3.1.14}). Similarly from (\ref{eq 4.1.6}) (later) there follow stronger results and a new kind of variation.

\begin{theorem}
\label{4.1.1}
For each $n$ let the real- or complex-valued $s_n(x)$ be integrable over $[a,b]$, and let $s_n(x) \rightarrow f(x)$ as $n \rightarrow \infty$, for all $x$ in $[a,b]$. Then a necessary and sufficient condition for (\ref{eq 4.1.4}) is that there are a closed interval $C$ (when all $s_n(x)$ are real-valued) of arbitrarily small length, or a closed circle $C$ on the complex plane of arbitrarily small radius, and positive functions $N(x)$, $\da(x)$ on $[a,b]$, such that for all positive integer-valued functions $n(x)\geq N(x)$ on $[a,b]$ and all $\da$-fine divisions $D$ of $[a,b]$,
\begin{equation}
\label{eq 4.1.6}
(D)\sum s_{n(x)}(x) (v-u) \in C.
\end{equation}
More exactly,there is a function $F$ on $[a,b]$ such that
\begin{equation}
\label{eq 4.1.7}
\left|(D)\sum s_{n(x)}(x) (v-u) - F(b) + F(a) \right| < \ve
\end{equation}
when $C$ has diameter $\ve>0$.
For $D'$ a subset of $D$,
\begin{equation}
\label{eq 4.1.8}
\left|(D')\sum \left(s_{n(x)}(x) (v-u) - F(v) + F(u)\right) \right| \leq \ve,
\end{equation}
the Saks-Henstock lemma here. The Kolmogorov-Henstock lemma here is 
\begin{equation}
\label{eq 4.1.9}
(D)\sum \left|s_{n(x)}(x) (v-u) - F(v) + F(u) \right| \leq K \ve,
\end{equation}
where $K=2$ for real-valued functions and $K=4$ for complex-valued functions.
\end{theorem}
\proof
For all $x$ in $[a,b]$, $s_n(x) \rightarrow f(x)$ as $n \rightarrow \infty$, and for an $N(x)>0$,
\begin{equation}
\label{eq 4.1.10}
 \left|s_{n(x)}(x) -f(x) \right| < \ve \;\;\;
 \mbox{ for all } n(x) \geq N(x) \mbox { in } [a,b].
\end{equation}
If $f$ is integrable to $F$ on $[a,b]$, there is a gauge $\da(x)>0$ on $[a,b]$ such that
\begin{equation}
\label{eq 4.1.11}
\left|(D)\sum f(x) (v-u) - F \right| \leq  \ve
\end{equation}
for all $\da$-fine divisions $D$ of $[a,b]$. By (\ref{eq 4.1.10}), (\ref{eq 4.1.11}) we have (\ref{eq 4.1.6}) from
\begin{equation}
\label{eq 4.1.12}
\left|(D)\sum s_{n(x)}(x) (v-u) - F \right| < \ve + \ve (b-a)\;\;\;\;\;(n(x) \geq N(x)).
\end{equation}
As the integral $F(x)$ of $f$ exists over $[a,x]$ ($a<x\leq b$), the integral of $f$ over $[u,v] \subseteq [a,b]$ being $F(v)-F(u)$, taking $F(a)=0$, then $F=F(b)-F(a)$ and (\ref{eq 4.1.12}) becomes (\ref{eq 4.1.7}).

For $D'$ a subset of the set $D$ of $([u,v],x)$ and $E$ the union of the finite number of $[u,v]$ with $([u,v],x)$ in $D \setminus D'$, $f$ is integrable over the disjoint intervals of $E$ and, given $\ve'>o$, we have (\ref{eq 4.1.10}), (\ref{eq 4.1.11}) for $\ve'$ replacing $\ve$, $E$ replacing $[a,b]$, and $D''$ replacing $D$, where $D''$ is a $\da$-fine division of $E$ and $D' \cup D''$ is a $\da$-fine division of $[a,b]$ to replace $D$. So
\begin{eqnarray*}
\left|(D'')\sum s_{n(x)}(x) (v-u) - G \right|& < & \ve' + \ve'(b-a)\vt &&(n(x) \geq N(x),\;n(x) \geq N_1(x)), \vt
\left|(D'\cup D'')\sum s_{n(x)}(x) (v-u) - F \right| &<&  \ve + \ve(b-a)\;\;\;\;\;(n(x) \geq N(x)), \vt
\left|(D')\sum s_{n(x)}(x) (v-u) - F(v) +F(u) \right|& < &\ve +\ve' + (\ve+\ve')(b-a)\,
\end{eqnarray*}
where $G$ is the integral of $f$ over $E$. As $D'$ is independent of $\ve'$ we can let $\ve'\rightarrow 0+$ and obtain (\ref{eq 4.1.8}). Having chosen $D$ and $n(x)$, the inequality in (\ref{eq 4.1.8}) stays the same, whatever subset $D'$ we choose. If all $s_n$ are real we choose for $D'$ all $([u,v],x)$ of $D$ with
\[
\begin{array}{l}
s_{n(x)}(v-u) -F(v)+F(u) \;\;\;\geq \;\;\;0,  \vt
\left|(D')\sum \left(s_{n(x)}(x) (v-u) - F(v)+F(u)\right) \right| \;\;\;=\;\;\; \vt
 \;\;\;\;\;\;\;\;=\;\;\;(D')\sum \left|s_{n(x)}(x) (v-u) - F(v) +F(u) \right|.
\end{array}
\]
In $D \setminus D'$ the differences are negative, and again the modulus of the sum is the sum of the moduli. Adding the two results we have (\ref{eq 4.1.9}) with $K=2$. When some $s_n$ are complex-valued we deal separately with the real parts and with the imaginary parts, and so obtain (\ref{eq 4.1.9}) with $K=4$. Of course $\ve$ has been replaced by $\ve + \ve(b-a)$; to obtain (\ref{eq 4.1.9}) exactly we need only begin $\ve/(1+b-1)>0$ instead of $\ve>0$. 

To show that (\ref{eq 4.1.6}) is \emph{sufficient} (i.e.~gives (\ref{eq 4.1.4}) when each $s_n(x)$ is integrable over $[a,b]$, with $s_n(x) \rightarrow f(x)$ as $n \rightarrow \infty$ for all $x$ in $[a,b]$, we let $n(x) \rightarrow \infty$ at each $x$, so that
\begin{equation}
\label{eq 4.1.13}
(D)\sum f(x)(v-u) \in C
\end{equation}
as $C$ is closed, for each $\da$-fine division $D$ of $[a.b]$. For some $C$ with diameter $\ve>0$ there
is a gauge $\da$ for which (\ref{eq 4.1.13}) is true. So for $D^*$ another $\da$-fine division of $[a,b]$,
\[
\left| (D)\sum f(x)(v-u) - (D^*)\sum f(x^*)(v^*-u^*) \right| \leq 2\ve
\]
and Theorem \ref{3.1.7} shows the integrability of $f$ over $[a,b]$. Similarly (\ref{eq 4.1.7}) and (\ref{eq 4.1.9}) are sufficient. \nproof

The later Example \ref{ex 4.1.1} shows that if 
the small diameter condition on $C$ is omitted, the weakened condition is not sufficient.

\begin{theorem}
\label{4.1.2}
Given (\ref{eq 4.1.6}) and Theorem \ref{4.1.1}, a necessary and sufficient condition for (\ref{eq 4.1.5}) is that for each $\ve>0$ there are a positive integer $N$ and a gauge $\da_n>0$ ($n=1,2, \ldots $) on $[a,b]$ such that for all $n\geq N$ and all $\da_n$-fine divisions $D$ of $[a,b]$,
\begin{equation}
\label{eq 4.1.14}
\left| (D)\sum s_n(x)(v-u) - F \right| < \ve,\;\;\;\;\;\;F \equiv \int_a^b f\,dx,
\end{equation}
giving an interval or closed circle $C_1$ with centre $F$ and length $2\ve$ or radius $\ve$, respectively.
\end{theorem}
\proof
Given the integrals $S_n$ of $s_n(x)$ and $F$ of $f(x)$, over $[a,b]$, we assume $S_n \rightarrow F$ ($n \rightarrow \infty$). For $\ve>0$ let $N$, $\da_n>0$ on $[a,b]$ satisfy
\[
\left|S_n-F\right| < \frac 12 \ve\;\;\;\;(n \geq N),\;\;\;\;\;\;\;\;
\left| (D_n)\sum s_n(x)(v-u) - S_n \right| < \frac 12\ve,
\]
for every $\da_n$-fine division $D_n$ of $[a,b]$. These two inequalities give (\ref{eq 4.1.14}). Conversely we can let $\da_n(x) \rightarrow 0$ suitably to obtain (\ref{eq 4.1.5}) from (\ref{eq 4.1.14}). \nproof

In (\ref{eq 4.1.14}) $n$ is constant over $[a,b]$, so that for any subset $D'$ of $D$, $n$ in $D \setminus D'$ is fixed and cannot tend to $+\infty$, as it seems that a Saks-Henstock lemma here would need, so the lemma seems unattainable here by an easy proof.

We now generalize Theorems \ref{4.1.1}, \ref{4.1.2} to the more general limit processes already mentioned. The first step is to change the integer $n$ to $y$ in $y \geq 0 $ with $y \rightarrow \infty$.

\begin{theorem}
\label{4.1.3}
For all $x$ in $[a,b]$ and all $y \geq 0$ let $s(x;y)$ be real- or complex-valued with $s(x;y) \rightarrow f(x)$ (finite) as $y \rightarrow \infty$, for all $x$ in $[a,b]$. For each $y\geq 0$ let $s(x;y)$ be integrable with respect to $x$ in $a,b]$. Then a necessary and sufficient condition for the integrability of $f$ over $[a,b]$ is that there are a closed interval $C_2$ of arbitrarily small length (when all $s(x;y)$ are real-valued) or a closed circle $C_2$ on the complex plane of arbitrarily small radius, and positive functions $Y$, and $\da$ on $[a,b]$, such that for all functions $y(x) \geq Y(x)$ on $[a,b]$ and all $\da$-fine divisions $D$ of $[a,b]$,
\begin{equation}
\label{eq 4.1.15}
(D)\sum s(x;y(x))(v-u) \in C_2.
\end{equation}
More exactly, there is a function $F$ on $[a,b]$ such that
\begin{equation}
\label{eq 4.1.16}
\left| (D)\sum s(x;y(x))(v-u) - F(b) + F(a)\right| < \ve \;\;\;\;\;\;(y(x) \geq Y(x))
\end{equation}
when $C_2$ has diameter $2\ve >0$. For $D'$ a subset of $D$,
\begin{equation}
\label{eq 4.1.17}
\left| (D')\sum \left(s(x;y(x))(v-u) - F(v) + F(u)\right)\right| \leq \ve ,
\end{equation}
\begin{equation}
\label{eq 4.1.18}
 (D)\sum \left|s(x;y(x))(v-u) - F(v) + F(u)\right| \leq K\ve ,
\end{equation}
where $K=2$ for real-valued functions, and $K=4$ otherwise.
\end{theorem}
The last two results are the Saks-Henstock and Kolmogorov-Henstock lemmas here. The proof of the theorem follows the proof of Theorem \ref{4.1.1}.
\begin{theorem}
\label{4.1.4}
Given $s(x;y)$ integrable in $x$ on $[a,b]$ for each $y\geq 0$, and $s(x;y) \rightarrow f(x)$ (finite) as $y \rightarrow \infty$, and (\ref{eq 4.1.15}), a necessary and sufficient condition for the integral of $s(x;y)$ to tend to the integral of $f$, both integrals over $[a,b]$, is that for some closed interval $C_3$ of arbitrarily small length, or some closed circle $C_3$ on the complex plane of arbitrarily small radius, with $C_3 \cap C_2$ not empty, a real number $Y \geq 0$, a gauge $\da(x;y)$ in $x$ on $[a,b]$ for each $y \geq Y$, and all $\da(x;y)$-fine divisions $D$ of $[a,b]$,
\begin{equation}
\label{eq 4.1.19}
(D) \sum s(x;y) (v-u) \in C_3\;\;\;\;\;\;(y \geq Y).
\end{equation}
\end{theorem}
For proof see the proof of Theorem \ref{4.1.2}; no Saks-Henstock not Kolmogorov-Henstock lemmas seem possible here.

We now specialize the last two theorems to continuity and differentiability under the integral sign, and we also consider the inversion of order of integration in double integrals. Note that Fubini's theorem needs integration in two dimensions at least.

For the continuity of the integral of $s(x;y)$ as $y$ tends to some number $w$ in an interval $[t,r]$ of $y$ let $s(x;y)$ be integrable in $x$ over $[a,b]$ for each fixed $y$ in $[t,r]$. Then in particular $s(x;w)$ is integrable in $x$ and we can ignore Theorem \ref{4.1.3}, change $f(x)$ to $s(x;w)$, and apply Theorem \ref{4.1.4} to obtain the necessary and sufficient condition for
\[
\lim_{y \rightarrow w}\int_a^b s(x;y)dx 
= \int_a^b \lim_{y \rightarrow w} s(x;y)dx .
\]
If we are given that $s(x;w)$ is integrable \emph{a priori} over $[a,b]$ we need only also invoke Theorem \ref{4.1.3}.

We differentiate $\int_a^b s(x;y)dx$ by $y$ on replacing $s(x;y)$, $f(x)$ by
\[
R(s;x,y,h) \equiv \frac{s(x;y+h) -s(x;y)}h\;\;\;\;\;(h \neq 0),\;\;\;\;\;\;\;\;\;\; \frac{\partial s(x;y)}{\partial y}
\]
respectively. We suppose that
\begin{equation}
\label{eq 4.1.20}
\lim_{h \rightarrow 0} R(s;x,y,h) \equiv \frac{\partial s(x;y)}{\partial y}
\end{equation}
exists (finite) in $a \leq x \leq b$ for every $y$ in $[t,r]$.

\begin{theorem}
\label{4.1.5}
Let the real- or complex-valued $s(x;y)$ be integrable with respect to $x$ in $[a,b]$ for each $y$ in $[t,r]$. If (\ref{eq 4.1.20}) exists for each $y$ in $[t,r]$, one-sided at $t$ and $r$, then a 
necessary and sufficient condition for the integrability with respect to $x$ of $\partial s(x;y)/ \partial y$ over $[a,b]$, is that for some closed interval $C_4$ of arbitrarily small length (when $s$ is real-valued) or some closed circle $C_4$ on the complex plane of arbitrarily small radius, two positive functions $\da$, $N$, all real $h(x)$ in $0<|h(x)| \leq N(x)$, and all $\da$-fine divisions $D$ in $[a,b]$,
\begin{equation}
\label{eq 4.1.21}
(D) \sum R(s;x,y,h(x))(v-u) \in C_4.
\end{equation}
More exactly, there is a function $F$ on $[a,b]$ such that
\begin{equation}
\label{eq 4.1.22}
(D)\left| \sum R(s;x,y,h(x))(v-u) -F(b) + F(a)    \right|< \ve\;\;\;\;\;(0<|h(x)| \leq N(x))
\end{equation}
when $C_4$ has diameter $\ve>0$. For $D'$ a subset of $D$,
\begin{equation}
\label{eq 4.1.23}
(D) \left|\sum \left( R(s;x,y,h(x))(v-u)  -F(v) + F(u)  \right) \right| \leq \ve,
\end{equation}
\begin{equation}
\label{eq 4.1.24}
(D) \sum \left| R(s;x,y,h(x))(v-u)   -F(v) + F(u)   \right|\leq K \ve,
\end{equation}
where $K=2$ for real-valued $s$, and $K=4$ otherwise
\end{theorem}
For proof we change $s(x;y)$ to $R(s;x,y,h)$ in Theorem \ref{4.1.3}.

\begin{theorem}
\label{4.1.6}
Given Theorem \ref{4.1.5} and (\ref{eq 4.1.22}) for each $\ve>0$, the necessary and sufficient condition that
\begin{equation}
\label{eq 4.1.25}
\frac d{dy} \int_a^b s(x;y)dy = \int_a^b \frac{\partial s}{\partial y} dx,
\end{equation}
is that for some closed interval $C_5$ of arbitrarily small length (when $s$ is real-valued), or some closed circle $C_5$ on the complex plane of arbitrarily small radius, with $C_4 \cap C_5$ not empty, a number $N>0$, a positive function $\da(\cdot, h)$ on $[a,b]$ for each $h$ in $0<|h| \leq N$, and all $\da(\cdot, h)$-fine divisions $D$ of $[a,b]$,
\begin{eqnarray}
\label{eq 4.1.26}
&&(D) \sum R(s;x,y,h)(v-u)\;\;=\vt
\nonumber &=& \frac 1h (D) \sum \left(s(x;y+h) - s(x;y)\right)(v-u) \;\;\;\in\;\;\; C_5\;\;\;\;\;(0<|h| \leq N)
\end{eqnarray}
as $h$ is the same for each fraction.
\end{theorem}
For proof, change $s(x;y)$ to $R(s;x,y,h)$ in Theorem \ref{4.1.4}.

In a final limit under the integral sign, Tonelli's theorem in Lebesgue integration extends to gauge integration, with necessary and sufficient conditions.

\begin{theorem}
\label{4.1.7}
The real- or complex-valued $s(x;y)$ is integrable over an interval $I$ with respect to $x$, for each fixed $y$, and over an interval $J$ with respect to $y$, for each fixed $x$. Then
\begin{equation}
\label{eq 4.1.27}
\int_I \left(\int_J s(x;y) \,dy \right) dx \equiv K
\end{equation}
exists if and only if there are a closed interval $C_6$ of arbitrarily small length (for $s$ real-valued) or a closed circle $C_6$ on the complex plane of arbitrarily small radius, some positive function $\da(x)$ on $I$, $\da(y;x)$ on $J$ for each $x$ in $I$, and all $\da$-fine divisions $D$ of $I$, all $\da(y;x)$-fine divisions $D(x)$ of $J$, with $([u,v],x)$ in $D$, $([\alpha, \beta],y)$ in $D(x)$, and
\begin{equation}
\label{eq 4.1.28}
(D)\sum (v-u) \left(\left(D(x)\right) \sum s(x;y) (\beta - \alpha)\right) \in C_6.
\end{equation}
For $C_6$ of diameter $\ve>0$ and $D'$ a subset of $D$, the Saks-Henstock lemma here is
\begin{equation}
\label{eq 4.1.29}
\left|(D')\sum (v-u) \left(\left(D(x)\right) \sum s(x;y) (\beta - \alpha)\right) 
- \int_u^v \left(\int_J s(x;t)\,dt \right) dx \right| \leq \ve,
\end{equation}
while the Weierstrass-Henstock lemma here is
\begin{equation}
\label{eq 4.1.30}
(D)\sum \left|(v-u) \left(\left(D(x)\right) \sum s(x;y) (\beta - \alpha)\right) 
- \int_u^v \left(\int_J s(x;t)\,dt \right) dx \right| \leq L \ve
\end{equation}
with $L=2$ for real-valued $s$, and otherwise $L=4$.
\end{theorem}
\proof
This is a more complicated limit theorem than Theorem \ref{4.1.1} and it deserves a more detailed proof. When (\ref{eq 4.1.27}) exists, given $\ve>0$, there is a gauge $\da$ on $I$ with
\begin{equation}
\label{eq 4.1.31}
\left|(D)\sum \left((v-u) \int_J s(x;y)\,dy\right) -K \right| < \ve.
\end{equation}
Then for $mI$ the length of $I$, suitable $\da(y;x)>0$, $\da(y;x))$-fine $D(x)$ over $J$, and (\ref{eq 4.1.31}), 
\begin{eqnarray*}
v-u &>&0, \vt
(D)\sum (v-u) &=& mI\;\;>\;\;0,\vt
\left| (D(x))\sum (\beta - \alpha) s(x;y) - \int_J s(x;y)\,dy\right| &<& \ve,\vt
\left|(D) \sum \left((v-u) (D(x))\sum s(x;y)(\beta - \alpha)\right)  - K\right| &<& \ve + \ve mI,
\end{eqnarray*}
and (\ref{eq 4.1.28}) is true with an interval or circle centre $K$ and diameter $2\ve + 2\ve mI$.  

To prove (\ref{eq 4.1.28}) sufficient we note that $s(x;y)$ is integrable with respect to $y$ over $J$ and that $C_6$ is closed. Let $\da(y;x) \rightarrow 0$ in order to give
\[
(D) \sum \left((v-u) \int_J s(x;y)\,dy\right) \in C_6
\]
with $C_6$ of arbitrarily small diameter. Using Theorem \ref{3.1.7}, $\int_J s(x;y)\,dy$ is integrable with respect to $x$ over $I$ and (\ref{eq 4.1.27}) exists. Further, (\ref{eq 4.1.29}) and (\ref{eq 4.1.30}) follow as in Theorem \ref{4.1.1}.  \nproof

\begin{theorem}
\label{4.1.8}
Given (\ref{eq 4.1.28}) and the conditions in Theorem \ref{4.1.7}, a necessary and sufficient condition that (\ref{eq 4.1.3}) holds is that for each $\ve>0$ there is a gauge $\da$ on $J$ such that for all $\da$-fine divisions $E$ of $J$, a gauge $\da(x;y)$ on $I$ and depending on $\da$ and $y$, all $\da(x;y)$-fine divisions $D(y)$ of $I$, and a closed interval $C_7$ of length less than $\ve$ (if $s$ is real-valued) or a closed circle $C_7$ of diameter less than $\ve$, such that $C_6 \cap C_7$ is not empty, and for $([\alpha, \beta],y)$ in $E$, $([u,v],x)$ in $D(y)$, we have
\begin{equation}
\label{eq 4.1.32}
(E) \sum \left((\beta - \alpha) (D(y))\sum s(x;y)(v-u) \right) \in C_7.
\end{equation}
For $E'$ a subset of $E$ we have Saks-Henstock and Kolmogorov-Henstock extensions, respectively,
\begin{equation}
\label{eq 4.1.33}
\left|(E') \sum \left(\left((\beta - \alpha) (D(y))\sum s(x;y)(v-u) \right)
- \int_\alpha^\beta \left(\int_I s(x;y)\,dx\right) dy\right) \right| \leq \ve,
\end{equation}
\begin{equation}
\label{eq 4.1.34}
(E) \sum \left|(\beta - \alpha)\left( (D(y))\sum s(x;y)(v-u) \right)
- \int_\alpha^\beta \left(\int_I s(x;y)\,dx\right) dy \right| \leq L\ve,
\end{equation}
with $L=2$ for real-valued $s$, and otherwise $L=4$
\end{theorem}
\proof
To show the necessity of (\ref{eq 4.1.32}) let (\ref{eq 4.1.3}) exist. Then given $\ve>0$, a gauge $\da$ on $J$ is such that for all $\da$-fine divisions $E$ of $J$,
\begin{eqnarray*}
\left| (E) \sum \left((\beta - \alpha) \int_I s(x;y)\,dy \right)- K\right| &<& \ve , \vt
(E) \sum \left((\beta - \alpha) \int_I s(x;y)\,dy\right)  &=& \int_I \left( (E)\sum s(x;y) (\beta - \alpha) \right) dx.
\end{eqnarray*}
Further, for each $y$ in $J$ and all $x$ in $I$, there is a gauge $\da(x;y)$ such that for $\da(x;y)$-fine divisions $D(y)$ of $I$,
\begin{eqnarray*}
\left| (D(y)) \sum s(x;y) (v-u)
- \int_I s(x;y)\,dy \right|
 &<& \ve , \vt
\left| (E) \sum (\beta - \alpha)  (D(y))\sum s(x;y) (v-u)  -K \right| 
&<& \ve + \ve mJ,
\end{eqnarray*}
giving (\ref{eq 4.1.32}) with $C_7 = [K-\ve -\ve mJ,\;K+\ve +\ve mJ]$, having $K$ in common with $C_6$.

To show that (\ref{eq 4.1.32}) is sufficient for (\ref{eq 4.1.3}), given the results of Theorem \ref{4.1.7}, we have $K$ in $C_6$ and $C_6 \cap C_7$ not empty. As $C_6$ is arbitrarily small in diameter and $C_7$ is closed, $K$ is also in $C_7$. In (\ref{eq 4.1.32}) we use gauge integration with $\da(x;y) \rightarrow 0$ for each fixed $y$ and
\[
(E) \sum (\beta - \alpha) \int_I s(x;y)\,dx \in C_7.
\]
Again, Theorem \ref{3.1.7} gives integrability, this time of
\[
\int_J\left(\int_I s(x;y)\,dx \right) dy \in C_7,\;\;\;\;\;\;K \in C_7,
\]
which give (\ref{eq 4.1.3}).  

Results (\ref{eq 4.1.28}) and (\ref{eq 4.1.32}) are very similar, with $x$ and $I$ interchanged with $y$ and $J$, so that Saks-Henstock and Kolmogorov-Henstock lemmas hold here in both theorems. Further, the limit function has been $\int_J s(x;y)\,dy$, taken under the integral sign $\int_I \cdots dx$, for which we have the two conditions (\ref{eq 4.1.28}), (\ref{eq 4.1.32}). 
\nproof

We could have taken the limit function to be $\int_I s(x;y)\, dx$ under the integral sign $\int_J \cdots dy$, in which case the two necessary and sufficient conditions are (\ref{eq 4.1.32}) and (\ref{eq 4.1.28}) in that order, which gives a check of the theory.

When $s(x;y)$ is gauge integrable in $I \times J$, (\ref{eq 4.1.30}) and (\ref{eq 4.1.34}) can be improved to
\begin{equation}
\label{4.1.35}
(D) \sum (D(x)) \sum \left| s(x;y)(v-u)(\beta-\alpha)
- \int_u^v \left(\int_\alpha^\beta s(x;y)\,dy \right) dx\right| \leq L\ve.
\end{equation}
See Theorem \ref{4.2}.
This seems unprovable using only the two one-dimensional integrations.

Bartle (1994--95) assumes the integrability of $s_n$ and $f$ in Theorems 
\ref{4.1.1}, \ref{4.1.2}, with another necessary and sufficient condition.

\begin{theorem}
\label{4.1.9}
For a sequence $(s_n)$ of real-valued functions on a closed interval $I$ of the real line, with $s_n$ and $f$ integrable over $I$, but not with $s_n \rightarrow f$ as $n \rightarrow \infty$, the necessary and sufficient conditions for (\ref{eq 4.1.5}) is that, given $\ve>0$, there is an integer $N$ such that if $n \geq N$, a gauge $\da_n$ on $I$ has the property that for every $\da_n$-fine division $D$ of $I$, 
\begin{equation}
\label{eq 4.1.36}
\left| (D)\sum s_n(x) (v-u) -(D)\sum f(x)(v-u) \right| < \ve.
\end{equation}
\end{theorem}
Condition (\ref{eq 4.1.36}) is like (\ref{eq 4.1.14}), the integral $F$ of $f$ being replaced by the approximating Riemann sum. For proof, follow the proof of Theorem \ref{4.1.2}.

The part of the theory involving closed bounded sets $C$, is taken from Henstock (1988), pp.~105--108, 127--128, 135--136, and 170--172. The various Saks-Henstock and Kolmogorov-Henstock lemmas in this section are new and involve a new kind of variation, which we could call \emph{limit variation} of functions $h_y(I,x)$ of the interval-point pairs $(I,x)$, together with a parameter $y$ that in some sense tends to a limit. When $y$ is nearer than $Y$ to the limit, in this sense, we write $y>>Y$. Thus
\[
\mbox{LV}\left(h_y(I,x);Y,\da;a,b\right) = \sup_{y>>Y,\; \da'\leq \da,\;D}(D)\sum \left|h_y(I,x)\right|
\]
where $D$ is a $\da'$-fine division of $[a,b]$, and the \emph{limit variation} is 
\[
\mbox{LV}\left(h_y(I,x);Y;a,b\right) =
\inf_{Y, \da} \,\mbox{LV}\left(h_y(I,x);Y,\da;a,b\right)
\equiv \limsup (D) \sum |h_y(I,x)|.
\]
In every example of limits up to this point the limit variation is $0$. In (\ref{eq 4.1.9}), (\ref{eq 4.1.18}), (\ref{eq 4.1.24}), (\ref{eq 4.1.30}), (\ref{eq 4.1.34}), respectively, with $I=[u,v]$,
\begin{eqnarray*}
h_y(I,x) &=& s_y(x)(v-u) -F(v)+F(u), \;\;\;\;y=n(x) \rightarrow \infty, \vt
h_y(I,x) &=& s_y(x)(v-u) -F(v)+F(u), \;\;\;\;y \rightarrow \infty,\vt
h_y(I,x) &=& R(s;x,y,k)(v-u) -F(v)+F(u), \;\;\;\;k(x) \rightarrow 0, \vt
h_y(I,x) &=& (v-u)(D(x))\sum s(x;t) (\beta - \alpha) - \int_u^v \left( \int_J s(x;t) \,dt\right) dx,
\end{eqnarray*}
where $y$ and $Y$ are gauges $\da''$ and $\eta$ on $J$, $\da''\leq \eta$, and $D(x)$ is $\da''$-fine, with $I^*=[u,v]$.
\[
h(I^*,y) = (\beta-\alpha) (D(y))\sum s(t;y)(v-u) - \int_\alpha^\beta \left(\int_I s(t;y)\,dt\right)dy,\;\;\;I^*=[\alpha,\beta],
\]
where $\xi$ and $X$ are gauges $\da_1$ and $\zeta$ on $I$, $\da_1 \leq \zeta$, and $D(y)$ is $\da_1$-fine over $J$.

As for the variation, we can consider the \emph{limit variation for sets} $X$ on the real line by using their indicator $\chi(X;x)$, giving
\[
\mbox{LV}\left(h_y(I,x);Y,\da;a,b;X\right),\;\;\;\;\;\;\;\;
\mbox{LV}\left(h_y(I,x);Y;a,b;X\right)
\]
on substituting $h_y \cdot \chi(X;x)$ for $h_y$ in the previous definitions. We can now have definitions and theorems analogous to Theorems \ref{3.1.12}, \ref{3.1.13}, \ref{3.1.14}, \ref{3.1.15} for the variations.

\begin{example} \label{ex 4.1.1}
Let 
\[
s_n(x) = \left\{
\begin{array}{ll}
1 & (0\leq x<1),\vt
(-1)^p p(p+1) & (2-\frac 1p \leq x<2 - \frac 1{p+1}\;\;\;\;\;(p=1,2, \ldots ,n), \vt
0 & (2 - \frac 1{n+1} \leq x \leq 2).
\end{array}
\right.
\]
Show that $s_n(x)$ is integrable over $[0,2]$ and give the integral. 

If $s_n(x) \rightarrow f(x)$ as $n \rightarrow \infty$, so that
\[
f(x) = \left\{
\begin{array}{ll}
1 & (0\leq x<1),\vt
(-1)^p p(p+1) & (2-\frac 1p \leq x<2 - \frac 1{p+1}\;\;\;\;\;(p=1,2, \ldots ), \vt
0 & (x=2),
\end{array}
\right.
\]
If $N(2) \geq 1$, and $N(x) \geq 1$, $p$, in the respective ranges, and if $n \geq N$, then $s_{n(x)}(x) = f(x)$.

Show that sums, over divisions, of $f(x)$ lie between $0$ and $1$, so that (\ref{4.1.6}) is satisfied for a $C$ that includes $[0,1]$. But
\[
\int_0^{2-\frac 1n} f(x)\,dx = 1 + \sum_{p=1}^n (-1)^p
\]
and does not tend to a limit as $n \rightarrow \infty$.

Thus the Cauchy limit property fails, $f$ is not integrable in $[0,2]$, and (\ref{4.1.6}) cannot be satisfied by a $C$ with arbitrarily small diameter.

\end{example}

\subsection{Special Limit Theorems for Sequences of Functions}\label{s4.2}
Following Lebesgue (1902), W.H.~Young (1910) gave a list of known results on limits of sequences of functions under the integral sign, and so for 80 years these results, including monotone and majorized (dominated) convergence, have been the mainstay of Lebesgue theory.

In this section we prove such sufficient conditions, with others, beginning with monotone convergence over an interval $[a,b]$ ($a<b$) with $b-a$ finite. We use repeatedly the fact that a monotone increasing sequence of real numbers either tends to a limit or tends to infinity.

\begin{theorem}
\label{4.2.1}

For each positive integer $n$ let $s_n(x)$ be gauge integrable over $[a,b]$, and for each such $n$ let 
\[
s_n(x)\leq s_{n+1}(x) \leq g(x) \;\;\;\;\;\;(a\leq x\leq b)
\]
where $g$ is some finite real-valued function, so that
\[
f(x) \equiv \lim_{n \rightarrow \infty} s_n(x)
\]
exists (finite) on $[a,b]$. If $S_n$ is the integral of $s_n(x)$ over $[a,b]$ and if $S_n$ is a bounded sequence with supremum $F$, then $f$ is gauge integrable over $[a,b]$ to the value $F$. 
\end{theorem}
\proof
The integral $S_n(I)$ of $s_n(x)$ exists over every interval $I \subseteq [a,b]$ (Theorem \ref{3.1.8}). As 
\[
\left(S_n\right) \equiv \left(S_n([a,b])\right)
\]
is given to be a bounded sequence with supremum $F$, and monotone increasing by Theorem \ref{3.1.5}, so that $(S_n)$ converges to $F$, then given $\ve>0$, there is an integer $N$ with
\begin{equation}
\label{eq 4.2.1}
F-\ve < S_N \leq S_n \leq F\;\;\;\;\;\;(n>N).
\end{equation}
Also, by Theorem \ref{3.1.10} there is a gauge $\da_n$ such that for all $\da_n$-fine divisions $D$ of $[a,b]$,
\begin{equation}
\label{eq 4.2.2}
(D) \sum \left| s_n(x)(v-u) - S_n([u,v]) \right| < \frac \ve{2^n} \;\;\;\;\;\;(n=1,2, \ldots ).
\end{equation}
Let $t(x) \geq N$ be the least integer for which
\begin{equation}
\label{eq 4.2.3}
f(x) - \ve <s_{t(x)}(x) \leq f(x),\;\;\;\;\;\;\da(x) \equiv \da_{t(x)}(x) >0,
\end{equation}
and let $D$ be a $\da$-fine division of $[a,b]$, so formed of a \textbf{finite} number of interval-point pairs $(I,x)$. Let $p,q$ be the least and greatest values of $t(x)$ for $(I,x)$ in $D$. Then $q\geq p\geq N$. As the integral is finitely additive, and monotone increasing in $n$ (Theorems \ref{3.1.8} and \ref{3.1.5}) and with (\ref{eq 4.2.1}),
\begin{eqnarray*}
F-\ve \;\;<\;\; S_p &=& (D) \sum S_p(I) \;\;\leq \;\;(D) \sum S_{t(x)}(I) \vt
&\leq & (D) \sum S_q(I) \;\;=\;\;S_q \;\;\leq \;\;F.
\end{eqnarray*}
Grouping together the $(I,x)$ in $D$ with equal $t(x)$ and using (\ref{eq 4.2.2}), (\ref{eq 4.2.3}) and monotonicity,
\begin{eqnarray*}
F-2\ve &<&  (D) \sum s_{t(x)}(x)(v-u) \;\;< \;\;F+\ve \vt
F-2\ve & < & (D) \sum f(x)(v-u) \;\;\leq\;\; (D) \sum s_{t(x)}(x)(v-u) +\ve(b-a) \vt
&<& F+\ve +\ve(b-a).
\end{eqnarray*}
As $\ve \rightarrow 0+$, the first and last values tend to $F$, so $f$ integrates to $F$ on the domain $[a,b]$.  \nproof

This is usually called the \emph{weak monotone convergence theorem}. 

For the \emph{strong monotone convergence theorem} we omit just the boundedness of $(s_n(x))$, proving that there is a set $X$ of
$x$ in which $(s_n(x))$ is unbounded, but for which 
\[
s_n^*(x) \equiv s_n(x) \chi(\setminus X;x)
\]
is integrable to the same value as $s_n(x)$, where $\setminus X$ is the \emph{complement} of $X$ and $\chi(\setminus X;x)$ is its \emph{indicator}, taking the value $1$ on $\setminus X$, and the value $0$ on $X$.
\begin{theorem}
\label{4.2.2}
In the conditions of Theorem \ref{4.2.1}, omitting that $(s_n(x))$ is bounded in $n$, and for $X$ the set where $(s_n(x))$ is unbounded in $n$,
\[
\lim_{n \rightarrow \infty}s_n^*(x)\]  is finite and integrable in $ [a,b]$,
\[
\lim_{n \rightarrow \infty}\int_a^b s_n(x)\,dx =
\lim_{n \rightarrow \infty}\int_a^b s_n^*(x)\,dx .
\]
\end{theorem}
\proof
By the monotonicity, $s_n(x) - s_1(x) \geq 0$. Replacing $s_n(x) - s_1(x)$ by $s_n(x)$ in the proof, we assume that $s_n(x)\geq 0$. Let $N$ be a fixed positive integer, with $u(x)$ the smallest integer with
\begin{equation}
\label{eq 4.2.4}
s_{u(x)}(x) \geq N \;\;\;\;(x\in X),\;\;\;\;\;\;\;\;u(x) =1\;\;\;\;(x \in \setminus X).
\end{equation}
Using (\ref{eq 4.2.2}), given $\ve>0$, let $\da_j$ be a gauge on $[a,b]$ such that
\begin{equation}
\label{eq 4.2.5}
(D)\sum\left| s_j(x) (v-u) - S_j([u,v]) \right| < \frac \ve{2^j}
\end{equation}
for each $\da_j$-fine division $D$ of $[a,b]$, and put 
\[
\da(x) = \da_{u(x)}(x)>0.
\]
A $\da$-fine division $D$ of $[a,b]$ has only a finite number of $(I,x)$, and so only a finite number of $u(x)$, which have a maximum, say $W$. By (\ref{eq 4.2.4}), (\ref{eq 4.2.5}) and monotonicity,
\begin{eqnarray*}
(D)\sum N(v-u) \chi(X;x) &\leq & (D) \sum s_{u(x)}(x) (v-u) 
\vt &\leq &(D) \sum S_{u(x)}([u,v]) + \ve
\vt
&\leq & (D) \sum S_{W}([u,v]) + \ve \vt &=& S_W([a,b]) + \ve \;\;\;\leq \;\;\; F+\ve, \vt
V(m(I);[a,b];X)&\leq & V(m(I); \da ; X)\;\;\;\leq \;\;\;\frac{F+\ve}N, \vt
V(m(I);[a,b];X)&=& 0,
\end{eqnarray*}
being true for all positive integers $N$. By Theorem \ref{3.1.12} (\ref{eq 3.1.22}),
\[
\int_a^b s_n(x)\,dx = \int_a^b s_n(x) \chi (\setminus x;x)\,dx,
\]
and the sequence $\left( s_n(x) \chi (\setminus X;x)\right)$ is bounded in $n$ for each $x$ in $[a,b]$, so that Theorem \ref{4.2.1} gives the result. \nproof

Now we can deal with absolutely convergent series in integration.

\begin{theorem}
\label{4.2.3}
Let $s_n(x)$ and $|s_n(x)|$ be integrable in $[a,b]$ ($a<b$). If
\begin{equation}
\label{eq 4.2.6}
\sum_{n=1}^\infty \int_a^b |s_n(x)|\,dx
\end{equation}
is convergent, then
\begin{equation}
\label{eq 4.2.7}
f(x) \equiv \sum_{n=1}^\infty s_n(x)
\end{equation}
is absolutely convergent almost everywhere and is integrable over $[a,b]$ with
\begin{equation}
\label{eq 4.2.8}
\int_a^b f(x)\,dx =\sum_{n=1}^\infty \int_a^b s_n(x)\,dx.
\end{equation}
\end{theorem}
\proof
Bt Theorem \ref{4.2.2} and (\ref{eq 4.2.6}), (\ref{eq 4.2.7}) is absolutely convergent a.e.~since
\begin{eqnarray*}
\int_a^b \sum_{n=1}^n |s_n(x)|\,dx &=& \sum_{n=1}^n \int_a^b  |s_n(x)|\,dx \;\;\leq \;\;\sum_{n=1}^\infty \int_a^b  |s_n(x)|\,dx;\vt
g(x) & \equiv & \sum_{n=1}^\infty |s_n(x)|
\end{eqnarray*}
is integrable over $[a,b]$, and we have (\ref{eq 4.2.8}) with $s_n(x)$ replaced by $|s_n(x)|$. If the $s_n(x)$ are real-valued, the theorem is true for the $|s_n(x)| + s_n(x) \geq 0$.
Subtracting the results for $|s_n(x)|$ gives (\ref{eq 4.2.8}) for real-valued $s_n(x)$. Complex-valued $s_n(x)$ are split into the real and imaginary parts using Theorem \ref{3.1.6}; (\ref{eq 4.2.8}) is true for the real parts and for the imaginary parts, and so for complex-valued $s_n(x)$, finishing the proof. \nproof

We now come to a test that uses bounded Riemann sums, from which we obtain the Arzela-Lebesgue majorized (dominated) convergence test.

We begin with the integrability of the minimum function of a finite collection of functions. Clearly here we need the $s_n(x)$ to be real-valued, and in part of the proof the expert will recognise a refinement integral.
\begin{theorem}
\label{4.2.4}
Let the real-valued $s_n(x)$ be integrable to $S_n$ over $[a,b]$, for $n=N, N+1, \ldots ,P$ ($P>N$). For a gauge $\da$ on $[a,b]$ and a number $M$ let every $\da$-fine division $D$ of $[a,b]$ and every function $n(I,x)$ of interval-point pairs $(I,x)$, integer-valued in the range $N \leq n\leq P$, give
\begin{equation}
\label{eq 4.2.9}
(D) \sum s_{n(I,x)}(x) (v-u) \geq M ; \;\;\;\;\;\;\;\;\;\mbox{ and }\;\;\;\min_{N \leq n(I,x) \leq P} s_{n(I,x)}(x)
\end{equation}
is then integrable over $[a,b]$.
\end{theorem}
\proof
Given the real numbers $s_N, s_{N+1}, \ldots , s_p$ and $t_N, t_{N+1}, \ldots , t_P$, we have
\begin{eqnarray*}
\min_{N \leq n \leq P} s_n &\leq & s_j \;\;=\;\;(s_j -t_j) +t_j \;\;\leq \;\;|s_j-t_j| + t_j \vt
&\leq & |s_N -t_N|+|s_{N+1} -t_{N+1}| + \cdots + |s_P -t_P| +t_j \vt
&&\;\;\;\;\;\;\;\;\;\;\;\;\;\;\;\;\;\;\;\;\;\;\;\;\;\;\;\;\;\;\;\;\;\;\;\;\;\;\;\;\;\;\;\;\;\;\;\;\;\;\;\;\;\;\;\;\;\;\;\;(\mbox{for }N\leq j \leq P).
\end{eqnarray*}
On the right we choose $t_j$ to be the least of the $t_N, \ldots , t_P$, followed by interchanging $s_j,t_j$:
\begin{equation}
\label{eq 4.2.10}
\left| \min_{N \leq n \leq P} s_n - \min_{N \leq m \leq P} t_m \right| 
\leq |s_N -t_N| + \cdots + |s_P -t_P| .
\end{equation}
In the theorem with $s_n(x)$ integrable to $S_n$ over $[a,b]$ ($n=N, \ldots ,P$), there is by Theorem \ref{3.1.11} (\ref{eq 3.1.16}), for $\ve>0$, a gauge $\da_n$ on $[a,b]$ such that for every $\da_n$-fine division $D_n$ of $[a,b]$, 
\begin{eqnarray}
\label{eq 4.2.11}
(D_n) \sum \left| s_n(x) (v-u) -S_n(u,v)\right| &<& \frac \ve{P-N+1},\nonumber \vt
S_n(u,v) &\equiv& \int_u^v s_n(x)\,dx.
\end{eqnarray}
Hence for $\da(x)$ the least of the finite number of gauges $\da_N(x), \ldots , \da_P(x)$, then $\da(x)>0$ in $[a,b]$, and for $D$ a $\da$-fine division of $[a,b]$, we have (\ref{eq 4.2.11}) for $D$ replacing $D_n$. By (\ref{eq 4.2.10}),
\begin{eqnarray}
\label{eq 4.2.12}
&&(D) \sum \left| \min_{N \leq n \leq P}s_n(x) (v-u) -
\min_{N \leq m \leq P}S_m(u,v)\right| \nonumber \vt
&&\leq \;\;\;
\sum_{n=N}^P
\left|  s_n(v-u) - S_n(u,v) \right| 
\;\;\;< \;\;\;\ve.
\end{eqnarray}
We now have to show that as $\da \rightarrow 0$, 
\[
(D) \sum \min_{N \leq n \leq P} S_n(u,v)
\]
tends to a limit. The expert will recognise that we are proving the existence of the refinement integral of the minimum.
\[
\min_{N \leq n \leq P} s_n + \min_{N \leq m \leq P} t_m \leq s_j + t_j \;\;\;\;\;\;(N \leq j \leq P)
\]
proves on taking the minimum of the right-hand side, that
\begin{equation}
\label{eq 4.2.13}
\min_{N \leq n \leq P} s_n + \min_{N \leq m \leq P} t_m \leq \min_{N \leq n \leq P}(s_n + t_n).
\end{equation}
If $a<u<b$ put $S_n(a,u)$ for $s_n$, $S_n(u,b)$ for $t_n$, to show that if $[a,b]$ splits into two intervals at $u$, the sum of the minima does not rise. Repeating for more points of division, we see that on subdividing (i.e.~\emph{refining}) a division, the sum of minima of the $S_n(u,v)$ falls or stays the same. By (\ref{eq 4.2.9}) the infimum $T$ of the sums of minima over divisions of $[a,b]$ exists and $T \geq M$. Given $\ve>0$ let $D^*$ be a division of $[a,b]$ for which the sum of minima lies in $[T,T+\ve)$. By Theorem \ref{3.1.10} there is a gauge $\da'$ such that every $\da'$-fine division $D'$ of $[a,b]$ refines $D^*$.  With the $\da$ for (\ref{eq 4.2.12}) and the gauge $\da'' \equiv \min (\da, \da')$, every $\da''$-fine division $D''$ of $[a,b]$ refines $D^*$ and so its sum of minima lies in $[T, T+\ve)$. By (\ref{eq 4.2.12}),
\[
T-\ve < (D'') \sum \min_{N\leq n \leq P} s_n(x) (v-u) < T+2\ve.
\]
As $\ve>0$ is arbitrarily small, 
\[
\min_{N \leq n \leq P} s_n(x)
\]
is integrable to $T$ over $[a,b]$.  \nproof

We can now turn to majorized or dominated convergence, beginning with a theorem that in Lebesgue theory leads to Fatou's lemma.

\begin{theorem}
\label{4.2.5}
For $a<b$ and the real-valued $s_n(x)$ integrable over $[a,b]$ to $S_n(a,b) \geq M$, with (\ref{eq 4.2.9}) ($m=1,2, \ldots$) true, if $\liminf_{n \rightarrow \infty} S_n(a,b)$ is not $+\infty$, $\liminf_{n \rightarrow \infty} s_n(x)$ exists almost everywhere in $[a,b]$, is integrable over $[a,b]$, and
\begin{equation}
\label{eq 4.2.14}
\int_a^b \liminf_{n \rightarrow \infty} s_n(x)\,dx \leq \liminf_{n \rightarrow \infty} \int_a^b s_n(x)\,dx \equiv \liminf_{n \rightarrow \infty} S_n(a,b).
\end{equation}
\end{theorem}
\proof
By Theorem \ref{4.2.4} the minimum of $s_n(x)$ for $N \leq n \leq P$, is integrable, and it is monotone decreasing in $P$ to the infimum, which is monotone increasing in $N$. BY (\ref{eq 4.2.9}), Theorem \ref{4.2.2} (strong monotone convergence), and $\liminf_{n \rightarrow \infty} S_n(a,b)$ finite, the infimum exists everywhere and is integrable and the $\liminf$ exists almost everywhere and is integrable since
\begin{eqnarray*}
\int_a^b \liminf_{n \rightarrow \infty} s_n(x) \, dx &=& \lim_{N \rightarrow \infty} \int_a^b \inf_{n\geq N} s_n(x)\,dx \vt
&=& \lim_{N \rightarrow \infty} \lim_{P \rightarrow \infty} \int_a^b \min_{N \leq n \leq P} s_n(x) \, dx \vt
& \leq & \lim_{N \rightarrow \infty} \lim_{P \rightarrow \infty} \min_{N \leq n \leq P} \int_a^b s_n(x)\,dx,
\end{eqnarray*}
given finite. Thus we also prove (\ref{eq 4.2.14}). \nproof

Multiplying by $-1$ gives us another theorem.

\begin{theorem}
\label{4.2.6}
For $a<b$ let the real-valued $s_n(x)$ be integrable over $[a,b]$ to $S_n(a,b)$. Let there be a gauge $\da$ on $[a,b]$ and a number $Q$, such that for every $\da$-fine division $D$ of $[a,b]$ and every integer-valued function $n(I,x) \geq 1$,
\begin{equation}
\label{eq 4.2.15}
(D) \sum s_{n(I,x)} (x) (v-u) \leq Q.
\end{equation}
Then \[\max_{N \leq n \leq P} s_n(x)\;\;\;\;\;\;\;\; (P>N)\]
is integrable over $[a,b]$. If also
\[
\limsup_{N \leq n \leq P} s_n(x) \;\;\;\;\;\;\;\;(P>N)
\]
is integrable over $[a,b]$. If also $\limsup_{n \rightarrow \infty} S_n(a,b)$ is finite then \[\limsup_{n \rightarrow \infty} s_n(x)\] exists almost
everywhere in $[a,b]$ and is integrable over $[a,b]$, and
\begin{equation}
\label{eq 4.2.16}
\int_a^b \limsup_{n \rightarrow \infty} s_n(x)\,dx \geq 
 \limsup_{n \rightarrow \infty} \int_a^b s_n(x)\,dx
 \equiv \limsup_{n \rightarrow \infty} S_n(a,b).
\end{equation}
\end{theorem}

Combining Theorems \ref{4.2.5}, \ref{4.2.6}, we have the majorized convergence theorem.

\begin{theorem}
\label{4.2.7}
For $a<b$ let the real-valued $s_n(x)$ be integrable over $[a,b]$ to $S_n(a,b)$. Let there be a gauge $\da$ on $[a,b]$ and numbers $P\leq Q$, such that for every $\da$-fine division $D$ of $[a,b]$ and every integer-valued function $n(I,x) \geq 1$,
\begin{equation}
\label{eq 4.2.17}
P \leq (D)\sum s_{n(I,x)}(x) (v-u) \leq Q.
\end{equation}
If $\lim_{n \rightarrow \infty}s_n(x)$ exists (finite) almost everywhere, then both sides exist (finite) below and
\begin{equation}
\label{eq 4.2.18}
\int_a^b \lim_{n \rightarrow \infty} s_n(x)\,dx =
 \lim_{n \rightarrow \infty} \int_a^b s_n(x)\,dx.
\end{equation}
\end{theorem}
\proof
In (\ref{eq 4.2.17}) take $n(I,x)=n$ (constant). Then  $P \leq S_n(a,b) \leq Q$ and both of 
\[
\limsup_{n \rightarrow \infty}S_n(a,b),\;\;\;\;\;\;\liminf_{n \rightarrow \infty}S_n(a,b)
\]
are finite. By (\ref{eq 4.2.14}), (\ref{eq 4.2.16}), we have (\ref{eq 4.2.18}) from
\begin{eqnarray*}
\int_a^b \lim_{n \rightarrow \infty} s_n(x)\,dx &\leq &
\liminf_{n \rightarrow \infty}\int_a^b  s_n(x)\,dx  \vt
& \leq &
\limsup_{n \rightarrow \infty}\int_a^b  s_n(x)\,dx \vt
&\leq &
\int_a^b \lim_{n \rightarrow \infty} s_n(x)\,dx,
\end{eqnarray*}
giving the result. \nproof

The more usual statement of the majorized (dominated) convergence theorem involves two functions $g(x), h(x)$ integrable over $[a,b]$, such that
\[
g(x) \leq s_n(x) \leq h(x)
\]
for all (or almost all) $x$. Then for a suitable gauge $\da$ and $\da$-fine divisions of $[a,b]$, (\ref{eq 4.2.17}) is true. Conversely,
\[
g(x) \equiv \inf_{n \geq 1} s_n(x),\;\;\;\;\;\;\;\;
h(x) \equiv \sup_{n \geq 1} s_n(x)
\]
are suitable, given (\ref{eq 4.2.17}), and by Theorem \ref{4.2.4} for the integrability.


\vspace{15pt}
\noindent
\textbf{References\footnote{\texttt{Items marked * are editor's additions to Henstock's tentative list.}} (tentative list):}

\begin{enumerate}
\item
* R.G.~Bartle, \emph{Book Review: A General Theory of Integration, by R.~Henstock, 1991}, 
Bulletin of the American Mathematical Society, Volume 29, Number 1, July 1993, pages 136--139.

\item 
R.G.~Bartle (1994--95), \emph{A convergence theorem for generalized Riemann integrals}, Real Analysis Exchange 20, no.~1, 119--124.

\item 
E.~Borel (1895) Annales sci.\'{e}c.normale (3) 12 (p.~51).

\item 
T.J.I'a Bromwich (1931) (1926) \texttt{(* Probably \emph{An Introduction to the Theory of Infinite Series}, Macmillan, London, 1908, 1926)}.

\item
G.~Cantor (1874) \emph{\"{U}ber eine Eigenschaft des Inbegriffes aller reelen algebraischen Zahlen}, J.~Reine Angew.~Math.~77 (pp.~258--260).

\item
G.~Cantor (1875) \emph{Ein Beitrag zur Mannigfaltigkeitslehre}, J. Reine An\-g\-ew. Math., 84, 242--258.

\item
A.L.~Cauchy (1821), \emph{Cours d'analyse de l'\'{E}cole Royale Polytechnique, Analyse Algebrique}, (Works 3(2), Gauthier-Villars, Paris, 1900).

\item
P.~Cousin (1895) \emph{Sur les fonctions de $n$ variables complexes}, Acta Math.~19, 1--62. Jbuch 26, 456.

\item
J.G.~Darboux (1875) \emph{Memoire sur les fonctions discontinues}, Ann. Sci. Ec. Norm. Sup. 4(2) 57--112.  

\item
R.~Dedekind (1909) \emph{Essay on the theory of numbers} (Translated by W.W.~Berman), Chicago.

\item
A.~Denjoy (1912) \emph{Une extension de l'int\'egrale de M.~Lebesgue}, C.R. Acad. Sci. Paris 154, 859--62.

\item
P.~Dienes (1931) The Taylor Series, Clarendon Press, Oxford.

\item
C.~Goffman (1877) \emph{A bounded derivative which is not Riemann integ\-rable}, American Math.~Monthly 84, no.~3, 205--206, MR54\#13000.

\item
S.~Haber and O.~Shisha (1974) \emph{Improper integrals, simple integrals and numerical quadratures}, Journal Approximation Theory 11, 1--15, MR50\#5309.

\item
A.~Harnack (1884) \emph{Die allgemeinen S\"atze \"uber den Zusammenhang der Funktionen einer reelen Variabeln mit ihren Ableitungen,II}, Math.~Annalen 24, 217--52.

\item
Heine

\item
R.~Henstock (1955) \emph{The efficiency of convergence factors for functions of a continuous real variable}, Journal London Math.~Soc.~30, 271--286 (see pp.~277--78), MR17-359.

\item
R.~Henstock (1961) \emph{Definitions of Riemann type of the variational integrals}, Proceedings London Math.~Soc.~3(11), 402--418, MR24\#A1994.

\item 
R.~Henstock (1963),\emph{Theory of Integration}, Butterworth, London.

\item
R.~Henstock (1968a) \emph{A Riemann-type integral of Lebesgue power}, Canadian Journal of Math.~20, 79--87, MR36\#2754.

\item 
R.~Henstock (1968b), \emph{Linear Analysis}, Butterworth, London, MR 34 \#7725. 

\item
R.~Henstock (1988),  \emph{Lectures on the Theory of Integration}, World Scientific, Singapore, MR91a:28001.

\item 
R.~Henstock (1991), \emph{The General Theory of Integration}, Clarendon, Oxford, MR92k:26011.

\item
R.~Henstock (1993--94), \emph{Measure spaces and division spaces}, Real Analysis Exchange 19, no.~1, 121--28.

\item
 * Henstock Archive, Library, University of Ulster, Coleraine, 2007.

\item
Kolmogorov (1933),  \emph{Foundations of the Theory of Probability}, 1933. 

\item
J.~Kurzweil (1957), \emph{Generalized ordinary differential equations and continuous dependence on a parameter}, Czech.~Math.~Journal 7(82), 418--49 (see 422--28), MR22\#2735.

\item
J.~Kurzweil (1980), \emph{Nichtabsolut Konvergente Integrale}, Leipzig, MR 82m: 26007.

\item
H.~Lebesgue (1902), \emph{Integr\'ale, Longueur, Aire}, Annali di Matematica Pura ed Applicata (3)7, 231--359, Jbuch 33, 307.

\item
J.T.~Lewis and O.~Shisha (1983), \emph{The generalized Riemann, simple, dominated and improper integrals}, Journal Approximation Theory 38, 192--99, MR84h:26014.

\item 
E.J.~McShane (1969), \emph{A Riemann type integral that includes Lebesgue-Stieltjes, Bochner and stochastic integrals}, Memoirs American Math. Soc. 88, MR42\#436.

\item 
E.J.~McShane (1973), \emph{A unified theory of integration}, American Math. Monthly 80, 349--59, MR47\#6981.

\item 
* P.~Muldowney (2012), \emph{A Modern Theory of Random Variation, with Applications in Stochastic Calculus, Financial Mathematics, and Feynman Integration}, Wiley, Hoboken, New Jersey, 2012.

\item  
* P.~Muldowney (2016), \newline\emph{Beyond dominated convergence: newer methods of integration},\newline
\texttt{https://sites.google.com/site/stieltjescomplete/\newline $\mbox{ }\;\;\;\;\;\;\;\;\;\;\;\;\;\;\;\;\;\;\;\;\;\;\;\;\;\;\;\;\;\;\;\;\;\;\;\;\;\;\;\;$
\hfill home/convergence-criteria}

\item 
C.F.~Osgood and O.~Shisha, \emph{The dominated integral}, Journal Approximation Theory 17, 150--65, MR54\#7128.

\item 
G.F.B.~Riemann (1868), \emph{\"Uber die Darstellbarkeit einer Funktion durch eine trigonometrische Reihe}, Abh. Gesell. Wiss. Gottingen 13, Math. Kl.  87--132, MR36\#4952.

\item 
S.~Saks (1927),  \emph{Sur les fonctions d'intervalle}, Fundamenta Math.~10, 211--24 (see 214), Jbuch 53, 233.

\item 
S.~Saks (1937), \emph{Theory of the Integral}, 2nd.~English edition, Zbl.~17, 300.

\item 
L.~Tonelli (1924), \emph{Sulla Nozione di Integrale}, Annali di Mat.~(IV)1, 105--45, Jbuch 50, 178.

\item 
C.~de la Vall\'ee Poussin (1892a), \emph{\'Etude des integr\'ales a limites infinies pour lesquelles la fonction sous le signe est continue}, Ann. Soc. Sci. Brux\-el\-les 16 (2nd.~part), 150--80.

\item 
C.~de la Vall\'ee Poussin (1892b),\emph{Recherches sur la convergence des int\'e\-g\-r\-a\-l\-es definies}, J. de Math. pures et appl. 8(4), Fasc. 4, 421--67 (see 453 et sqq.).

\item 
C.~de la Vall\'ee Poussin (1934), Int\'egrales de Lebesgue, Fonctions d'En\-s\-em\-ble, Classes de Baire, 2nd.~edition, Paris, Zbl.~9, 206.

\item 
V.~Volterra (1881), \emph{Sui principii del calcolo integrale}, Giorn.~Mat.~Batt\-ag\-li\-ni, 19, 333--72.

\item 
W.H.~Young (1910), \emph{On semi-integrals and oscillating successions of functions}, Proc.~London Math.~Soc.~(2)9, 286--324.

\end{enumerate}

\end{document}